\documentclass[12pt,psfig]{article}
\usepackage{graphicx}
\def\l {\label}

\def\be {\begin{equation}}
\def\ee {\end{equation}}
\def\ba {\begin{eqnarray}}
\def\ea {\end{eqnarray}}
\newcommand{\dss}{\displaystyle}
\newcommand{\tss}{\textstyle}

\begin{document}
\baselineskip=15pt
\renewcommand {\thefootnote}{\dag}
\renewcommand {\thefootnote}{\ddag}
\renewcommand {\thefootnote}{ }

\pagestyle{empty}

\begin{center}
                \leftline{}
                \vspace{ 0.0 in}

{\LARGE \bf Four limit cycles from perturbing quadratic \\ [0.5ex] 
integrable systems by quadratic polynomials} \\ [0.3in]

{\large Pei\, Yu\,$^{a,\,b\,*}$}
\footnote{$^*\,$Corresponding author. Fax: (519) 661-3523; 
Email: pyu@uwo.ca}
$\!\!$and \ {\large Maoan\, Han\,$^a$}
\\

\vspace{0.15in}
{\small {\em 
$^a$Department of Mathematics, Shanghai 
Normal University, Shanghai, China \ 200234
}}

\vspace{0.05in}
{\small {\em
$^b$Department of Applied Mathematics, The University of Western 
Ontario \\[-0.5ex]
London, Ontario, Canada \ N6A 5B7
}}


\noindent
\end{center}

\baselineskip=15pt 
\vspace{0.0in}

\noindent 
\rule{6.5in}{0.012in}

\vspace{-0.2in}
\section*{Abstract}
In this paper, we give a positive answer to the open question: 
Can there exist $4$ limit cycles in quadratic near-integrable polynomial 
systems?  
It is shown that when a quadratic integrable system has two centers 
and is perturbed by quadratic polynomials, 
it can generate at least $4$ limit cycles 
with $(3,1)$ distribution. The method of Melnikov function is used.

\vspace{0.1in} 

\noindent 
{\it Keywords}: Hilbert's 16th problem, quadratic near-integrable system, 
limit cycle, 

\hspace{0.50in}reversible system, Hopf bifurcation, Poincar\'{e} bifurcation, 
Melnikov function  

\vspace{0.10in} 
\noindent
{\it MSC}: 34C07; 34C23 

\noindent 
\rule{6.5in}{0.012in}

\section{Introduction}\label{intro}

\setcounter{equation}{0}
\renewcommand{\theequation}{1.\arabic{equation}}

The well-known Hilbert's 16th problem is remained unsolved since 
Hilbert~\cite{hilbert1900} proposed the 23 mathematical problems at the 
Second International Congress of Mathematics in 1990. 
Recently, a modern version of
the second part of the 16th problem was formulated by Smale~\cite{smale1998},
chosen as one of the 18 challenging mathematical problems for the 21st
century. To be more specific, consider the following planar system: 
\begin{eqnarray}
\frac{dx}{dt} = P_n(x,y), \quad \frac{dy}{dt} = Q_n(x,y),
\label{b1}
\end{eqnarray}
where $\, P_n(x,y) \,$ and $\, Q_n(x,y) \,$ represent $\, n^{\rm th}$-degree
polynomials of $\, x\,$ and $\, y$. The second part of Hilbert's 16th 
problem is to find the upper bound $\, H(n) \leq n^q\,$ on the number
of limit cycles that the system can have, where $\, q\,$
is a universal constant, and $H(n)$ is called Hilbert number.
In early 90's of the last century, 
Ilyashenko~\cite{Ilyashenko1991}
and \'{E}calle~\cite{Ecalle1992} proved the finiteness theorem pioneered by 
Dulac, for given planar polynomial vector fields.  
In general the finiteness problem has not been solved even for 
quadratic systems.  
A recent survey article~\cite{li2003} (and more references therein)
has comprehensively discussed this problem and reported the recent progress.

If the problem is restricted to
the neighborhood of isolated fixed points, then the
question is reduced to studying degenerate Hopf bifurcations, which
give rise to fine focus points.
In the past six decades, many researchers have considered the local
problem and obtained many results
(e.g., see~[6--12]). 
In the last 20 years, much progress
on finite cyclicity near a fine focus point or a homoclinic loop has been
achieved. Roughly speaking, the so-called finite cyclicity means that
at most a finite number of limit cycles can exist in some neighborhood
of focus points or homoclinic loop under small perturbations on the
system's parameters.

\pagestyle{myheadings}
\markright{{\footnotesize {\it P. Yu} \& {\it M. Han
\hspace{1.2in} {\it 
4 limit cycles in near-quadratic systems}}}}

In this paper, we particularly consider bifurcation of limit cycles 
in quadratic systems. Early results can be found in a survey article
by Ye~\cite{Ye1982}. Some recent progress has been reported in 
a number of papers (e.g., see~\cite{Roussarie1998,RoussarieSchlomiuk2002}). 
For general quadratic system (\ref{b1}) ($n=2$), 
in 1952, Bautin~\cite{bautin} proved that there exist $3$ small limit cycles 
around a fine focus point or a center. After $30$ years, 
until the end of 1970's, 
concrete examples were given to show that general quadratic systems can have 
$4$ limit cycles~\cite{ChenWang,Shi80}, around two foci
with $(3,1)$ configuration.  
Since then, many researchers have paid attention to integrable 
quadratic systems, and a number of results have been obtained. 
A question was naturally raised: 
Can near-integrable quadratic systems have $4$ limit cycles?
A quadratic system is called near-integrable 
if it is a perturbation of a quadratic integrable system by quadratic 
polynomials. On one hand, it is reasonable to believe that the answer should 
be positive since general quadratic systems have at least $4$ limit cycles; 
while on the other hand, near-integrable quadratic systems have limitations on 
their system parameters and thus it is more difficulty to find $4$ limit cycles 
in such systems. In fact, this is still an open problem after 
another $30$ years since the finding of $4$ limit cycles in general 
quadratic systems.   

The study of bifurcation of limit cycles for near-integrable systems is 
related to the so called weak Hilbert's 16th problem~\cite{Arnold1977}, 
which is transformed to finding the maximal number of 
isolated zeros of the Abelian integral or Melnikov function: 
\be 
M(h,\delta) = \dss\oint_{H(x,y)=h} Q_n \, dx - P_n \, dy, 
\l{b2} 
\ee 
where $\, H(x,y), \, P_n \,$ and $\, Q_n \,$ are all real polynomials 
of $\,x \,$ and $\, y\,$ with $\, {\rm deg} H = n+1$, 
and $\, \max\{ {\rm deg} P_n, \, {\rm deg} Q_n \} \le n$. 
The weak Hilbert's 16th problem is a very important problem, closely 
related to the maximal number of limit cycles of the following 
near-Hamiltonian system~\cite{Han2006}: 
\be  
\dss\frac{dx}{dt} = \dss\frac{\partial H(x,y)}{\partial y}
+ \varepsilon \, p_n(x,\,y), \quad 
\dss\frac{dy}{dt} =- \,  \dss\frac{\partial H(x,y)}{\partial x} 
+ \varepsilon \, q_n(x,\,y), 
\l{b3} 
\ee
where $\, H(x,y)$, $ p_n (x,y) \,$ and $\, q_n(x,y) \,$ are 
polynomials of $\, x \,$ and $\, y$, and $\, 0 < \varepsilon \ll 1 \,$ 
is a small perturbation.

General quadratic systems with one center have been classified by 
\.{Z}ol\c{a}dek~\cite{Zoladek1994} using a complex analysis on the 
condition of the center, as four systems: 
$Q_3^{LV}$ -- the Lotka-Volterra system; 
$Q_3^H$ -- Hamiltonian system; $Q_3^R$ -- reversible system; and 
$Q_4$ -- codimension-4 system.  
In 1994, 
Horozov and Iliev~\cite{HI1994} proved that in 
quadratic perturbation of 
generic quadratic Hamiltonian vector fields with one center and 
three saddle points there can appear at most two limit cycles, and this 
bound is exact. Later, Gavrilov~\cite{Gavrilov2001}
extended Horozov and Iliev's method to give a fairly complete analysis on 
quadratic Hamiltonian systems with quadratic perturbations. 
Quadratic Hamiltonian systems, with  at most four singularities, 
can be classified as three 
cases~\cite{Gavrilov2001}: (i) one center and three saddle points; 
(ii) one center and one saddle point; and (iii) two centers and 
two saddle points.  
In~\cite{Gavrilov2001}, Gavrilov showed that like case (i), cases (ii) 
and (iii) can also have at most two limit cycles. Therefore, 
generic quadratic Hamiltonian systems with quadratic perturbations
can have maximal two limit cycles, and this case has been completely solved.

For the $Q_3^R$ reversible system, there have been many results published. 
For example, Dumortier {\it el al.}~\cite{DumortierLiZhang1997} 
studied a case of $Q_3^R$ system 
with two centers and two unbounded heteroclinic loops, 
and presented a complete analysis of quadratic $3$-parameter unfolding. 
It was proved that $3$ is the maximal number of limit cycles surrounding a 
single focus, and only the $(1,1)$-configuration can occur in case of 
simultaneous nests of limit cycles. That is, $3$ is the maximal number 
of limit cycles for the system they studied~\cite{DumortierLiZhang1997}. 
Later, Peng~\cite{Peng2002} considered a similar case with a homoclinic 
loop and showed that $2$ is the maximal number of limit cycles which 
can bifurcate from the system.   
Around the same time, Yu and Li~\cite{YuLi2002} investigated a similar 
case as Peng considered but with a varied parameter in a 
certain interval, and obtained the same conclusion as Peng's. Later, 
Iliev~{\it et al.}~\cite{IlievLiYu2005} re-investigated the same case 
but for the varied parameter in a different interval (which yields two centers) 
and got the same conclusion as that of~\cite{DumortierLiZhang1997}, i.e., 
$3$ is the maximal number of limit cycles which can be obtained from 
this case.   
Recently, Li and Llibre~\cite{LL2006} considered 
a different case of $Q_3^R$ system which can exhibit
the configurations of limit cycles: $(0,0)$, $(1,0)$, $(1,1)$ and 
$(1,2)$. Again, no $4$ limit cycles were found.  
In order to explain why the above authors did not find $4$ limit cycles 
from the $Q_3^R$ reversible system, 
consider the $Q_3^R$ system with quadratic perturbations, 
which can be described by~\cite{DumortierLiZhang1997} 
\be 
\begin{array}{ll} 
\dot{x} = -\, y + a\, x^2 + b\, y^2 
+ \varepsilon \, (\mu_1 \, x + \mu_2 \, x\,y), \\[1.0ex] 
\dot{y} = x\, (1 + c\, y) + \varepsilon\, \mu_3 \, x^2,  
\end{array} 
\l{b4} 
\ee 
where  $\, a, \, b, \, c \,$ are real parameters, $\, \mu_i, \ i=1,\,2,\,3\,$ 
are real perturbation parameters, and $\, 0 < \varepsilon \ll 1$. 
When $\, \varepsilon = 0$, system (\ref{b4})$_{\varepsilon=0}$ is a reversible 
integrable system. It has been noted that in all 
the cases considered 
in~\cite{DumortierLiZhang1997,Peng2002,YuLi2002,IlievLiYu2005}, 
the parameters 
$\,a \,$ and $\, c \,$ were chosen as $\, a = -\,3, \ c = -\,2$, but with 
$\, b=1 \,$ in~\cite{DumortierLiZhang1997}; 
$\, b=-\,1\,$ in~\cite{Peng2002},  
$\, b \in (-\,\infty, -1) \cup (-1,0)\,$ in~\cite{YuLi2002}, 
and $\, b \in (0,\,2)\,$ in~\cite{IlievLiYu2005}.  
In these papers, complete analysis on the perturbation parameters 
was carried out with the 
aid of Poincar\'{e} transformation and the Picard-Fuchs equation,
but it needed to fix all (or most of) the parameters $\,a, \, b \,$ and $\,c$.
This way it may miss opportunity 
to find more limit cycles, such as possible existence of $4$ limit cycles.   
As a matter of fact, for the cases considered 
in~\cite{YuLi2002,IlievLiYu2005}, a simple scaling on the parameter $b$ 
($b \ne 0$) can be used to eliminate $b$. 
So, suppose the non-perturbed system (\ref{b4})$_{\varepsilon=0}$ has two 
free parameters and let us consider the 2-dimensional parameter plane. 
Then, all the cases studied in the above mentioned articles are special cases, 
represented by just a point or a line segment in the 2-dimensional 
parameter plane (see more details in Section 2). 
It has been noted that a different method was used in~\cite{LL2006} with 
Melnikov function up to second order, but no more 
limit cycles were found.

It should be mentioned that Zhang~\cite{PGZhang2002} has proved 
that the possible cycle distributions in 
general quadratic systems with two foci must be $(0,1)$-distribution 
or $(1,i)$-distribution, $\, i=0,\,1,\,2,\,3, \cdots$. 
So far, no results have been obtained for $\, i \ge 4$.
This result also rules out the possibility of 
$\,(2,2) $-distribution. 
It is conjectured that at most $3$ limit cycles 
can exist around one focus point.    
The problem of bifurcation of $3$ limit cycles near 
an isolated homoclinic loop is still open. 

In this paper, we turn to a different angle to consider bifurcation of 
limit cycles in quadratic near-integrable systems with two centers. 
We shall leave more free parameters in the integrable systems, so that 
we will have more chances to find more limit cycles. 
The basic idea is as follows: we first consider bifurcation of 
multiple limit cycles from Hopf singularity, which does not need to 
fix any parameters, and use expansion of Melnikov function near 
centers to get such limit cycles as many as possible. 
This leads to determination of a maximal number of parameters. Then,
for the remaining undetermined parameters, we compute the global 
Melnikov function to look for possible large limit cycles.  
Indeed, although, due to the complex 
integrating factor in the analysis, we are not able to give a complete 
analysis for classifying the perturbation unfolding, we do get a positive 
answer to the open question of existence of $4$ limit cycles in 
quadratic near-integrable systems. 
In particular, we will show that perturbing a 
reversible, integrable quadratic system 
with two centers can have at least $4$ limit cycles, with $(3,1)$ 
distribution, bifurcating from 
the two centers under quadratic perturbations. 

The rest of paper is organized as follows. 
In Section 2, we give a different classification in real domain 
for quadratic systems with one center, and compare it with that 
given by \.{Z}ol\c{a}dek~\cite{Zoladek1994}. 
Also, we use our classification to present a simple summary on 
some of the existing results for the reversible near-integrable system. 
Section 3 is devoted to the analysis on bifurcation of small limit cycles 
from Hopf singularity. In Section 4, we show how to find large 
limit cycles bifurcating from closed orbits to obtain a total of 
$4$ limit cycles. Finally, conclusion is drawn in Section 5.

\section{Classification of generic quadratic systems with at least one center}

\setcounter{equation}{0}
\renewcommand{\theequation}{2.\arabic{equation}}

In this section, we give a different classification in real domain 
for quadratic systems with a center, 
which is consistent with the Hamiltonian systems considered 
in~\cite{HI1994,Gavrilov2001}. 
We start from the following general quadratic system: 
\vspace{0.0in} 
\be  
\begin{array}{ll} 
\dss\frac{d z_1}{dt} = c_{100} + c_{110}\, z_1 + c_{101}\, z_2 
+ c_{120}\, z_1^2 + c_{111}\, z_1\, z_2 + c_{102}\, z_2^2, \\ [1.5ex]  
\dss\frac{d z_2}{dt} = c_{200} + c_{210}\, z_1 + c_{201}\, z_2 
+ c_{220}\, z_1^2 + c_{211}\, z_1\, z_2 + c_{202}\, z_2^2, 
\end{array}
\l{b5} 
\ee 
where $\, c_{ijk}$'s are real constant parameters. 
It is easy to show that this system has at most four singularities, 
or more precisely, it can have $0$, $2$ or $4$ singularities in real domain.
In order for system (\ref{b5}) to have limit cycles, the system must 
have some singularity.
In this paper, we assume that system (\ref{b5}) has at least two singularities. 
Without loss of 
generality, we may assume that one singular point is located at 
the origin $(0,0)$, which implies $\, c_{100} = c_{200} = 0$, 
and the other at $(p,q)$ ($p^2 + q^2 \ne 0$). 
Further assume the origin is a linear center. 
Then introducing a series of linear transformations, parameter 
rescaling and time rescaling to system (\ref{b5}) yields 
the following general quadratic system: 
\be  
\begin{array}{ll} 
\dss\frac{dx}{dt} = y + a_1\, x\, y + a_2\, y^2, \\ [1.5ex]  
\dss\frac{dy}{dt} = -\, x 
+ x^2 + a_3\, x\, y + a_4\, y^2, 
\end{array}
\l{b6} 
\ee 
which has a linear center at the origin $(0,0)$ and another singularity 
at $ (1,0)$.

In order to have the origin of system (\ref{b6}) being a center, we 
may calculate the focus values of system (\ref{b6}) and find 
four cases under which $(0,0)$ is a center, listed in the 
following theorem (here we use \.{Z}ol\c{a}dek's 
notation in our classification). 

\vspace{0.10in} 
\noindent  
{\bf Theorem 1.1} {\it The origin of (\ref{b6}) is a center if and only 
if one of the following conditions is satisfied: 

\vspace{0.10in} 
\noindent 
\underline{$Q_3^R$\,--\,Reversible system}: $\, a_3 = a_2 = 0$, under which 
system (\ref{b6}) becomes 
\be 
\begin{array}{ll} 
\dss\frac{d x}{dt} = y + a_1\, x\, y, \\[1.0ex]  
\dss\frac{d y}{dt} =  -\, x + x^2 + a_4 \, y^2 , 
\end{array}
\l{b7} 
\ee 
with 
$$  
(1,0) \ \ {being \ a} \ \ 
\left\{ 
\begin{array}{ll} 
{center} \quad & {if} \ \ a_1 < -\,1, \\ [1.0ex] 
{saddle \ point} \quad & {if} \ \ a_1 > -\,1. 
\end{array} 
\right.
$$

\vspace{0.10in} 
\noindent 
\underline{$Q_3^H$\,--\,Hamiltonian system}: $\, a_3 = a_1 + 2\, a_4 = 0$, 
under which system (\ref{b6}) is reduced to 
\be 
\begin{array}{ll} 
\dss\frac{d x}{dt} =  y + a_1 \, x\, y  + a_2 \, y^2 , \\[1.0ex] 
\dss\frac{d y}{dt} = -\,x + x^2 - \dss\frac{1}{2}\, a_1\,  y^2 , 
\end{array} 
\l{b8} 
\ee 
with 
$$  
(1,0) \ \ {being \ a} \ \ 
\left\{ 
\begin{array}{ll} 
{center} \quad & {if} \ \ a_1 < -\,1, \\ [1.0ex] 
{saddle \ point} \quad & {if} \ \ a_1 > -\,1. 
\end{array} 
\right.
$$

\vspace{0.10in} 
\noindent 
\underline{$Q_3^{LV}$\,--\,Lokta-Volterra system}: $\, a_2 = 
1+a_4 = 0$, under which system (\ref{b6}) becomes 
\be 
\begin{array}{ll} 
\dss\frac{d x}{dt} =  y + a_1 \, x\, y , \\[1.0ex] 
\dss\frac{d y}{dt} = -\,x + x^2 + a_3 \, x\, y - y^2 , 
\end{array} 
\l{b9} 
\ee 
with
$$  
(1,0) \ \ {being \ a} \ \ 
\left\{ 
\begin{array}{ll} 
{focus} \quad & {if} \ \ a_1 < -\,( 1 + \frac{1}{4}\, a_3^2 \,), 
\\ [1.0ex] 
{node} \quad & {if} \ \ -\,( 1 + \frac{1}{4}\, a_3^2 \,) < a_1 < -\,1, 
\\ [1.0ex] 
{saddle \ point} \quad & {if} \ \ a_1 > -\,1. 
\end{array} 
\right.
$$

\vspace{0.10in} 
\noindent 
\underline{$Q_4$\,--\,Codimension-4 system}:
\be 
a_3 - 5\, a_2 = a_1 - (5 + 3\, a_4) = a_4 + 2 (1+a_2^2) = 0, 
\l{b10} 
\ee 
under which system (\ref{b6}) can be rewritten as 
\be 
\begin{array}{ll} 
\dss\frac{d x}{dt} =  y - (1 + 6\, a_2^2) \, x\, y + a_2 \, y^2 , \\[1.0ex] 
\dss\frac{d y}{dt} = -\,x + x^2 + 5\, a_2 \, x\, y - 2\, (1 + a_2^2)\, y^2 , 
\end{array} 
\l{b11} 
\ee 
with $\, (1,0)\,$ being a node for $\, a_2 \ne 0$. 
} 

\vspace{0.10in} 
\noindent 
{\bf Remark 1.2.}\ 
There is one more case found from the above process, 
defined by the following conditions:   
\be 
a_3 - 5\, a_2 = a_1 - (5 + 3\, a_4) = 
3\,(a_4+2)\,(a_4+1)^2 - (5\,a_4+6)\, a_2^2 = 0. 
\l{b12} 
\ee 
We will show later in this section, when we compare our above real 
classification with the complex classification given by 
\.{Z}ol\c{a}dek~\cite{Zoladek1994}, that the case defined by 
(\ref{b12}) actually belongs to the $\, Q_3^R$-reversible system.

\vspace{0.10in} 
\noindent 
{\bf Proof.} Necessity is easy to be verified by computing the focus values 
of system (\ref{b6}) associated with the origin. 
Some focus values will not equal zero if the condition is not satisfied. 

For sufficiency, we find an integrating factor for each case when the 
condition holds. For the $Q_3^H$\,-\,Hamiltonian system (\ref{b8}), 
we know that the 
integrating factor is $\,1$, and the Hamiltonian is given by 
\be 
H(x,y) = \dss\frac{1}{2}\, (x^2 + y^2) - \dss\frac{1}{3}\, x^3 
+ \dss\frac{1}{2}\, a_1 \, x\, y^2 + \dss\frac{1}{3}\, a_2 \, y^3, 
\l{b13} 
\ee 
which is exactly the same as that given in~\cite{HI1994,Gavrilov2001}. 

For the $Q_3^R$\,-\,reversible system (\ref{b7}), the integrating factor is 
\be
\gamma = |1+a_1 x|^{- \frac{a_1+2 a_4}{a_1}},  
\l{b14} 
\ee 
and the first integral of the system is given by 
\be 
F(x,y) = \dss\frac{1}{2}\, 
{\rm sign}(1+a_1 x) \, |1+a_1 x|^{- \frac{2 a_4}{a_1}}  
\left[ y^2 + \dss\frac{(1+a_1-a_4)\, (1 + 2\, a_4\, x)}{a_4\, (a_1 - a_4)\, 
(a_1 - 2\, a_4)} - \dss\frac{x^2}{a_1 - a_4} \right].  
\l{b15} 
\ee

For the $Q_3^{LV}$-\,Lokta-Volterra system (\ref{b9}), 
we find the integrating factor to be 
\be
\gamma = | g(x,y) |^{-1}, \quad 
{\rm where} \ \ 
g(x,y) = (1+a_1 x)\, \Big[(x-1)^2 + a_3 \, (x-1) \, y - (1+a_1)\, y^2 \Big], 
\l{b16} 
\ee 
and the first integral of the system is 
\be 
F(x,y) = \!  
\left\{ \!\!  
\begin{array}{ll} 
-\, \dss\frac{{\rm sign} (g(x,y))} {2\, a_1 (1+a_1)}  
{\LARGE \Big\{} 2\, \ln |1\!+\!a_1 x|  \!+\! 
a_1 \ln \left| (1\!+\!a_1) y^2 \!-\! a_3  y (x\!-\!1) \!-\! (x\!-\!1)^2 
\right|   \\[1.0ex]  
\hspace{1.2in} 
+\, \frac{ 2\,a_1 \,a_3\,(x-1)}{
\sqrt{[a_3^2 + 4\, (1+a_1)]\,(x-1)^2}}\, 
\tanh^{-1} \! \Big[ 
\frac{ a_3 \,(x-1) - 2 ( 1+a_1) \,y}
{\sqrt{[a_3^2 + 4\, (1+a_1)]\,(x-1)^2}}
\Big] {\LARGE \Big\}}, \\ [2.5ex] 
\qquad \qquad \qquad  {\rm when} \quad  a_3^2 + 4\,(1+a_1) > 0,   
\\[2.0ex]
-\, \dss\frac{{\rm sign} (1\!+\!a_1 x)} {2\, a_1 (1+a_1)}  
{\LARGE \Big\{} 2\, \ln |1\!+\!a_1 x|  \!+\! 
a_1 \ln \left[ (1\!+\!a_1) y^2 \!-\! a_3  y (x\!-\!1) \!-\! (x\!-\!1)^2 
\right]   \\[1.0ex]  
\hspace{1.2in} 
-\, \frac{ 2\,a_1 \,a_3\,(x-1)}{
\sqrt{[-\,a_3^2 - 4\, (1+a_1)]\,(x-1)^2}}\, 
\tan^{-1} \! \Big[ 
\frac{ a_3 \,(x-1) - 2 ( 1+a_1) \,y}
{\sqrt{[-\,a_3^2 - 4\, (1+a_1)]\,(x-1)^2}}
\Big] {\LARGE \Big\}}, \\ [2.5ex] 
\qquad \qquad \qquad  {\rm when} \quad  a_3^2 + 4\,(1+a_1) < 0.   
\end{array} 
\right.      
\l{b17} 
\ee 

Finally, for the $Q_4$\,-\,codimension-4 system (\ref{b11}), we have 
\be
\gamma = | g(x,y) |^{-5/2}, \quad {\rm where}  \ \ 
g(x,y) =  1 - 2\,( 1 + 2\, a_2^2)\,x - 2\, a_2 \, y
+ (1+4\, a_2^2)\, (x + a_2\, y)^2, 
\l{b18} 
\ee 
and the first integral of the system is equal to   
\be 
F(x,y) = \dss\frac{{\rm sign} (g(x,y))}{12\, a_2^6} \, 
| g(x,y)|^{-3/2} f(x,y), 
\l{b19} 
\ee
 
\vspace{-0.05in}  
\noindent 
where 
\be 
\begin{array}{rl}   
f(x,y) =\!\!& -\,(1+a_2^2) 
+ 3\,( x + a_2 \, y + 2\, a_2^2\, x) 
\, \Big[ 1+a_2^2 - (1+3\,a_2^2)\, ( x+ a_2\, y) \Big] \nonumber \\ [1.5ex] 
\!\!& +\, (1+3\,a_2^2)\, (1+4\,a_2^2)\, (x + a_2\,y)^3. 
\end{array}  
\l{b20} 
\ee 
The proof is complete. 
\put(10,0.5){\framebox(6,7.5)}

Note that among the four classifications of the integrable system (\ref{b6}), 
the first three classified systems (\ref{b7}), (\ref{b8}) and (\ref{b9}) 
have two free parameters, while the last system (\ref{b11}) 
has only one free parameter.

\vspace{0.10in} 
\noindent 
{\bf Remark 1.3.} We now show that our classification in Theorem 1.1 is 
equivalent to that given by \.{Z}ol\c{a}dek~\cite{Zoladek1994}.
The general quadratic system considered in~\cite{Zoladek1994} is given 
in the complex form: 
\be 
\dss\frac{d z}{dt} = (i + \lambda)\, z + A \, z^2 + B \, z \, \bar{z} 
+ C\, \bar{z}^2,  
\l{b21} 
\ee 
where $\, z = x + i\, y$, and $\, A, \ B\,$ and $C\,$ are complex 
coefficients. It has been shown in~\cite{Zoladek1994} that 
the point $\, z = 0\,$ is a center if and only if one of the 
following conditions is fulfilled: 
\be 
\begin{array}{ll} 
Q_3^{LV}: & \lambda = B = 0, \\[1.5ex]  
Q_3^H: & \lambda = 2\,A + \bar{B} = 0, \\[1.5ex]  
Q_3^R: &  \lambda = {\rm Im}(AB) = {\rm Im} (\bar{B}^3 C) 
= {\rm Im} (A^3 C) = 0,  \\[1.5ex] 
Q_4: &  \lambda = A - 2\,\bar{B} = |C| - |B| = 0. 
\end{array} 
\l{b22} 
\ee

In the following, 
we first use real differential equation to give a brief proof 
(different from \.{Z}ol\c{a}dek's~\cite{Zoladek1994}), 
and then show that our classification is equivalent to \.{Z}ol\c{a}dek's 
when system (\ref{b21}) is assumed to have a non-zero singularity. 
To prove this, let
$$ 
A = A_1 + i\, A_2, \quad B = B_1 + i\, B_2, \quad C = C_1 + i\, C_2, 
\quad (i^2 = -\,1), 
$$ 
and then rewrite the complex equation (\ref{b21}) in the real form:   
\be  
\begin{array}{ll} 
\dss\frac{d x}{dt} = \lambda\,x + y 
+ (A_1 + B_1 + C_1)\, x^2 
+ 2\, (A_2 - C_2)\, x\, y 
- (A_1 - B_1 + C_1)\, y^2 , \\[1.0ex] 
\dss\frac{d y}{dt} = -\,x + \lambda\,y 
- (A_2 + B_2 + C_2)\, x^2 
+ 2\, (A_1 - C_1)\, x\, y 
+ (A_2 - B_2 + C_2)\, y^2 , 
\end{array} 
\l{b23} 
\ee 
where $\, y \rightarrow -\, y \,$ has been used. 
Letting $\, \lambda = 0 \,$ yields the focus value $\, v_0 = 0$. 
Then, it is easy to find the first focus value (or the first Lyapunov constant) 
as 
\be 
v_1 =  -\,A_1\, B_2 - B_1\, A_2 = -\, {\rm Im}(AB). 
\l{b24} 
\ee 
Letting $\, v_1 = 0 \,$ results in $\, {\rm Im}(AB) = 0$, 
which gives 
\be 
B_2 = -\, \dss\frac{B_1\, A_2}{A_1}, \quad 
{\rm under \ the \ assumption \ of} \ \ A_1 \ne 0. 
\l{b25}
\ee 
(The degenerate case $A_1 = 0$ can be similarly analyzed and the 
details are omitted here.) 
Then, we apply our Maple program (e.g., see~\cite{Yu1998}) to 
system (\ref{b23}), with the conditions $\lambda =0$ and (\ref{b25}),
to obtain 
$$
v_2 = \dss\frac{-\,f\, (A_1 - 2\, B_1)}{3 A_1^3} , \ \ \ 
v_3 = \dss\frac{-\,f\, f_3}{216 A_1^5}, \ \ \ 
v_4 = \dss\frac{-\,f\, f_4}{9720 A_1^7}, \ \ \ 
v_5 = \dss\frac{-\,f\, f_5}{466560 A_1^9}, \ \ \cdots 
$$
where 
$$ 
f = B_1\, (2\,A_1 +B_1)\, ( C_2\, A_1^3 + 3\, C_1\,A_2\, A_1^2 - 
3\, A_2^2 \, C_2 \, A_1 - C_1 \, A_2^3), 
$$ 
and $\, f_3, \, f_4$, etc. are polynomials of 
$\, A_1,\,A_2,\,C_1,\,C_2\,$ and $\,B_1$. 
Letting $\, f =0$, i.e., 
$$ 
B_1 = 0 \quad {\rm or} \quad 
2\,A_1 +B_1 = 0 \quad {\rm or} \quad 
C_2 \, A_1^3 + 3 \,C_1 \,A_2 \,A_1^2 - 
3 \,A_2^2\, C_2 \,A_1 - C_1 \,A_2^3 = {\rm Im}(A^3 C) = 0  
$$ 
yields $\, v_2 = v_3 = \cdots =0$.
 
Indeed, $\, B_1 = 0\,$ implies $\, B_2 = 0$ due to the condition (\ref{b25}), 
and so $\, B = 0$. 
Thus, we obtain $\, \lambda =  B = 0$, corresponding to the $Q_3^{LV}\,$ case. 

For the condition $\,2\,A_1 +B_1 = 0$, it follows from (\ref{b25}) that 
$\, 2\,A_2 - B_2 = 0$, i.e., $\, 2\, A + \bar{B} =0$, which 
plus the condition $\,\lambda =0 \,$ gives the $\, Q_3^H \,$ case.  

The third condition $\, {\rm Im}(A^3 C) = 0 $, with $\, \lambda = 0 \,$ 
and $\, {\rm Im}(AB) = 0$, corresponds to the $Q_3^R \,$ case. 
Further, it is easy to show that under the condition $\, {\rm Im}(AB) = 0$, 
$\, {\rm Im}(A^3 C) = 0 \,$ and $\, {\rm Im}(\bar{B}^3 C) = 0 \,$ 
are equivalent. Thus, the conditions $\, \lambda = {\rm Im}(AB) = 
{\rm Im}(\bar{B}^3 C) = 0 \,$ also applicable for this case. 
So for this case, either $\,  {\rm Im}(A^3 C) = 0\,$ 
or $\, {\rm Im}(\bar{B}^3 C) = 0 \,$ 
is needed, but not both of them. 
In the following, we show one more case to join this case, leading to 
both the two conditions being needed.

Note that there is one more condition $\, A_1 = 2\, B_1 \,$ which renders 
$\, v_2 = 0$. Letting $\, A_1 = 2\, B_1$, and so 
$\, A_2 = -\, 2\, B_2 \,$ (see (\ref{b25})), implying that 
$\, A - 2\, \bar{B} = 0$. 
Under the condition $\, A = 2\, \bar{B} $, 
$\, v_1 = v_2 = 0$, and the other focus values become
\begin{eqnarray*} 
v_3 &\!\!\!=\!\!\!& 
\dss\frac{25}{8} \, 
(C_1^2+C_2^2-B_1^2-B_2^2)
( C_2 B_1^3 -3 C_1 B_1^2 B_2 -3 C_2 B_1 B_2^2 + C_1 B_2^3 ), \\[0.5ex] 
v_4  &\!\!\!=\!\!\!& 
\dss\frac{v_3}{45}\,
\Big[45 B_1^2
+585 B_2^2
+60 (B_1 C_1+B_2 C_2)
-196 (C_1^2+C_2^2)
\Big], \\[0.5ex]  
v_5  &\!\!\!=\!\!\!& 
\dss\frac{v_3}{6480}\,
\Big[ 
648 (7 B_1^4 \!+\! 124 B_1^2 B_2^2 \!+\! 1557 B_2^4) 
\!-\! 3 (961 B_1^2 C_1^2 \!-\! 7680 B_1 B_2 C_1 C_2 \!+\! 202345 B_2^2 C_2^2)  
\ \  \\ 
&\!\!\!  \!\!\!&
\qquad \quad +\, 576 B_1 C_1 (106 B_1^2 + 307 B_2^2) 
+ 288 B_2 C_2 (371 B_1^2 + 773 B_2^2) 
\\ 
&\!\!\!  \!\!\!&
\qquad \quad -\, 3 (4801 B_1^2 C_2^2 + 206185 B_2^2 C_1^2) 
-80688 (C1^2+C2^2) (B_1 C_1 + B_2 C_2) 
\\ 
&\!\!\!  \!\!\!&
\qquad \quad +\, 86144 (C1^2+C2^2)^2 
\Big], \\[-1.0ex] 
&\!\!\! \vdots \!\!\!& 
\end{eqnarray*} 
Hence, under the conditions $\, \lambda = A - 2\, \bar{B} = 0$, 
there are two possibilities such that 
$\, v_3 = v_4 = \cdots = 0$. The first possibility is 
$$ 
C_1^2 + C_2^2 - B_1^2 - B_2^2, \quad {\rm i.e.,} \quad 
|C| - |B| = 0, 
$$ 
which is one of the conditions given for the $Q_4$ case 
(see (\ref{b22})). 

The second possibility is given by the condition: 
\be  
C_2 \, B_1^3 -3 \,C_1 \,B_1^2 \,B_2 -3\,C_2\, B_1\, B_2^2 + C_1\, B_2^3 = 
{\rm Im} (\bar{B}^3 C) =  \tss\frac{1}{8}\, {\rm Im} (A^3 C) = 0,
\l{b26} 
\ee 
due to $\, A = 2 \ \bar{B}$. 
Since these conditions can be included in the conditions 
$\, \lambda = {\rm Im}(AB) = {\rm Im} (\bar{B}^3 C) = {\rm Im} (A^3 C) = 0$, 
this possibility belongs to the $\, Q_3^R \,$ case.

The remaining task is to show that the conditions classified in 
(\ref{b22}) are sufficient. This can be done by finding an integrating factor
for each case. For brevity, we only list these integrating factors 
below (while the lengthy expressions of the first integrals are omitted):  
\be 
\gamma = \left\{ 
\begin{array}{ll} 
\Big| 1+4 \,(A_2\, x-A_1 \,y) + 4 \,(A_1 C_2+A_2 C_1-2 A_1 A_2) \,x\, y 
\\ [1.0ex]
\quad +\,  \left[ (A_1+C_1)\,(A_1-3 C_1)+(A_2+C_2)\,(5 A_2-3 C_2)\right] x^2 
\\ [1.0ex] 
\quad +\, \left[(A_2+C_2)\,(A_2-3 C_2)+(A_1+C_1)\,(5 A_1-3 C_1)\right] y^2
\\ [1.0ex] 
\quad +\, 
2\, (A_1^2+A_2^2-C_1^2-C_2^2) \left[ (A_2+C_2)\, x^3 - (A_1+C_1)\, y^3 
\right. \\ [0.5ex] 
\hspace{1.0in} -\left. (A_1-3 C_1)\, x^2 y  + (A_2-3 C_2) \,x y^2  \right] 
\Big|^{-1}, & {\rm for} \ \ Q_3^{LV}, \\[2.0ex]  
1, & {\rm for} \ \ Q_3^{H}, \\[1.5ex] 
\left| 1 - 2\, (A_1 - C_1)\, y \right|^{-\, \frac{2\,A_1 + B_1}{A_1-C_1}}, 
& {\rm for} \ \ Q_3^{R}, \\[2.0ex]
\Big| 1 - 4 ( B_2\, x + B_1\, y) 
+ 2 (B_1^2+B_2^2)\, (x^2+y^2)  \\ [0.0ex] 
\quad +\, 2 (B_1 C_1 + B_2 C_2)\,(x^2-y^2) 
+ 4 \,(B_1 C_2 - B_2 C_1)\, x\, y \Big|^{-5/2}, & {\rm for} \ \
Q_4. 
\end{array}
\right.  
\l{b27}
\ee 
For the integrating factors of degenerate cases (e.g., 
$\, A_1 - C_1 = 0$), one can easily find them.

Next, compare the classification listed in (\ref{b22}) with ours 
given in Theorem 1.1. 
First, consider the $\, Q_3^{LV}\,$ case. Letting 
$\, \lambda = B_1 = B_2 = 0\,$ in (\ref{b23}) yields 
\be 
\begin{array}{ll} 
\dss\frac{d x}{dt} =  y + (A_1 + C_1)\, x^2 + 2\, (A_2 - C_2)\, x\, y 
- (A_1 + C_1)\, y^2 , \\[1.0ex] 
\dss\frac{d y}{dt} = -\,x - (A_2 + C_2)\, x^2 + 2\, (A_1 - C_1)\, x\, y 
+ (A_2 + C_2)\, y^2 . 
\end{array} 
\l{b28} 
\ee 
Then, let 
\be
k = \tan(\theta), \quad {\rm and \ so} \ \ 
\sin(\theta) = \dss\frac{k}{\sqrt{1+k^2}}, \quad 
\cos(\theta) = \dss\frac{1}{\sqrt{1+k^2}}, 
\l{b29} 
\ee 
where $\, k\,$ is solved from the following cubic polynomial: 
\be 
P_1(k) = (A_2+C_2)\, k^3 + (A_1 - 3\, C_1)\, k^2 + (A_2 - 3\,C_2)\, k 
+ A_1 + C_1 = 0.  
\l{b30} 
\ee  
This cubic polynomial at least has one real solution for $k$, which gives the 
slope of the line on which a second fixed point is located. 
$\, k = 0 \,$ if $\, A_1 + C_1 =0$, otherwise, $k \ne 0$.  
Let $\overline{k}$ be a real root of $P_1(k)$, i.e., $P_1(\overline{k})=0$. 

Further, introducing the linear transformation (rotation): 
\be 
x = \cos(\theta)\, u - \sin(\theta)\, v, \quad 
y = \sin(\theta)\, u + \cos(\theta)\, v, 
\l{b31} 
\ee 
into (\ref{b28} yields 
\be 
\begin{array}{ll} 
\dss\frac{d x}{dt} = y + m_{120}\, x^2 + m_{111}\, x\, y + m_{102}\, y^2, \\[1.0ex]
\dss\frac{d y}{dt} =-\,x + m_{220}\, x^2 + m_{211}\, x\, y + m_{202}\, y^2, 
\end{array} 
\l{b32} 
\ee 
where 
$$ 
\begin{array}{ll} 
m_{120} = -\, m_{102} = (1+ \overline{k}^2)^{-3/2}\, 
P_1(\overline{k})=0, 
\\[1.5ex] 
m_{220} = -\, m_{202} = (1+ \overline{k}^2)^{-3/2}\, \Big[ 
(A_1+C_1)\, \overline{k}^3 - (A_2-3\,C_2)\, \overline{k}^2 
+ (A_1 - 3\, C_1)\, \overline{k}  
-A_2 - C_2 \Big], \\[1.5ex] 
m_{111} = -\,2\,(1+ \overline{k}^2)^{-3/2}\, \Big[ 
(A_1-C_1)\, \overline{k}^3 - (A_2+3\,C_2)\, \overline{k}^2 
+ (A_1 + 3\, C_1)\, \overline{k}  
-A_2 + C_2 \Big], \\[1.5ex]  
m_{211} = 2\,(1+ \overline{k}^2)^{-3/2}\, \Big[ 
(A_2-C_2)\, \overline{k}^3 + (A_1+3\,C_1)\, \overline{k}^2 
+ (A_2 + 3\, C_2)\, \overline{k}  
+A_1 - C_1 \Big]. 
\end{array} 
$$ 
Suppose $\, m_{220} \ne 0$. Then, introducing 
$\, \overline{x} = m_{220}\, x, \ \overline{y} = m_{220}\, y \,$ into 
(\ref{b32}) results in 
\be 
\begin{array}{ll} 
\dss\frac{d \overline{x}}{dt} = \overline{y} + \dss\frac{m_{111}}{m_{220}}\, 
\overline{x}\, \overline{y}, \\[2.0ex]
\dss\frac{d \overline{y}}{dt} =-\, \overline{x} + \overline{x}^2 
+ \dss\frac{m_{211}}{m_{220}}\, x\, y - \overline{y}^2, 
\end{array} 
\l{b33} 
\ee  
which is identical to (\ref{b9}) as long as letting 
$\, a_1 = \frac{m_{111}}{m_{220}} \,$ 
and $\, a_3 = \frac{m_{211}}{m_{220}}$. 
This shows that the four parameters 
$\, A_1, \, A_2, \, C_1\, $ and $\, C_2 \,$ are not independent. Thus, 
alternatively, we may simply take $ \overline{k} = 0 $ 
(which renders the second singularity of (\ref{b28}) on the $x$-axis), 
yielding $\, C_1 = -\,A_1 $. Thus, (\ref{b28}) becomes 
$$ 
\begin{array}{ll} 
\dss\frac{d x}{dt} =  y + 2\, (A_2 - C_2)\, x\, y , \\[1.0ex] 
\dss\frac{d y}{dt} = -\,x - (A_2 + C_2)\, x^2 + 4\, A_1 \, x\, y 
+ (A_2 + C_2)\, y^2 . 
\end{array} 
$$
Suppose $\, A_2 + C_2 \ne 0$. Introducing 
$\, \overline{x} = -\,(A_2+C_2)\, x , \ 
\overline{y} = -\,(A_2+C_2)\, y \, $ into the above equations 
we obtain 
\be 
\begin{array}{ll} 
\dss\frac{d \overline{x}}{dt} =  \overline{y} 
- \tss\frac{2\,(A_2 - C_2)}{A_2 + C_2} \, \overline{x}\, \overline{y} , 
\\[1.5ex] 
\dss\frac{d \overline{y}}{dt} = -\, \overline{x} + \overline{x}^2 
- \tss\frac{4\, A_1}{A_2 + C_2} \, \overline{x}\, \overline{y} 
- \overline{y}^2 , 
\end{array} 
\l{b34} 
\ee 
which is identical to (\ref{b9}) if letting 
$\,a_1 = \frac{-\,2\,(A_2 - C_2)}{A_2 + C_2}  \,$ 
and $\, a_3 = \frac{-\,4\, A_1}{A_2 + C_2}$.  
In the following, we will use this simple approach  for other cases.

For the $\, Q_3^H $ case, substituting $\, \lambda = 0$, 
$\,B_1 = -\,2\,A_1 \,$ and 
$\, B_2 = 2\, A_2 \,$ into system (\ref{b23}) results in 
$$
\begin{array}{ll} 
\dss\frac{d x}{dt} =  y - (A_1 - C_1)\, x^2 + 2\, (A_2 - C_2)\, x\, y 
- (3\,A_1 + C_1)\, y^2 , \\[1.0ex] 
\dss\frac{d y}{dt} = -\,x - (3\, A_2 + C_2)\, x^2 + 2\, (A_1 - C_1)\, x\, y 
- (A_2 - C_2)\, y^2 . 
\end{array} 
$$
Further, taking $\, C_1 = A_1\, $ 
in the above equations gives another singularity on the $x$-axis, 
and introducing $\, \overline{x} = -\,(3\,A_2+C_2)\, x , \ 
\overline{y} = -\,(3\,A_2+C_2)\, y \, $ into the resulting equations 
yields 
\be 
\begin{array}{ll} 
\dss\frac{d \overline{x}}{dt} =  \overline{y} 
- \tss\frac{2\, (A_2 - C_2)}{3\, A_2 + C_2} \, \overline{x}\, \overline{y} 
+ \tss\frac{4\, A_1}{3\, A_2 + C_2} \overline{y}^2, \\[1.5ex] 
\dss\frac{d \overline{y}}{dt} = -\,\overline{x} + \overline{x}^2 
+ \tss\frac{A_2 - C_2}{3\, A_2 + C_2} \, \overline{y}^2 , 
\end{array} 
\l{b35} 
\ee 
which is identical to (\ref{b8}) if we set 
$\, a_1 = \frac{-\,2\, (A_2 - C_2)}{3\, A_2 + C_2}\,$ 
and $\, a_2 = \frac{4\, A_1}{3\, A_2 + C_2}$. 

For the $Q_3^R$ reversible case, it follows from~\cite{Zoladek1994} that 
all the coefficients $A, \, B \,$ and $ C\, $ are real, and thus we obtain 
the following real form from the complex system (\ref{b21})
\be 
\begin{array}{ll} 
\dss\frac{d x}{dt} =  -\, y + a \, x^2  + b \, y^2 , \\[1.0ex] 
\dss\frac{d y}{dt} = x + c\, x\, y, 
\end{array} 
\l{b36} 
\ee  
where 
$$
a = A_1 + B_1 + C_1, \quad  b = B_1 - A_1 - C_1 , \quad 
c = 2\, A_1 - 2\, C_1. 
$$ 
Suppose $\,b \ne 0$. Then, introducing 
$\, \overline{x} = b\, y, \ \overline{y} = b\, x \,$ 
into (\ref{b36}) results in 
\be 
\begin{array}{ll} 
\dss\frac{d \overline{x}}{dt} = \overline{y} + \dss\frac{c}{b}\, 
\overline{x}\, \overline{y}, \\[1.0ex]  
\dss\frac{d \overline{y}}{dt} =  -\, \overline{x} + \overline{x}^2 
+ \dss\frac{a}{b} \, \overline{y}^2, 
\end{array}
\l{b37} 
\ee
which is identical to (\ref{b7}) if  
$$ 
a_1 = \dss\frac{c}{b} = \dss\frac{2\,(A_1 - C_1)}{B_1 - A_1 -C_1} \quad  
{\rm and} \quad 
a_4 = \dss\frac{a}{b} = \dss\frac{A_1 + B_1 + C_1}{B_1 - A_1 -C_1}. 
$$

For the last $\, Q_4$ case, 
under the condition $\, \lambda = A - 
2\, \bar{B} = 0 $, by setting 
$\, C_1 = -\,3 \, B_1 \, $  (which renders a non-zero singularity 
on the $x$-axis) in (\ref{b23}) we obtain 
$$ 
\begin{array}{ll}
\dss\frac{d x}{dt} = y - 2\,(2\,B_2 +C_2) \, x\, y + 2\, B_1 \, y^2,
\\[1.0ex]
\dss\frac{d y}{dt} =  -\, x + (B_2-C_2)\, x^2 + 10 \, B_1 \, x\, y
- ( 3\, B_2 - C_2) \, y^2,
\end{array}
$$
Suppose $\, B_2 - C_2 \ne 0$.
Then, introducing $\, \overline{x} = (B_2 - C_2)\, x, \ 
 \overline{y} = (B_2 - C_2)\, y \,$ into the above
equations yields
\be
\begin{array}{ll}
\dss\frac{d \overline{x}}{dt} = \overline{y} 
- \tss\frac{2\,(2\,B_2 +C_2)}{B_2-C_2} \, \overline{x}\, \overline{y}
+ \tss\frac{2\, B_1}{B_2-C_2} \, \overline{y}^2,
\\[1.5ex]
\dss\frac{d \overline{y}}{dt} =  -\, \overline{x} + \overline{x}^2 
+ \tss\frac{10 \, B_1}{B_2-C_2} \, \overline{x}\, \overline{y}
- \tss\frac{ 3\, B_2 - C_2}{B_2-C_2} \, \overline{y}^2,
\end{array}
\l{b38}
\ee
Comparing the coefficients of the above system (\ref{b38}) 
with our system (\ref{b6}) results in 
\be 
a_1 = -\,\tss\frac{2\,(2\,B_2 +C_2)}{B_2-C_2}, \quad
a_2 = \tss\frac{2\, B_1}{B_2-C_2}, \quad
a_3 =  \tss\frac{10 \, B_1}{B_2-C_2}, \quad
a_4 = -\,  \tss\frac{ 3\, B_2 - C_2}{B_2-C_2},
\l{b39}
\ee
which in turn implies that
$\,  a_3 - 5\, a_2= a_1 - (5+3\,a_4) = 0$,
and
$$ 
a_4 + 2\, (1+a_2^2) = \tss\frac{8\,B_1^2 + C_2^2 - B_2^2}{(B_2-C_2)^2} 
= \tss\frac{C_1^2 + C_2^2 - B_1^2 - B_2^2}{(B_2-C_2)^2} 
= 0, \quad {\rm for} \ \ |C| - |B|=0. 
$$
The above conditions are the exact conditions given in (\ref{b10}) for
the $\, Q_4\, $ case. 

Finally, we turn to the conditions 
given in (\ref{b12}). It follows from (\ref{b39}) that 
\be 
3\,(a_4+2)\,(a_4+1)^2 - (5\,a_4+6)\, a_2^2 = 
-\,\tss\frac{ 4}{(B_2-C_2)^3} \, 
( 3\,B_2^3+3\,B_2^2\,C_2-C_1^2\,B_2-B_1^2\,C_2 ). 
\l{b40} 
\ee 
On the other hand, under the condition $\, C_1 = -\,3\,B_1 $,
the condition (\ref{b26}) for the second possibility becomes 
$$ 
\begin{array}{rl} 
C_2 B_1^3 -3 C_1 B_1^2 B_2 -3 C_2 B_1 B_2^2 + C_1 B_2^3 
= \!\!&C_2 B_1^3 + C_1^2 B_1 B_2 -3 C_2 B_1 B_2^2 - 3 B_1 B_2^3 
\\ [1.5ex]  
= \!\!& -\, B_1 \, ( 3\,B_2^3+3\,B_2^2\,C_2-C_1^2\,B_2-B_1^2\,C_2 ) 
= 0 ,
\end{array}
$$ 
which implies, by Eq.~(\ref{b40}), 
$\, 3\,(a_4+2)\,(a_4+1)^2 - (5\,a_4+6)\, a_2^2 = 0 \,$ for $\, B_1 \ne 0$. 
Hence, according to \.{Z}ol\c{a}dek's classification (see (\ref{b22})), 
this case should be included in the $Q_3^R \,$ case. 
However, one can not prove this by directly using the conditions in (\ref{b12}) 
as well as that for the $Q_3^R $ case (see Theorem 1.1). One must 
trace back to the original system coefficients. 

In~\cite{Zoladek1994}, \.{Z}ol\c{a}dek used Bautin's system to verify his 
classification. Bautin's system is described by~\cite{bautin} 
\be  
\begin{array}{ll} 
\dss\frac{d x}{dt} = \lambda_1 \,x - y 
+ \lambda_3 \, x^2 
+ (2\, \lambda_2 + \lambda_5)\, x\, y 
+ \lambda_6 \, y^2 , \\[1.0ex] 
\dss\frac{d y}{dt} = x + \lambda_1 \,y 
+ \lambda_2 \, x^2 
+ (2\,\lambda_3 + \lambda_4)\, x\, y 
- \lambda_2 \, y^2 . 
\end{array} 
\l{b41}
\ee  
It is seen from (\ref{b23}) and (\ref{b41}) that Bautin's system 
has only $6$ parameters, while  \.{Z}ol\c{a}dek's system 
has $7$ (in real domain) parameters. This indicates that 
\.{Z}ol\c{a}dek's system has one redundant parameter. 
In fact, putting Bautin's system in \.{Z}ol\c{a}dek's complex form gives 
the following expressions: 
$$ 
\begin{array}{ll}  
\lambda = \lambda_1, \quad 
A = \tss\frac{1}{4}\,(\lambda_3 + \lambda_4 - \lambda_6 - i\, \lambda_5), 
\quad B = -\,\tss\frac{1}{2}\, ( \lambda_3 - \lambda_6), \\[1.0ex] 
C = \tss\frac{1}{4}\, \Big[ - ( 3\, \lambda_3 + \lambda_4 + \lambda_6) 
+ i\, (4\,\lambda_2 + \lambda_5) \Big]. 
\end{array}
$$ 
Then, applying the formulas given in (\ref{b23}) will immediately 
generate the centers conditions obtained by Bautin~\cite{bautin}. 
The above expressions clearly show that $\, B_2 = 0$. As a matter of factor, 
the integral factor for the system, corresponding to the second possibility, 
i.e., when $\lambda = A - 2 \, \bar{B} = {\rm Im} (\bar{B}^3 C ) = 0$, is 
given by 
$$ 
\left| 1 + 2 \left[ \tss\frac{ C_1 (B_1^2+B_2^2)} {B_1 (B_1^2-3 B_2^2)}
-2 \right]  (B_2\, x+ B_1\, y) \right|^{\frac{5 B_1 (B_1^2-3 B_2^2) }
{C_1 (B_1^2+B_2^2)-2 B_1 (B_1^2-3 B_2^2)}}. 
$$
For $\, B_2 = 0$, the above expression is reduced to 
$$
\Big| 1 - 2\, ( 2\, B_1 - C_1 )\, y  \Big|^{ \frac{5 \, B_1}{C_1 - 2\, 
B_2}}
=  \Big| 1 - 2\, (A_1 - C_1)\, y \Big|^{-\, \frac{2\,A_1 + B_1}{A_1-C_1}} 
\quad ({\rm due \ to} \ \, A_1 = 2\, B_1),   
$$  
which is the integrating factor for the $Q_3^R \,$ system, as shown in 
(\ref{b27}). 

\begin{figure}[!h]
\vspace{-1.45in}
\hspace{-0.50in} 
\resizebox{1.2 \textwidth}{!}{\includegraphics{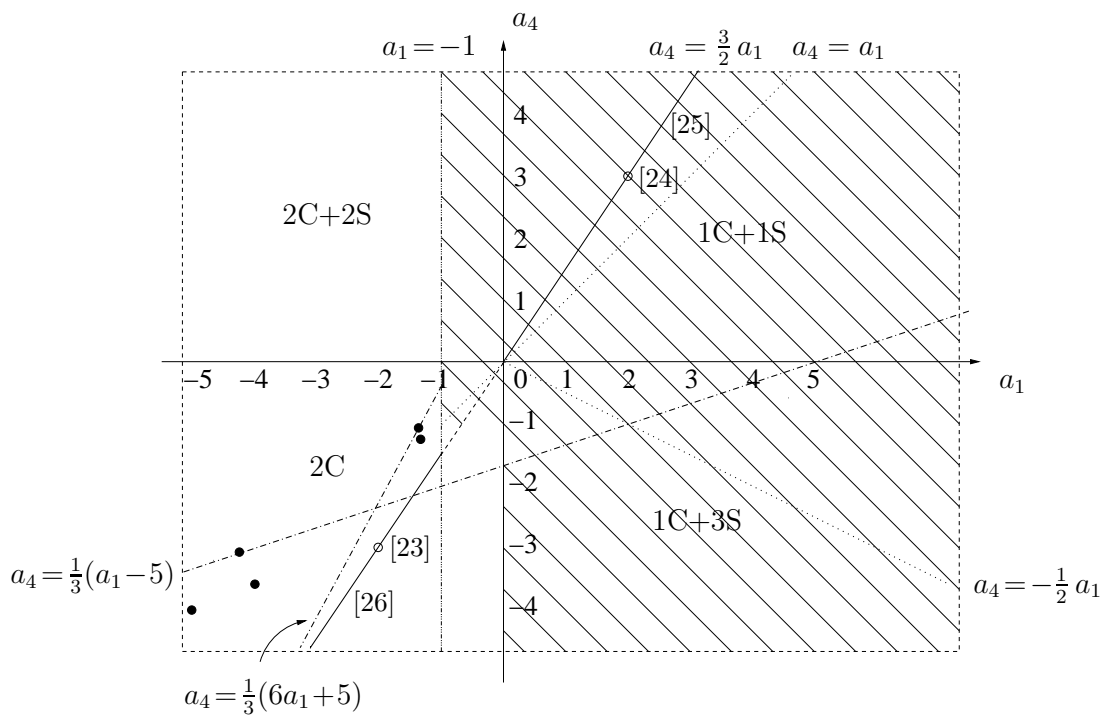}}

\vspace{-5.90in} 
\caption{Case studies for the $Q_3^R$ reversible system.} 
\label{fig1} 
\vspace{0.10in}
\end{figure}

\vspace{0.10in} 
Now we return to system (\ref{b6}). 
Among the four classifications, 
the Hamiltonian system ($Q_3^H$) has been completely studied
in~\cite{HI1994,Gavrilov2001}: the system can have maximal two limit cycles. 
In this paper, we will concentrate on the $Q_3^R$\,-\,reversible case. 
Special cases for the reversible system have been investigated by 
a number of authors 
(e.g., see~\cite{DumortierLiZhang1997,Peng2002,YuLi2002,IlievLiYu2005,LL2006}). 
It is easy to see that system (\ref{b7}) 
is invariant under the mapping $(t,\,y) \to (-t,\, -y)$, 
where $\, a_1 \,$ and $\, a_4 \,$ can be considered as perturbation 
parameters. 
The singular point $(1,0)$ of (\ref{b7}) is a center when $a_1 < -\,1$; 
but a saddle point when $\, a_1 > -\,1 $.
$a_1 = -\,1 \,$ gives a degenerate singular point at $(1,0)$. 
Further, it is easy to verify that  
when $\, (a_1 +1)\,a_4 >0$, there are no more singularity; 
while when $\, (a_1 + 1)\, a_4 < 0 $, there exist additional two 
saddle points, given by
$$ 
(x^*\!,y^*) = \Big(- \frac{1}{a_1}, \ 
\pm \tss\frac{ \sqrt{ -\, a_4 \, ( a_1 + 1)}}{a_1 \, a_4} \Big).  
$$ 
$a_4 = 0\,$ is a critical value, yielding the two additional saddle 
points at infinity: $ (x^*\!, y^*) = ( -\, \frac{1}{a_1},\, \pm \infty)$. 
In summary, the distribution of singularity of the reversible 
system (\ref{b7}) has the following possibility
(see Fig.~\ref{fig1}, where 1C+1S stands for one center and one 
saddle point, similar meaning applies to 2C, 2C+2S and 1C+3S): 
\be  
\begin{array}{ll} 
{\rm Two \ centers \ when} \ \ a_1 < -\,1 \ \, {\rm and} \, \ a_4 < 0;  \\[1.0ex] 
{\rm Two \ centers \ and \ two \ saddle \ points \ when} \ \  
a_1 < -\,1 \, \ {\rm and} \, \ a_4 > 0; \\[1.0ex]  
{\rm One \ center \ and \ one \ saddle \ point \ when} \ \ 
a_1 > -\,1 \, \ {\rm and} \, \ a_4 > 0;& 
\\[1.0ex] 
{\rm One \ center \ and \ three \ saddle \ points \ when} \ \ 
a_1 > -\,1 \, \ {\rm and} \, \ a_4 < 0.& 
\end{array}
\l{b43} 
\ee 

In this paper, we pay particular attention to 
$\, a_1 < -\,1, \ a_4 < 0$, for which system (\ref{b7}) has only two 
singularities at $(0,0)$ and $(1,0)$, both of them are centers. 

By adding quadratic perturbations to system (\ref{b7}) we obtain  
the following perturbed quadratic system: 
\be  
\begin{array}{rl} 
\dss\frac{dx}{dt} = & y \, (1 + a_1\, x) 
+ \varepsilon \, P(x,y)  \\ [1.0ex] 
= & y\, (1 + a_1\, x) 
+ \varepsilon \, (a_{10}\, x + a_{01}\, y 
+a_{20}\, x^2 + a_{11}\, x\, y + a_{02}\, y^2), \\ [0.5ex]  
\dss\frac{dy}{dt}= & -\, x + x^2 + a_4\, y^2 + \varepsilon \, Q(x,y) \\ [1.0ex] 
= & -\, x + x^2 + a_4\, y^2 
+ \varepsilon \, (b_{10}\, x + b_{01}\, y 
+ b_{20}\, x^2 + b_{11}\, x\, y + b_{02}\, y^2), 
\end{array}
\l{b44} 
\ee 
where $\, 0 < \varepsilon \ll 1 $, $a_{ij}$'s and $b_{ij}$'s are perturbation 
parameters.

\vspace{0.10in} 
\noindent 
{\bf Remark 1.4.} The special system considered in~\cite{DumortierLiZhang1997}
is the system (\ref{b4}) with 
$$ 
a = -\,3, \quad c = -\,2, \quad b = 1. 
$$ 
This is equivalent to our system when $\, a_1 = -\,2 \,$ and
$\, a_4 = -\,3$ for which the system has only two centers at $(0,0)$ and 
$(1,0)$.  
Consider the $a_1$-$a_4$ parameter plane, as shown in Fig.~\ref{fig1}. 
It can be seen that the case considered in~\cite{DumortierLiZhang1997} is just a 
point, $(a_1, a_4) = (-\,2,\, -\,3) $, in the parameter plane, 
marked by a blank circle in the third quadrant 
on the line $\, a_4 = \frac{3}{2}\, a_1$ (see Fig.~\ref{fig1}).  

The special system studied in~\cite{Peng2002}
is the system (\ref{b4}) with
$$ 
a = -\,3, \quad c = -\,2, \quad b = -\,1. 
$$
This is equivalent to our system when $\, a_1 = 2 \,$ and
$\, a_4 = 3$, for which the system has one center at $(0,0)$ 
and one saddle point at $(1,0)$.
Thus, this  case considered in~\cite{Peng2002} is 
again a point, $\, (a_1, a_4) = (2,\, 3) $, 
in the $a_1$-$a_4$ parameter plane, 
marked by another blank circle in the first quadrant 
on the line $\, a_4 = \frac{3}{2}\, a_1$ (see Fig.~\ref{fig1}).

The cases considered in~\cite{YuLi2002,IlievLiYu2005} 
correspond to the system (\ref{b4}) with $\, a=-\,3, \ c=-\,2$, 
and $\, b \in (-\,\infty, -1) \cup (-1,0)\,$ in~\cite{YuLi2002}, 
and $\, b \in (0, 2) $ in~\cite{IlievLiYu2005}. 

When $\, \varepsilon = 0$ in system (\ref{b4}), one can use the 
following transformation: 
$$ 
x = \dss\frac{\tilde{y}}{b}, \quad y = \dss\frac{\tilde{x}}{b}, 
$$ 
to transform system (\ref{b4})$_{\varepsilon=0}$ to 
\be 
\begin{array}{ll}  
\dss\frac{d \tilde{x}}{dt} = 
\tilde{y} \, \Big( 1 + \dss\frac{c}{b} \, \tilde{x} \Big), \\[1.5ex] 
\dss\frac{d \tilde{y}}{dt} = -\, \tilde{x} + \tilde{x}^2 
+ \dss\frac{a}{b}\, \tilde{y}^2,  
\end{array} 
\l{b45} 
\ee 
which is our system (\ref{b7}) with 
\be 
a_1 = \dss\frac{c}{b}, \quad a_4 = \dss\frac{a}{b}. 
\l{b46} 
\ee
Equation (\ref{b46}) yields  
\be  
a_4 = \dss\frac{a}{c} \ a_1 \qquad (b \ne 0), 
\l{b47} 
\ee
which represents a line in the $a_1$-$a_4$ parameter 
plane, passing through the origin with the slope $\, \frac{a}{c}$.  
In particular, the parameter values: 
$\, a = -\,3, \ c = -\,2, \ 
b \in (- \infty, -1) \cup (-1, 0) \cup (0, 2)$, yielding 
$a_1 = -\,\frac{2}{b} \,$ and $a_4 = -\, \frac{3}{b}$, 
correspond to a part of the line, described by 
\be 
a_4 = \frac{3}{2}\, a_1  \qquad  
\forall \, a_1 \in (- \infty, -1) \cup (0, \infty), 
\l{b48} 
\ee
as shown in Fig.~\ref{fig1}, 
where the dotted line for $ a_1 \in [-1,0]$ is excluded from the 
studies~\cite{YuLi2002,IlievLiYu2005}. 
  
It should be noted that when $\, a=-\,3, \ c=-\,2$, 
the point $(0, \frac{1}{b})$ is 
a saddle point if and only if 
$$ 
1 + \dss\frac{c}{b} = 1 - \dss\frac{2}{b} > 0 \quad 
\Longrightarrow \quad  b \in (-\infty, 0) \cup (2, +\infty). 
$$
Thus, the case considered in~\cite{YuLi2002} has 
one center and one saddle point; while the case studied in~\cite{IlievLiYu2005}
has two centers. 
But even these two studies together do not cover the whole line 
$\, a_4 = \frac{3}{2}\, a_1$ (the missing part is denoted by a 
dotted line segment in Fig.~\ref{fig1}).

Another alternative form for a special case of our system (\ref{b7}) 
considered in~\cite{Han1997} is described by 
\be  
\begin{array}{ll} 
\dss\frac{d \overline{x}}{dt} = \overline{y}\, \Big[ 1 + 2\, ( 1-e)\,
\Big( \overline{x} + \dss\frac{1}{d} \Big) \Big], \\ [1.5ex]  
\dss\frac{d \overline{y}}{dt} = \overline{x} + d\, \overline{x}^2 + e\, 
\overline{y}^2, 
\end{array}
\l{b49} 
\ee
where $\, e \,$ and $\, d \ (\ne 0) \,$ are parameters.  
This system has a saddle point at the origin and a center at 
$\, (\overline{x},\, \overline{y}) = ( -\frac{1}{d},\, 0)$. 
Based on the two parameters, seven cases are classified~\cite{Han1997}. 
We can apply the following transformation: 
$$ 
\overline{x} = \dss\frac{1}{d} \, ( x -1), \quad 
\overline{y} = \dss\frac{1}{d}\, y, 
$$ 
to system (\ref{b13}), yielding 
\be  
\begin{array}{ll} 
\dss\frac{d x}{dt} = y\, \Big[ 1 + \dss\frac{2\, ( 1-e)}{d} \, x
 \Big], \\ [2.0ex]  
\dss\frac{d y}{dt} = -\, x + x^2 + \dss\frac{e}{d} \, y^2 , 
\end{array}
\l{b50} 
\ee 
which has a center at the origin and a saddle point at $(1,0)$. 
Then, setting 
\be  
a_1 = \dss\frac{2\,(1-e)}{d}, \quad a_4 = \dss\frac{e}{d}, 
\l{b51} 
\ee 
in system (\ref{b50}) leads to our system (\ref{b7}). 
Equation (\ref{b51}) denotes a line, given by 
\be 
a_4 = \dss\frac{e}{2\,(1-e)} \ a_1, 
\l{b52} 
\ee 
in the $a_1$-$a_4$ parameter plane, passing through the origin with the slope 
$\, \frac{e}{2\,(1-e)}$. However, it is easy to see that 
using our system (\ref{b7}) in analysis is simpler than using 
system (\ref{b49}). In fact, all the seven cases classified 
in~\cite{Han1997} together denote a region in Fig.~\ref{fig1}, 
see the shaded area in this figure. 
This area covers most of the region, defined by $a_1 > -\,1$.  
But the study given in~\cite{Han1997} for the seven cases is restricted 
to local analysis on the bifurcation of limit cycles near a homoclinic loop, 
except the two lines (see Fig.~\ref{fig1}): 
\be  
a_4 = a_1 \qquad \forall \, a_1 \in (-1, 0) \cup (0, \infty), 
\l{b53} 
\ee  
which corresponds to the parameter value $\, e = \frac{2}{3}$, and 
\be  
a_4 = -\, \dss\frac{1}{2}\, a_1 \qquad \forall \, a_1 \in (0, \infty), 
\l{b54} 
\ee  
which corresponds to $\, e \rightarrow \pm \infty$. 
It has been shown~\cite{Han1997} that except the above two lines,  
for the parameter values in the shaded area, system (\ref{b49}) 
can have at most $2$ limit cycles near a homoclinic loop
under quadratic perturbation. 

Figure~\ref{fig1} shows 
the $a_1$-$a_4$ parameter plane associated with the reversible 
system (\ref{b7}), where the above mentioned case studies 
are indicated on the line $\, a_4 = \frac{3}{2}\, a_1 \,$ as well 
as in the shaded area.
More precisely, a complete global analysis given in~\cite{YuLi2002}, 
which includes the result in~\cite{DumortierLiZhang1997} as a special case, 
shows that corresponding to each point on the line segment 
$\, a_4 = \frac{3}{2}\, a_1 \ (a_1 > 0)$, the system has one center and 
one saddle point, and has maximal $2$ limit cycles. In~\cite{Han1997} 
it is shown for each point in the shaded area (except the two line segments
$a_4 = a_1 \ ( a_1 > -1)$ and $ a_4 = -\, \frac{1}{2}\,a_1 \ (a_1 >0)$), 
which contains the above line segment, the system has one center and 
one (or three) saddle(s), and has maximal $2$ limit cycles, but restricted to 
local analysis near one homoclinic loop.  
Similarly, a global analysis given in~\cite{IlievLiYu2005}, 
which contains the result in~\cite{DumortierLiZhang1997} as a special case, 
proves that corresponding to each point on the line segment
$\, a_4 = \frac{3}{2}\, a_1 \ (a_1 < -\,1)$, the system has two centers, 
and exhibits maximal $3$ limit cycles around one center. 
The technique of Poincar\'{e} transformation and Picar-Puchs equation, 
used for the above mentioned global analysis on parameter unfolding, 
seems not possible to be generalized to consider general situation for 
arbitrary points in the $a_1$-$a_4$ parameter plane. 
The two particular dash-dotted lines: $a_4 = \frac{1}{3}\,(a_1 - 5) \ 
\forall \, a_1 \! \in \! (-\infty, \, - 1) \cup (-1, \infty)$, 
and $ a_4 = \frac{1}{3} \,(6\,a_1 + 5) \ \forall \, a_1 \! \in \! 
(-\infty,\,-1)$, as well as the five dark circles correspond 
to our results, presented in the next two sections. In particular, 
we will show that there exist $3$ small limit cycles on 
the two dash-dotted lines, and at least $4$ limit cycles for the 
parameter values marked by the five dark circles.

\vspace{0.10in} 
In the following, we will use the perturbed quadratic system (\ref{b44}) 
for our study on 
bifurcation of limit cycles. 
Without loss of generality, we may assume (e.g., 
see~\cite{DumortierLiZhang1997}) 
that $\, a_{01} = a_{20} = a_{11} = a_{02} = b_{10} = b_{20} = b_{02} = 0$.  
Thus, system (\ref{b44}) is reduced to 
\be  
\begin{array}{rl} 
\dss\frac{dx}{dt} = & y \, (1 + a_1\, x) 
+ \varepsilon \, a_{10}\, x , \\ [1.5ex]  
\dss\frac{dy}{dt}= & -\, x + x^2 + a_4\, y^2 
+ \varepsilon \, ( b_{01}\, y + b_{11}\, x\, y ), 
\end{array}
\l{b55} 
\ee 
where $\, a_1 < -\,1 \,$ and $\, 0 < \varepsilon \ll 1$.

\section{Hopf bifurcation associated with the two centers}

\setcounter{equation}{0}
\renewcommand{\theequation}{3.\arabic{equation}}

In this section, we study Hopf bifurcation of system (\ref{b55}) 
from two centers $(0,0)$ and $(1,0)$, leading to bifurcation of multiple 
limit cycles. The result is summarized in the following theorem. 
 
\vspace{0.1in}
\noindent
{\bf Theorem 2.1.} \ {\it
When $\, a_1 < -\,1$, the quadratic near-integrable system (\ref{b55}) 
can have small limit 
cycles bifurcating from the two centers $(0,0)$ and $(1,0)$ with 
distributions: $(3,0)$, $(0,3)$, $(2,0)$, $(0,2)$ and $(1,1)$. 
$(2,1)$- or $(1,2)$-distribution does not exist. 
}

\vspace{0.10in}
\noindent
{\bf Proof.} \
Consider system (\ref{b55}) for $\, a_1 < -\,1$. 
The system (\ref{b55})$_{\varepsilon=0}$ is a reversible integrable system. 
In order to compute the Melnikov function near the two 
centers $(0,0)$ and $(1,0)$, we need 
transform system (\ref{b55})$_{\varepsilon=0}$ to a Hamiltonian system.  
The integrating factor $\,\gamma \,$ is given in (\ref{b14}). 
Now, introducing $\, d t = \gamma\, d \tau\,$ into (\ref{b55}) yields 
the perturbed Hamiltonian system: 
\be  
\begin{array}{rl} 
\dss\frac{d x}{d \tau} = & \gamma\, (y + a_1\, x\, y ) 
+ \varepsilon \, \gamma \, a_{10}\, x , \\[1.5ex] 
\dss\frac{d y}{d \tau} = & \gamma \, (-\, x + x^2 + a_4\, y^2 ) 
+ \varepsilon \, \gamma \, ( b_{01}\, y + b_{11}\, x\, y ), 
\end{array}
\l{b56} 
\ee 
with the Hamiltonian of (\ref{b56})$_{\varepsilon=0}$, given by 
\be 
H(x,y) = \dss\frac{1}{2}\, 
{\rm sign}(1+a_1 x) \, |1+a_1 x|^{- \frac{2 a_4}{a_1}}  
\left[ y^2 + \dss\frac{(1+a_1-a_4)\, (1 + 2\, a_4\, x)}{a_4\, (a_1 - a_4)\, 
(a_1 - 2\, a_4)} - \dss\frac{x^2}{a_1 - a_4} \right],  
\l{b57} 
\ee
for $\, a_4 \ne 0, \ a_1 \ne a_4 , \ a_1 \ne 2\, a_4$. 
The cases $\, a_4 = 0$, $a_1 = a_4 \,$ or $\, \ a_1 = 2\, a_4\, $ 
will not be considered in this paper.

\begin{figure}[!h]
\vspace{-0.60in}
\hspace{0.00in}
\resizebox{1.10\textwidth}{!}{\includegraphics{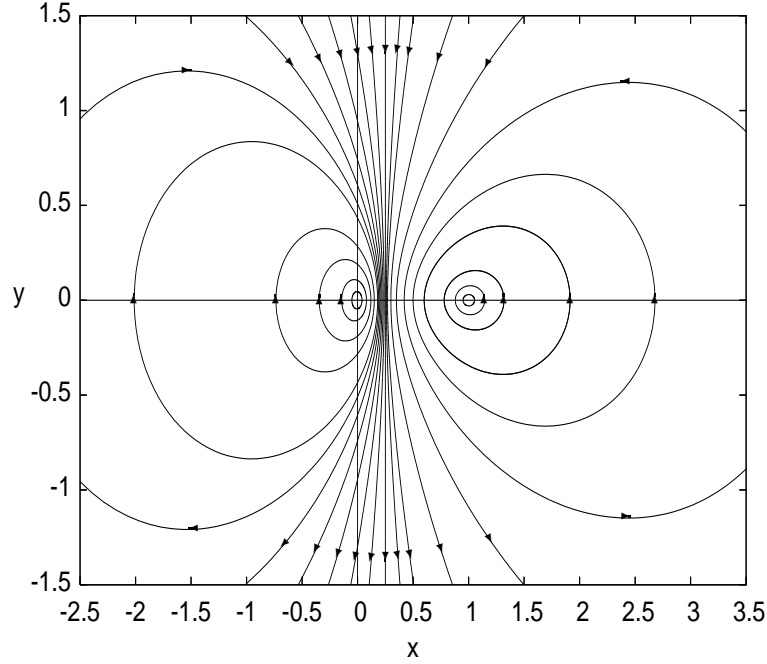}} \hspace{-0.3in}

\vspace{-6.20in} 
\caption{A phase portrait of the reversible system (\ref{b7}) with two 
centers for $a_1=-3,\ a_4 = -\, \frac{8}{3}$.} 
\label{fig2} 
\vspace{0.10in} 
\end{figure}

Note that 
\be 
\begin{array}{ll}  
h_{00} = 
H(0,0) = \dss\frac{1+a_1-a_4}{2\, a_4\, (a_1-a_4)\,(a_1 - 2 a_4)},  
& {\rm for} \ \ 1+a_1\,x > 0, \\[2.0ex] 
h_{10} = H(1,0) = -\, \dss\frac{(a_1+1)\,(a_4+1)}
{2\, a_4\, (a_1-a_4)\,(a_1 - 2 a_4)}\, (-1-a_1)^{- \frac{2 a_4}{a_1}}, 
\quad & {\rm for} \ \ 1 + a_1\,x < 0. 
\end{array}
\l{b58} 
\ee   
Since in this paper, we concentrate on the case that system 
(\ref{b55})$_{\varepsilon=0}$ has only two centers, we assume 
$\, a_1 < -\,1, \ a_4 < 0$. Thus, 
$$ 
\lim_{x \rightarrow -\frac{1}{a_1}^-} H(x,y) = + \,\infty \quad {\rm and} 
\quad \lim_{x \rightarrow -\frac{1}{a_1}^+} H(x,y) = -\,\infty. 
$$ 
It is easy to see from system (\ref{b55}) that the trajectories of 
(\ref{b55})$_{\varepsilon=0}$ rotate around the center $(0,0)$ in 
the clock-wise direction, while rotate around the center $(1,0)$ in 
the counter clock-wise direction, as shown in Fig.~\ref{fig2}. 
Thus, the values of $h$ in $H(x,y) = h$ are taken from the two intervals: 
$\, h \in (h_{00}, \infty) \,$ for $\, 1 + a_1\,x > 0$, and 
$\, h \in (-\infty, h_{10}) \,$ for $\, 1 + a_1\,x < 0$. 
It should be noted that $\, h_{00} \,$ is not necessarily larger than 
$\, h_{10}$. The analyses on the two half-plane in the 
$x$-$y$ plane (see Fig.~\ref{fig2}), 
divided by the singular line $\, 1+a_1 \,x =0$, are independent.

Next, introduce  \vspace{-0.08in}  
\be 
L_{h}: H(x,y) = h \ \left\{ 
\begin{array}{ll} 
h \in (h_{00}, \infty), & {\rm for} \ \ 1+a_1\,x > 0, \\[1.0ex] 
h \in (-\infty, h_{10}), \quad & {\rm for} \ \ 1+a_1\,x < 0, 
\end{array}
\right. 
\l{b59} 
\ee 
and define the Melnikov function: \vspace{-0.05in}  
\be  
\begin{array}{ll}  
M (h, a_{ij}, b_{ij}) = 
\dss\oint_{L_{h}} q(x,y, b_{ij}) \, dx -  p(x,y, a_{ij}) \, dy, 
\end{array} 
\l{b60} 
\ee 
where $\, p(x,y,a_{ij}) = \gamma\, a_{10} \, x \,$ 
and $\, q(x,y,b_{ij}) = \gamma \, ( b_{01} + b_{11}\, x)\,y $. 
Using the results in~\cite{Han2006,Han2000,HC2000}, we can 
expand $\, M \,$ near $\, h = h_{00} \,$ and $\, h= h_{10} \,$ as 
\be
\begin{array}{rl}  
M_0(h, a_{ij}, b_{ij}) = & \mu_{00} \, ( h- h_{00}) 
+ \mu_{01} \, (h-h_{00})^2 
+ \mu_{02} \, (h-h_{00})^3 
\\ [1.5ex]  
& +\, \mu_{03} \, (h-h_{00})^4 
+ O((h-h_{00})^5), \quad  {\rm for} \ \ 0 < h - h_{00} \ll 1, 
\\[1.5ex] 
M_1(h,  a_{ij}, b_{ij}) = & \mu_{10} \, ( h_{10} -h) 
+ \mu_{11} \, (h_{10}-h)^2 
+ \mu_{12} \, (h_{10}-h)^3 
\\ [1.5ex]  
& +\, \mu_{13} \, (h_{10}-h)^4 
+ O((h_{10}-h)^5), \quad  {\rm for} \ \ 0 < h_{10}-h \ll 1, 
\end{array} 
\l{b61}
\ee
where the coefficients $\, \mu_{ij}, \ i=0,\,1; \ j=0,\,1,\,2,\, \cdots \,$ 
can be obtained by using 
the Maple programs developed in~\cite{HanYangYu2010} as follows: 
\begin{eqnarray*} 
\mu_{00} &\!\!\!=\!\!\!& 2\, \pi \, ( a_{10} + b_{01} ),  \\[-0.0ex] 
\mu_{01} &\!\!\!=\!\!\!& 
\dss\frac{\pi}{12} \Big[ 
( 10 - 13\,a_1 - 14 \, a_4 + 13 \, a_1^2 + 7 \, a_1 \, a_4 
- 20 \, a_4^2)\, a_{10}  \\ [-0.5ex] 
&\!\!\! \!\!\!& 
\quad \ 
+\, ( 10 - a_1 + 10 \, a_4 + a_1^2 - 5 \, a_1 \, a_4 + 4 \, a_4^2)\, b_{01}  
+ 12\, (1+a_4)\, b_{11} \Big], \\[-0.0ex] 
\mu_{02} &\!\!\!=\!\!\!& 
\dss\frac{\pi}{864} \Big[ 
(1540-980 a_1-280 a_4+861 a_1^2-1512 a_1 a_4-3948 a_4^2-626 a_1^3 
+ 1566 a_1^2 a_4 \\ [-1.0ex]  
&\!\!\! \!\!\!& 
\qquad \ + 1620 a_1 a_4^2 \!-\! 4432 a_4^3 \!+\! 313 a_1^4 
\!-\!  1018 a_1^3 a_4 \!-\! 279 a_1^2 a_4^2 \!+\! 3080 a_1 a_4^3 
\!-\! 2096 a_4^4) \, a_{10} \qquad \\ [-0.0ex] 
&\!\!\! \!\!\!& 
\qquad +\, (1540+700 a_1+3080 a_4+21 a_1^2+168 a_1 a_4+2772 a_4^2-2 a_1^3 
+ 126 a_1^2 a_4 \\  [-0.0ex] 
&\!\!\! \!\!\!& 
\qquad \quad \ -\, 828 a_1 a_4^2 +1424 a_4^3 +a_1^4 -58 a_1^3 a_4 
+369 a_1^2 a_4^2 - 712 a_1 a_4^3 +400 a_4^4)\, b_{01}  \\ [-0.0ex] 
&\!\!\! \!\!\!& 
\qquad \ + \, 24 \,b_{11} \, (1+a_4)\,
(70+35 a_1+70 a_4+a_1^2-17 a_1 a_4+52 a_4^2) \, b_{11} \Big],  \\ 
\mu_{03} &\!\!\!=\!\!\!& 
\dss\frac{\pi}{622080} \Big[  
(3403400 \!-\! 300300 a_1 \!+\! 3003000 a_4 \!+\! 690690 a_1^2 
\!-\! 4984980 a_1 a_4 \!-\! 7327320 a_4^2 \qquad 
\\ [-1.0ex] 
&\!\!\! \!\!\!& 
\qquad \qquad -\,500885 a_1^3 +3314850 a_1^2 a_4 -4430580 a_1 a_4^2 
-17811640 a_4^3 +323121 a_1^4 
\\[-0.0ex] 
&\!\!\! \!\!\!& 
\qquad \qquad 
-\, 2444439 a_1^3 a_4 \!+\! 4201218 a_1^2 a_4^2 \!+\! 5794692 a_1 a_4^3 
\!-\! 18033936 a_4^4 \!-\! 168603 a_1^5
\\ [-0.0ex] 
&\!\!\! \!\!\!& 
\qquad \qquad 
+\, 1420500 a_1^4 a_4 -3253551 a_1^3 a_4^2 -1296282 a_1^2 a_4^3
+12107904 a_1 a_4^4 
\\ [-0.0ex] 
&\!\!\! \!\!\!& 
\qquad \qquad 
-\, 10462368 a_4^5 +56201 a_1^6 -520311 a_1^5 a_4 +1471287 a_1^4 a_4^2
-407053 a_1^3 a_4^3 
\\ [-0.0ex] 
&\!\!\! \!\!\!& 
\qquad \qquad 
-\, 4589772 a_1^2 a_4^4
+7149264 a_1 a_4^5
-3159616 a_4^6
)\, a_{10} 
\\ [-0.0ex] 
&\!\!\! \!\!\!& 
\qquad \quad \ \   
+\, (3403400 +3303300 a_1 +10210200 a_4 +690690 a_1^2 +5825820 a_1 a_4 
\\ [-0.5ex] 
&\!\!\! \!\!\!& 
\qquad \qquad \ \ \,  
+ 14294280 a_4^2 \!+\! 11935 a_1^3 \!+\! 404250 a_1^2 a_4 
\!+\!  2721180 a_1  a_4^2 \!+\! 12236840 a_4^3 
\\ [-0.0ex] 
&\!\!\! \!\!\!& 
\qquad \qquad \ \ \,
-699 a_1^4 -11379 a_1^3 a_4 +262458 a_1^2 a_4^2 -1891308 a_1 a_4^3
+6994704 a_4^4 
\\ [-0.0ex] 
&\!\!\! \!\!\!& 
\qquad \qquad \ \ \,
+417 a_1^5
+1380 a_1^4 a_4 -149091 a_1^3 a_4^2 +1121838 a_1^2 a_4^3 -2964576 a_1  a_4^4
\\ [-0.0ex] 
&\!\!\! \!\!\!& 
\qquad \qquad \ \ \,
+2670432 a_4^5 
-139 a_1^6
-291 a_1^5 a_4
+46227 a_1^4 a_4^2
-366193 a_1^3 a_4^3 
\\ [-0.0ex] 
&\!\!\! \!\!\!& 
\qquad \qquad \ \ \,
+1076988 a_1^2  a_4^4
-1335216 a_1 a_4^5
+578624 a_4^6
)\, b_{01} 
\\ [-0.0ex] 
&\!\!\! \!\!\!& 
\qquad \quad \ \  
+\, ( 3603600 +3603600 a_1 +10810800 a_4 +790020 a_1^2 + 6597360 a_1 a_4
\\ [-0.0ex] 
&\!\!\! \!\!\!& 
\qquad \qquad \ \ \,
+15024240 a_4^2 \!+\! 12600  a_1^3 \!+\! 480060 a_1^2 a_4 
\!+\! 3764880 a_1 a_4^2 \!+\! 12514320 a_4^3 
\\ [-0.0ex] 
&\!\!\! \!\!\!& 
\qquad \qquad \ \ \,
+180 a_1^4 -10800 a_1^3 a_4 +11340 a_1^2 a_4^2 -618480  a_1 a_4^3
+6566400 a_4^4 
\\ [-0.0ex] 
&\!\!\! \!\!\!& 
\qquad \qquad \ \ \,
+180 a_1^4 a_4
\!-\! 23400 a_1^3 a_4^2
\!+\! 321300 a_1^2 a_4^3
\!-\! 1389600 a_1 a_4^4 
\!+\! 1869120 a_4^5
) \, b_{11}  
\Big], 
\\ [-0.5ex] 
&\!\!\! \vdots \!\!\!&  
\end{eqnarray*} 
and 
\begin{eqnarray*} 
\mu_{10} &\!\!\!=\!\!\!& 2\, \pi \, (-1-a_1)^{3/2} 
\Big[ (1-2\,a_4)\, a_{10} + (1+a_1)\, (b_{01} + b_{11} ),  
\\[-0.0ex] 
\mu_{11} &\!\!\!=\!\!\!& 
\dss\frac{\pi}{12} \ (-1-a_1)^{-\, \frac{2\,(a_1-a_4)}{a_1}} \\ 
&\!\!\! \!\!\!& \times 
\Big[ 
(10 +33 a_1 -6 a_4 +36 a_1^2 -21 a_1 a_4 -24 a_1^2 a_4
+30 a_1 a_4^2 -8 a_4^3)\, a_{10}  
\\ [-0.5ex] 
&\!\!\!  \!\!\!&   
\quad  \ 
+ \, (1+a_1)\, (10 +21 a_1 -10 a_4 + 12 a_1^2 -15 a_1 a_4 +4 a_4^2) \, b_{01}  
\\ [-0.0ex] 
&\!\!\!  \!\!\!&   
\quad  \  
- \, (1+a_1)\,(1+a_4)\, (2 +3\, a_1 -4 \,a_4) \, b_{11} 
\Big],
\\[-0.0ex] 
\mu_{12} &\!\!\!=\!\!\!& 
\dss\frac{\pi}{864} \ (-1-a_1)^{-\, \frac{(5\,a_1-8\,a_4)}{2\,a_1}} \\ 
&\!\!\! \!\!\!& \times 
\Big[ (1540
+7140 a_1
-2800 a_4
+13041 a_1^2
-11592 a_1 a_4
+2212 a_4^2
+11448 a_1^3 
\\ [-0.5ex] 
&\!\!\!  \!\!\!&   
\quad  \  \
-18072 a_1^2 a_4
\!+\!8628 a_1 a_4^2
\!-\! 1112 a_4^3
\!+\! 752 a_4^4
\!-\! 12024 a_1^3 a_4
\!+\! 12213 a_1^2 a_4^2 
\!-\! 5232 a_1 a_4^3 \qquad 
\\ [-0.0ex] 
&\!\!\!  \!\!\!&   
\quad  \  \
+4320 a_1^4
-1728 a_1^4 a_4
+6192 a_1^3 a_4^2
-7938 a_1^2 a_4^3
+4272 a_1 a_4^4
-800 a_4^5
)\, a_{10} 
\\ [-0.0ex] 
&\!\!\!  \!\!\!&   
\quad  \  
+ \, (1\!+\!a_1)\,( 1540
\!+\! 5460 a_1
\!-\! 3080 a_4
\!+\! 7161  a_1^2
\!-\! 9072 a_1 a_4
\!+\! +2772 a_4^2
\!+\! 4104 a_1^3 
\\ [-0.0ex] 
&\!\!\!  \!\!\!&   
\hspace{1.00in} 
-\,9030 a_1^2 a_4
\!+\! 6372 a_1 a_4^2
\!-\! 1424 a_4^3
\!+\! 864 a_1^4
\!-\! 3096 a_1^3 a_4
\!+\! 3969 a_1^2 a_4^2 
\\ [-0.0ex] 
&\!\!\!  \!\!\!&   
\hspace{1.00in} 
-\,2136 a_1 a_4^3
+400 a_4^4
)\, b_{01} 
\\ [-0.0ex] 
&\!\!\!  \!\!\!&   
\quad  \  
- \, (1\!+\!a_1) (1\!+\!a_4) \,
(140
+420 a_1
-420 a_4
+423 a_1^2
-996 a_1 a_4
+576 a_4^2 
\\[-0.0ex]  
&\!\!\!  \!\!\!&   
\hspace{1.45in} 
+\,144 a_1^3
-633 a_1^2 a_4
+888 a_1 a_4^2
-400 a_4^3
)\, b_{11} 
\Big],  
\\[-0.0ex] 
\mu_{13} &\!\!\!=\!\!\!& 
\dss\frac{\pi}{1244160} \ (-1-a_1)^{-\, \frac{(2\,(a_1-3\,a_4)}{a_1}} \\ 
&\!\!\! \!\!\!& \times 
\Big[ (3403400
\!+\! 20720700 a_1
\!-\! 9809800 a_4
\!+\! 53243190 a_1^2
\!-\! 54234180 a_1 a_4
\!+\! 13093080 a_4^2  
\\[-1.0ex] 
&\!\!\! \!\!\!&  
\quad \  
+\,74334645 a_1^3
-123735150 a_1^2 a_4
+65571660 a_1 a_4^2
-10776920 a_4^3
+60023916 a_1^4
\\[-0.0ex] 
&\!\!\! \!\!\!&  
\quad \  
-\,147900519 a_1^3 a_4
\!+\! 131934978 a_1^2 a_4^2
\!-\! 49682268 a_1 a_4^3
\!+\! 6439744 a_4^4
\!+\! 27002160 a_1^5
\\[-0.0ex] 
&\!\!\! \!\!\!&  
\quad \  
-\,95460120 a_1^4 a_4
\!+\! 132380865 a_1^3 a_4^2
\!-\! 89408610 a_1^2 a_4^3
\!+\! 29027880 a_1 a_4^4
\!-\! 3527040 a_4^5
\\[-0.0ex] 
&\!\!\! \!\!\!&  
\quad \  
+\,5443200 a_1^6
\!-\! 28946160 a_1^5 a_4
\!+\! 63998532 a_1^4 a_4^2
\!-\! 74879613 a_1^3 a_4^3
\!+\! 48498336 a_1^2 a_4^4
\\[-0.0ex] 
&\!\!\! \!\!\!&  
\quad \  
-\,16296336 a_1 a_4^5
+2181248 a_4^6
-1555200 a_1^6 a_4
+9603360 a_1^5 a_4^2
-24061752 a_1^4 a_4^3
\\[-0.0ex] 
&\!\!\! \!\!\!&  
\quad \  
+\,31232358 a_1^3 a_4^4
-22072536 a_1^2 a_4^5
+8011296 a_1 a_4^6
-1157248 a_4^7
) \, a_{10} 
\\ [-0.0ex] 
&\!\!\!  \!\!\!&   
\quad    
+ \, (1\!+\!a_1)\,(3403400 
\!+\! 17117100 a_1
\!-\! 10210200 a_4
\!+\! 35225190 a_1^2
\!-\! 45225180 a_1 a_4 
\\ [-0.0ex] 
&\!\!\!  \!\!\!&   
\hspace{0.90in} 
+14294280 a_4^2
+37785825 a_1^3
-79202970 a_1^2 a_4
+54455940 a_1 a_4^2
\\ [-0.0ex] 
&\!\!\!  \!\!\!&   
\hspace{0.90in} 
-12236840 a_4^3
+22125636 a_1^4
-68371209 a_1^3 a_4
+77864598 a_1^2 a_4^2
\\ [-0.0ex] 
&\!\!\!  \!\!\!&   
\hspace{0.90in} 
-38601828 a_1 a_4^3
+6994704 a_4^4
+6629040 a_1^5
-28984608 a_1^4 a_4
\\ [-0.0ex] 
&\!\!\!  \!\!\!&   
\hspace{0.90in} 
+49687587 a_1^3 a_4^2
-41614974 a_1^2 a_4^3
+16953984 a_1 a_4^4
-2670432 a_4^5
\\ [-0.0ex] 
&\!\!\!  \!\!\!&   
\hspace{0.90in} 
+777600 a_1^6
-4801680 a_1^5 a_4
+12030876 a_1^4 a_4^2
-15616179 a_1^3 a_4^3
\\ [-0.0ex] 
&\!\!\!  \!\!\!&   
\hspace{0.90in} 
+11036268 a_1^2 a_4^4
-4005648 a_1 a_4^5
+578624 a_4^6
) \,b_{01} 
\\ [-0.0ex] 
&\!\!\!  \!\!\!&   
\quad    
- \, (1\!+\!a_1)(1\!+\!a_4)\,(
200200
\!+\! 900900 a_1
\!-\! 800800 a_4
\!+\! 1600830 a_1^2
\!-\! 3132360 a_1 a_4
\\ [-0.0ex] 
&\!\!\!  \!\!\!&    
\hspace{1.40in} 
+1530760 a_4^2
+1397655 a_1^3
-4596480 a_1^2 a_4
+5008500 a_1 a_4^2
\\ [-0.0ex] 
&\!\!\!  \!\!\!&    
\hspace{1.40in} 
-1808240 a_4^3
+594864 a_1^4
-3001266 a_1^3 a_4
+5594022 a_1^2 a_4^2
\\ [-0.0ex] 
&\!\!\!  \!\!\!&    
\hspace{1.40in} 
-4568112 a_1 a_4^3
+1379936 a_4^4
+97200 a_1^5
-736776 a_1^4 a_4
\\ [-0.0ex] 
&\!\!\!  \!\!\!&    
\hspace{1.40in} 
+2162079 a_1^3 a_4^2
\!-\! 3080268 a_1^2 a_4^3
\!+\! 2136528 a_1 a_4^4
\!-\! 578624 a_4^5 
)\, b_{11} 
\Big], 
\\[-0.5ex]  
&\!\!\! \vdots \!\!\!&   
\end{eqnarray*}

\noindent 
{\bf Remark 2.2.} \ The coefficients $\, \mu_{0j}\,$ listed above 
are applicable as long as $(0,0)$ is a center, 
and the coefficients $\, \mu_{1j}\,$ are applicable 
as long as $(1,0)$ is a center, regardless the number 
and distribution of the system's singularities. 
Therefore, for each point on the whole line $\, a_4 = \frac{1}{3}\,(a_1 - 5)$ 
(see Fig.~\ref{fig1}), 
there always exist $3$ small limit cycles bifurcating 
from the center $(0,0)$, no matter whether the system has two centers, 
or one center and three saddle points, or one center and one saddle point. 
For each point on the line segment $\, a_4 =  \frac{1}{3}\,(6\,a_1 + 5) 
\ (a_1 < -\,1)$, the system can have $3$ limit cycles bifurcating from 
the center $(1,0)$. This indicates that the results given 
in~\cite{DumortierLiZhang1997,Peng2002,Han1997} showing that the reversible 
near-integrable systems with one center and one saddle point can have 
maximal $2$ limit cycles is conservative, since on the part of 
the line $\,a_4 = \frac{1}{3}\,(a_1 - 5) \,$ in the first quadrant 
($a_1 > 5$) such a system can have at least $3$ limit cycles.

\vspace{0.10in}
First, we consider the maximal number of limit cycles which can 
bifurcate from the center $(0,0)$. Setting $\, \mu_{00} = 0\,$ yields 
\vspace{-0.10in}  
\be 
b_{01} = -\, a_{10}, 
\l{b62} 
\ee 
 
\vspace{-0.10in}  
\noindent 
and then we have \vspace{-0.05in}  
\be  
\mu_{01} = \pi \, \Big[ (a_1-1-a_4)(a_1+2 a_4) a_{10} +(1+a_4)\, b_{11} \Big].  
\l{b63} 
\ee 
In order to have $\, \mu_{01} = 0$, we suppose $\, a_4 \ne -\,1 \,$ and 
choose 
\be  
b_{11} = -\, \dss\frac{ (a_1-1-a_4)(a_1+2 a_4)} {1+a_4} \ a_{10} . 
\l{b64} 
\ee 
Then, $\, \mu_{02} \,$ and $\, \mu_{03} \,$ are simplified to 
\be 
\begin{array}{ll}  
\mu_{02} = \dss\frac{\pi}{3}\, a_1\, (a_1-a_4)\, (a_1+ 2\,a_4)  
(a_1 - 3 \,a_4 - 5) \, a_{10},  \\[1.5ex]  
\mu_{03} = -\, \dss\frac{\pi}{144}\, a_1 (a_1\!-\!a_4) (a_1\!+\!2 a_4)
( 770
+105 a_1
+1400 a_4
+42 a_1^2
-434 a_1 a_4
\\ [1.5ex] 
\hspace{2.30in} +1274 a_4^2
-13 a_1^3
+128 a_1^2 a_4
-415 a_1 a_4^2
+444 a_4^3
)\,a_{10}. 
\end{array} 
\l{b65} 
\ee
There are five choices for $\, \mu_{02} = 0 $. Except the 
choice $\, a_1 - 3\, a_4 - 5 = 0$, all other choices lead to 
$\, \mu_{0i} = 0, \ i=3,\,4, \cdots $. 
Thus, letting 
\be  
a_4 = \dss\frac{1}{3}\, ( a_1 - 5), 
\l{b66} 
\ee 
which implies $\, a_1 \ne 2 \,$ when $\, a_4 \ne -\,1$. 
Since we assume $\, a_1 < -\,1$, for this case (i.e., when 
the condition (\ref{b66}) holds), $\,a_4 \ne -\,1 \,$ is guaranteed. 
Then, we have 
\begin{eqnarray*}
\mu_{03} 
&\!\!\!=\!\!\!&
- \dss\frac{25 \pi}{162} \, a_1\, (a_1+1)\,(a_1-2)^2 \, 
(2\,a_1 + 5)\, a_{10}, \\ [0.0ex] 
\mu_{04} 
&\!\!\!=\!\!\!&
- \dss\frac{5 \pi}{8748}\, a_1\, (a_1+1)\,(a_1-2)^2 \, 
(2\,a_1 + 5)\, (a_1+4)\, (17\,a_1+518) \, a_{10} \\[0.0ex] 
&\!\!\! \vdots \!\!\!& \\ [0.0ex]  
\mu_{10} 
&\!\!\!=\!\!\!&
-\, \dss\frac{10 \pi}{3} \, (-1-a_1)^{-3/2} \, 
a_1 \, (2\,a_1 + 5) \, a_{10}, \\[0.0ex] 
\mu_{11} 
&\!\!\!=\!\!\!&
\dss\frac{25 \pi}{324} \,  (-1-a_1)^
{-\,\frac{2\,(2\,a_1+5)}{3\, a_1}} a_1 \, (a_1-2)^2 \, 
(2\,a_1 + 5)\, a_{10}, \\[0.0ex]  
&\!\!\! \vdots \!\!\!&
\end{eqnarray*}
implying that in addition we need \vspace{-0.10in}  
\be 
(2\,a_1 + 5)\, a_{10} \ne 0. 
\l{b67} 
\ee 

Under the above conditions (\ref{b62}), (\ref{b64}), (\ref{b66}) and 
(\ref{b67}), we obtain $\, \mu_{00} = \mu_{01} = \mu_{02} = 0$, 
but $\, \mu_{03} \ne 0 , \ \mu_{10} \ne 0$. Hence, at most $3$ small 
limit cycles can bifurcate from the center $(0.0)$ with no limit cycles 
bifurcating from the center $(1,0)$. Further, giving proper perturbations
to the parameters $a_4 $ (or $a_1$), $b_{11}$ and $b_{01}$, we can obtain 
$3$ small limit cycles bifurcating from the origin.  
This shows that the conclusion is true for the case of 
$(3,0)$-distribution. 

Next, consider the $(0,3)$-distribution. 
Similarly, letting $\, \mu_{10} = 0\,$ yields 
\be  
b_{01}= -\, b_{11} + \dss\frac{ 2\,a_4-1}{1+ a_1} \, a_{10}. 
\l{b68} 
\ee 
Then, $\, \mu_{11} \, $ becomes 
\be  
\mu_{11} = \pi (-1 - a_1)^{-\, \frac{2\,(a_1-a_4)}{a_1}} 
\Big[ ( a_1 + 2\, a_4) ( 2\,a_1 - a_4 + 1)\, a_{10} 
- (1+a_1)^2 \, (a_1- a_4 + 1)\, b_{11} \Big].   
\l{b69} 
\ee  
Hence, we set 
\be  
b_{11} = \dss\frac{ ( a_1 + 2\, a_4) ( 2\,a_1 - a_4 + 1)}
{(1+a_1)^2 \, (a_1 - a_4 + 1)} \, a_{10}, \qquad 
( a_1 - a_4 + 1 \ne 0 ),   
\l{b70} 
\ee 
to yield $\, \mu_{11} = 0$, and 
\be 
\begin{array}{ll}  
\mu_{12} = \dss\frac{\pi}{3}\, 
(-1-a_1)^{- \frac{5 a_1 - 8 a_4}{2\,a_1} }
\, a_1\, (a_1-a_4)\, (a_1+ 2\,a_4)  
(6\, a_1 - 3 \,a_4 + 5) \, a_{10},  \\[1.5ex]  
\mu_{13} = \dss\frac{\pi}{288}\, (-1-a_1)^{- \frac{2 (a_1 - 3 a_4)}{a_1} } 
a_1 (a_1\!-\!a_4) (a_1\!+\!2 a_4)
\,( 770
+2205 a_1
-1400 a_4
+2142 a_1^2
\\ [1.5ex] 
\hspace{1.20in} 
-\,3234 a_1 a_4
+1274 a_4^2
-720 a_1^3
-1962 a_1^2 a_4
+1689 a_1 a_4^2
-444 a_4^3
)\,a_{10}. 
\end{array} 
\l{b71} 
\ee
The only choice for $\, \mu_{12} = 0 \,$ is 
$\, 6\, a_1 - 3\, a_4 + 5 = 0$, from which we have 
\be  
a_4 = \dss\frac{1}{3}\, ( 6\, a_1 + 5). 
\l{b72} 
\ee 
This implies that $\, a_1 - a_4 + 1 = -\,(a_1 + \frac{2}{3}) > 0 \,$ 
for $\, a_1 < -\,1$. Further, we obtain 
\begin{eqnarray*} 
\mu_{13} &\!\!\!=\!\!\!& -\, \dss\frac{25 \pi}{324} \, 
(-1-a_1)^{\frac{10 + 11\, a_1}{a_1} }
a_1\, (3\,a_1+2)^2 \, (3\,a_1 + 5)\, a_{10}, \\ [0.0ex] 
\mu_{14} &\!\!\!=\!\!\!& - \dss\frac{5 \pi}{17496}\,
(-1-a_1)^{\frac{80 + 87\, a_1}{6\, a_1} }
 a_1\, (3\,a_1+2)^2\, (3\,a_1 + 5)\, 
(3\, a_1+4) \, (501\,a_1 + 518)\, a_{10} \\[0.0ex] 
&\!\!\! \vdots \!\!\!& \\ [0.0ex]  
\mu_{00} &\!\!\!=\!\!\!& \dss\frac{10 \pi}{3\, (1+a_1)^2} \, 
a_1 \, (3\,a_1 + 5) \, a_{10}, \\[0.0ex] 
\mu_{01} &\!\!\!=\!\!\!& -\, \dss\frac{25 \pi}{324\, (1+a_1)^2} 
\, a_1 \, (3\,a_1+5) \, 
(3\,a_1 + 2)^2 \, a_{10}, \\[0.0ex]  
&\!\!\! \vdots \!\!\!& 
\end{eqnarray*} 
implying that in addition we require \vspace{-0.05in}  
\be 
(3\,a_1 + 5)\, a_{10} \ne 0. 
\l{b73} 
\ee 

Under the above conditions (\ref{b68}), (\ref{b70}), (\ref{b72}) and 
(\ref{b73}), we have $\, \mu_{10} = \mu_{11} = \mu_{12} = 0$, 
but $\, \mu_{13} \ne 0 , \ \mu_{00} \ne 0$. Further, by properly 
perturbing the parameters $a_4 $ (or $a_1$), $b_{11}$ and $b_{01}$, 
we can obtain $3$ small limit cycles bifurcating from the 
center $(1,0)$, but no limit cycles from the origin. 
This proves the case of $(0,3)$-distribution. 

For the case of $(2,0)$-distribution, 
it follows from the conditions (\ref{b62}) 
and (\ref{b64}), and $\, a_4 \ne -\,1\,$ that  
$\, \mu_{00} = \mu_{01} =0$, and 
$$
\begin{array}{ll} 
\mu_{02} = \dss\frac{\pi}{3}\, a_1\, (a_1-a_4)\, (a_1+ 2\,a_4)  
(a_1 - 3 \,a_4 - 5) \, a_{10},  \\[1.5ex]   
\mu_{10} = -\, \dss\frac{2\,\pi}{(1+a_4)\, (-1-a_1)^{3/2}} \, 
a_1 \, (a_1 - a_4)\, ( a_1 + 2\, a_4)\, a_{10}. 
\end{array}
$$  
Thus, $\, \mu_{02} \ne 0 \,$ implies $\, \mu_{10} \ne 0$, 
indicating that the conclusion holds for the case of $(2,0)$-distribution. 
if $\, a_4 = -\,1$. 

When $\, a_4 = -\,1$, (\ref{b63}) becomes 
$$ 
\mu_{01} = \pi \, a_1 \, ( a_1 - 2) \, a_{10} \ne 0 \qquad 
{\rm for} \ \ a_1 < -\,1 \quad {\rm and} \quad a_{10} \ne 0. 
$$ 
Under the conditions $\, b_{01} = -\, a_{10} \,$ and 
$\, a_4 = -\,1$, $\, \mu_{10} \,$ and $\, \mu_{11} \,$ becomes 
\be 
\begin{array}{ll}  
\mu_{10} = - \, 2\, \pi\, (-1-a_1)^{-3/2} \, \Big[ 
(a_1 -2 )\, a_{10} - (1+a_1)\, b_{11} \Big], \\[1.5ex] 
\mu_{11} = \pi \, (-1-a_1)^{-\, \frac{2+ a_1)}{a_1}} \, 
a_1 \, (a_1 -2) \, a_{10} ,  
\end{array}   
\l{b74} 
\ee 
which shows that $\, \mu_{11} \ne 0 \,$ for $\, a_1 < -\,1 \,$ 
and $\, a_{10} \ne 0$. But we can choose 
$$ 
b_{11} = \dss\frac{ a_1 - 2}{ 1+a_1} \, a_{10} 
$$ 
to obtain $\, \mu_{10} =0$. Thus, for this case we have 
a $\, (1,1)$-distribution.

Similarly, for the $(0,2)$-distribution, we use the conditions 
(\ref{b68}) and (\ref{b70}) to obtain 
$$ 
\begin{array}{ll}  
\mu_{12} = \dss\frac{\pi}{3}\, 
(-1-a_1)^{- \frac{5 a_1 - 8 a_4}{2\,a_1} }
\, a_1\, (a_1-a_4)\, (a_1+ 2\,a_4)  
(6\, a_1 - 3 \,a_4 + 5) \, a_{10},  \\[1.5ex]  
\mu_{00} = \dss\frac{ 2 \, \pi}{(1+a_1)^2 \, (a_1 - a_4 + 1)} \, 
a_1\, (a_1 - a_4) \, ( a_1 + 2\, a_4) \, a_{10}. 
\end{array} 
$$ 
This indicates that $\, \mu_{12} \ne 0 \,$ 
implies $\, \mu_{00}  \ne 0$, and so the conclusion for the 
case of $(0,2)$-distribution is also true if $\, a_1 - a_4 + 1 \ne 0$. 

When $\, a_1 - a_4 + 1 = 0$, i.e., $a_4 = a_1 + 1 < 0$, (\ref{b69}) 
is reduced to 
$$ 
\mu_{11} = \pi (-1 - a_1)^{-\, \frac{2\,(a_1-a_4)}{a_1}} 
\, a_1 ( 3\, a_1 + 2)\, a_{10} \ne 0 \quad {\rm for} \ \ 
a_1 < -\,1 \quad {\rm and} \quad a_{10} \ne 0 , 
$$
and $\, \mu_{00} \,$ and $\, \mu_{01} \,$ become 
\be  
\begin{array}{ll}  
\mu_{00} = \dss\frac{2\,\pi}{1+a_1} \Big[ 
(a_1 + 2)\, a_{10} - (1+a_1)\, b_{11} \Big], \\[2.5ex]  
\mu_{01} = -\, \dss\frac{\pi}{1+a_1} \, a_1 \, (3\,a_1 + 2)\, a_{10},  
\end{array} 
\l{b75} 
\ee  
which clearly shows that $\, \mu_{01} \ne 0 \,$ for 
$\, a_1 < -\,1 \,$ and $\, a_{10} \ne 0$. However, we may choose 
$$ 
b_{11} = \dss\frac{ a_1 + 2}{ 1+a_1} \, a_{10} 
$$ 
to obtain $\, \mu_{00} =0$. Thus, for $\, a_1 - a_4 +1 =0$, 
we have a $(1,1)$-distribution. 

Finally, suppose the condition given in (\ref{b62}) is satisfied, i.e., 
$\, b_{01} = -\, a_{10}$, then substituting this into $\, \mu_{10} \,$ 
to solve $\, b_{11} $ to obtain 
\be  
b_{11} = \dss\frac{a_1 + 2\, a_4}{1 + a_1}. 
\l{b76} 
\ee 
Then, under the conditions (\ref{b62}) and (\ref{b76}), 
we obtain 
\be  
\begin{array}{ll}  
\mu_{01} = \dss\frac{\pi}{1+a_1} \, a_1 \, (a_1 - a_4)\, ( a_1 + 2\, a_4) \, 
a_{10}, \\[1.5ex] 
\mu_{11} = -\, \pi \, (-1-a_1)^{-\, \frac{2 \, (a_1-a_4)}{a_1} } 
\, a_1 \, (a_1 - a_4)\, ( a_1 + 2\, a_4) \,
a_{10}, 
\end{array} 
\l{b77} 
\ee 
which shows that $\, \mu_{01} \ne 0 \,$ implies 
$\, \mu_{11} \ne 0$, and thus in general the conclusion is true for 
the case of $(1,1)$-distribution. 

As we have seen in the above analysis, 
if the condition (\ref{b66}), $\, a_4 = \frac{1}{3}\,(a_1 - 5)$, 
is not used, then we can only have $2$ limit cycles 
bifurcating from the origin, but no limit cycles can occur 
from the center $(1,\,0)$. In other words, we can obtain one more limit 
cycle, by using the condition $\, a_4 = \frac{1}{3}\,(a_1 - 5)$,  
only bifurcating from the center $ (0,0)$.
Similarly, if the condition (\ref{b72}), $a_4 =\frac{1}{3}\,(6 a_1 + 5)$, 
is not used, then we can have only $2$ limit cycles
bifurcating from the center $(1,0)$, but no limit cycles can bifurcate 
from the origin. Then, condition $a_4 =\frac{1}{3}\,(6 a_1 + 5)$ can be only 
used to get one more limit cycle around the center $(1,0)$, rather than 
the origin. 
Therefore, $(2,1)$- or $(1,2)$-distribution is not possible.

This completes the proof of Theorem 2.1. 
\put(10,0.5){\framebox(6,7.5)}

\section{Limit cycles bifurcating from closed orbits}

\setcounter{equation}{0}
\renewcommand{\theequation}{4.\arabic{equation}}

In this section, based on the results of the small limit cycles 
obtained in the previous section, 
we want to investigate the possibility of existence of large limit cycles 
by applying the Melnikov function, defined in (\ref{b60}). 
We have the following result.

\vspace{0.1in}
\noindent
{\bf Theorem 4.1.} \ {\it
For the case of bifurcation of small limit cycles from the two centers 
$(0,0)$ and $(1,0)$ with $(3,0)$-distribution 
(respectively, $(0,3)$-distribution) there exists at least one 
large limit cycle near $L_h$ for some 
$\, h \in (- \infty, h_{10})$ (respectively for 
some $\, h \in (h_{00}, \infty)$). 
For the case of limit cycles with $(2,0)$-distribution 
(respectively, $(0,2)$-distribution) there exist at least two 
large limit cycles, one near $L_{h_1}$ for some $\, h_1 \in (- \infty, h_{10})$ 
and one near $L_{h_2}$ for some $\, h_2 \in (h_{00}, \infty)$.   
The corresponding values of the parameters $\,a_1 \,$ and $\, a_4\,$ 
for the existence of $4$ limit cycles  
can appear at least in some regions in the $a_1$-$a_4$ parameter plane. 
}

\vspace{0.1in}
\noindent
{\bf Remark 4.2.} \ Theorem 4.1 gives a positive answer to the open question 
of existence of limit cycles in near-integrable quadratic systems: at least 
$4$ limit cycles can exist. For the case of $(1,1)$-distribution, so far no more 
large limit cycles have been found.    

\vspace{0.1in}
\noindent
{\bf Proof.} \ It follows from (\ref{b60}) with
$$
p(x,y,a_{ij}) = |1+a_1\,x|^{- \frac{a_1+2 a_4}{a_1}} 
\, a_{10} \, x, \quad 
q(x,y,b_{ij}) = |1+a_1\,x|^{- \frac{a_1+2 a_4}{a_1}} 
 \, ( b_{01} + b_{11}\, x)\, y  
$$ 
that 
\ba  
&\!\!\! \!\!\!& M (h,a_1,a_4,a_{10},b_{01},b_{11}) 
\nonumber \\
&\!\!\!=\!\!\!&
\dss\oint_{L_{h}} q(x,y, b_{ij}) \, dx -  p(x,y, a_{ij}) \, dy  
\nonumber \\
&\!\!\!=\!\!\!&
\dss\oint_{L_{h}} q(x,y, b_{ij}) \, dx 
-  \dss\oint_{L_{h}} p(x,y, a_{ij}) \, dy  
\nonumber \\
&\!\!\!=\!\!\!&
\dss\oint_{L_{h}} q(x,y, b_{ij}) \, dx 
+ \dss\oint_{L_{h}} y\, p_x (x,y, a_{ij}) \, dx \qquad  
\nonumber \\
&\!\!\!=\!\!\!&
\dss\oint_{L_{h}} \Big[ q(x,y, b_{ij}) 
+ y\, p_x (x,y, a_{ij}) \Big] \, dx 
\nonumber \\ 
&\!\!\!=\!\!\!&
\dss\oint_{L_{h}}  \left[ |1+a_1\,x|^{- \frac{a_1+2 a_4}{a_1}}\, 
( b_{01} + b_{11}\,x)  
+ |1+a_1\,x|^{- \frac{a_1+2 a_4}{a_1}}\, 
a_{10} \left( 1 \mp \dss\frac{(a_1+2\,a_4)\,x}{  |1+a_1\,x|} \right) 
\right] y \, dx \qquad   \quad 
\nonumber \\
&\!\!\!=\!\!\!&
\dss\oint_{L_{h}}  |1+a_1\,x|^{- \frac{a_1+2 a_4}{a_1}} 
\left[ ( a_{10} + b_{01} ) 
+ b_{11}\, x - {\rm sign}(1 + a_1 \,x) \, 
a_{10}\, (a_1 + 2\, a_4) \, \dss\frac{x}{|1+a_1\,x|} 
\right] y \, dx  
\nonumber \\ [1.0ex]  
&\!\!\!=\!\!\!&
( a_{10} + b_{01} ) \, I_0(h,a_1,a_4) 
+ b_{11} \, I_1(h,a_1,a_4) 
+ a_{10}\, I_2(h,a_1,a_4) 
\nonumber \\ [1.5ex] 
&\!\!\! \equiv \!\!\!&
\left\{  \begin{array}{ll} 
M_0(h,a_1,a_4,a_{10}, b_{01}, b_{11}) \quad 
{\rm for} \ \ h\in (h_{00},\infty),  & {\rm when} \ \ 
1 + a_1\,x > 0, \\[1.5ex] 
M_1(h,a_1,a_4,a_{10}, b_{01}, b_{11}) \quad
{\rm for} \ \ h\in (-\infty, h_{10}), \ \ & {\rm when} \ \
1 + a_1\,x < 0,
\end{array} 
\right. 
\l{b78} 
\ea 
where 
\begin{eqnarray*} 
I_0(h,a_1,a_4) &\!\!\!\!=\!\!\!\!& \dss\oint_{L_h} 
|1+a_1\,x|^{- \frac{a_1+2 a_4}{a_1}}\, y \, dx \\[1.0ex] 
 &\!\!\!\!=\!\!\!\!& \left\{ \begin{array}{lll} 
\hspace{0.05in} 
2 \dss\int_{x_{\min}}^{x_{\max}} 
(1+a_1\,x)^{- \frac{a_1+2 a_4}{a_1}}\, y_+\, dx, 
& \forall \ h \in (h_{00},\infty), \!\!& {\rm when} \ \ 1 + a_1 x > 0, \\[2.5ex] 
\!\!\! -2 \dss\int_{x_{\min}}^{x_{\max}} 
(-1 \!-\! a_1 x)^{- \frac{a_1+2 a_4}{a_1}}\, y_+\, dx, 
& \forall \ h \in (-\infty, h_{10}), \!\!&{\rm when} \ \ 1 + a_1 x < 0; 
\end{array} 
\right.  
\\[1.0ex] 
I_1(h,a_1,a_4)  &\!\!\!\!=\!\!\!\!& \dss\oint_{L_h} 
|1+a_1\,x|^{- \frac{a_1+2 a_4}{a_1}}\, x\, y \, dx \\[0.0ex] 
 &\!\!\!\!=\!\!\!\!& \left\{ \begin{array}{lll} 
\hspace{0.05in} 
2 \dss\int_{x_{\min}}^{x_{\max}} 
(1+a_1\,x)^{- \frac{a_1+2 a_4}{a_1}} x\,y_+ dx, 
& \forall \ h \in (h_{00},\infty), \!\!& {\rm when} \ \ 1 \!+\! a_1 x > 0, 
\\[2.5ex] 
\!\!\! -2 \dss\int_{x_{\min}}^{x_{\max}} 
(-1 \!-\! a_1 x)^{- \frac{a_1+2 a_4}{a_1}} x\,y_+ dx, 
& \forall \ h \in (-\infty, h_{10}),\!\!& {\rm when} \ \ 1 \!+\! a_1 x < 0; 
\end{array} 
\right. 
\qquad  \qquad 
\\[1.0ex] 
I_2(h,a_1,a_4)  &\!\!\!\!=\!\!\!\!& -\, (a_1 + 2\, a_4) \, 
{\rm sign}(1+a_1 x) \dss\oint_{L_h} 
|1+a_1\,x|^{- \frac{2 (a_1+ a_4)}{a_1}}\, x\, y \, dx \\[1.0ex] 
 &\!\!\!\!=\!\!\!\!& \left\{ \begin{array}{lll} 
\!\!\! 
-2 (a_1 \!+\! 2 a_4) \! \dss\int_{x_{\min}}^{x_{\max}} 
\! (1\!+\!a_1 x)^{- \frac{2 (a_1+a_4)}{a_1}} x y_+ dx, 
& \hspace{-0.05in} \forall \, h \! \in \! (h_{00},\infty), \!\!\!
& \hspace{-0.1in} {\rm when} \ 1 \!+\! a_1 x \!>\! 0, 
\\[2.5ex] 
\!\!\! 2 (a_1 \!+\! 2 a_4) \!\! \dss\int_{x_{\min}}^{x_{\max}} 
\!\! (-1 \!-\! a_1 x)^{- \frac{2(a_1+a_4)}{a_1}} x y_+ dx, 
& \!\!\!\! \forall \, h \! \in \! (-\infty, h_{10}),\!\!\!
& {\rm when} \ 1 \!+\! a_1 x \!<\! 0. 
\end{array} 
\right. 
\end{eqnarray*} 
Here, 
\be  
y_+ = \left[ \dss\frac{x^2}{a_1 - a_4} 
- \dss\frac{(1+a_1-a_4)\,(1+2\,a_4\,x)}{a_4\, (a_1-a_4)\,(a_1 - 2\,a_4)} 
+ 2\, h\, {\rm sign}(1 + a_1 \, x) |1+a_1\,x|^{ \frac{2 \, a_4}{a_1}} 
\right]^{1/2},  
\l{b79} 
\ee 
and $\, x_{\min} \,$ and $\, x_{\max} \,$ are solved from the 
equation, $\, y_+ = 0$,  
for $\, h \in (h_{00}, \infty) \,$ when $\, 1 + a_1 \,x > 0$, 
and for $\, h \in (-\infty, h_{10})\,$ when $\, 1 + a_1 \,x < 0$.

Since one can not find the closed form of the integrals 
$\, I_i(h,a_1,a_4), \ i =0,\,1,\,2$, for 
general $\, a_1 \,$ and $\, a_4$, nor the technique of Picard-Fuchs equation 
can be applied here, we shall choose some values for 
$\, a_1 \,$ and $\, a_4\,$ and then find numerical values of the 
integral. We first use the results given in the previous section to determine 
$\, b_{01} $, $\, b_{11}$, and $\, a_4$, and then choose proper values for 
$\, a_1 \,$ to find more limit cycles.

\begin{figure}[!h]
\vspace{0.05in}
\begin{center}
\hspace{0.00in}

\resizebox{0.70\textwidth}{!}{\includegraphics{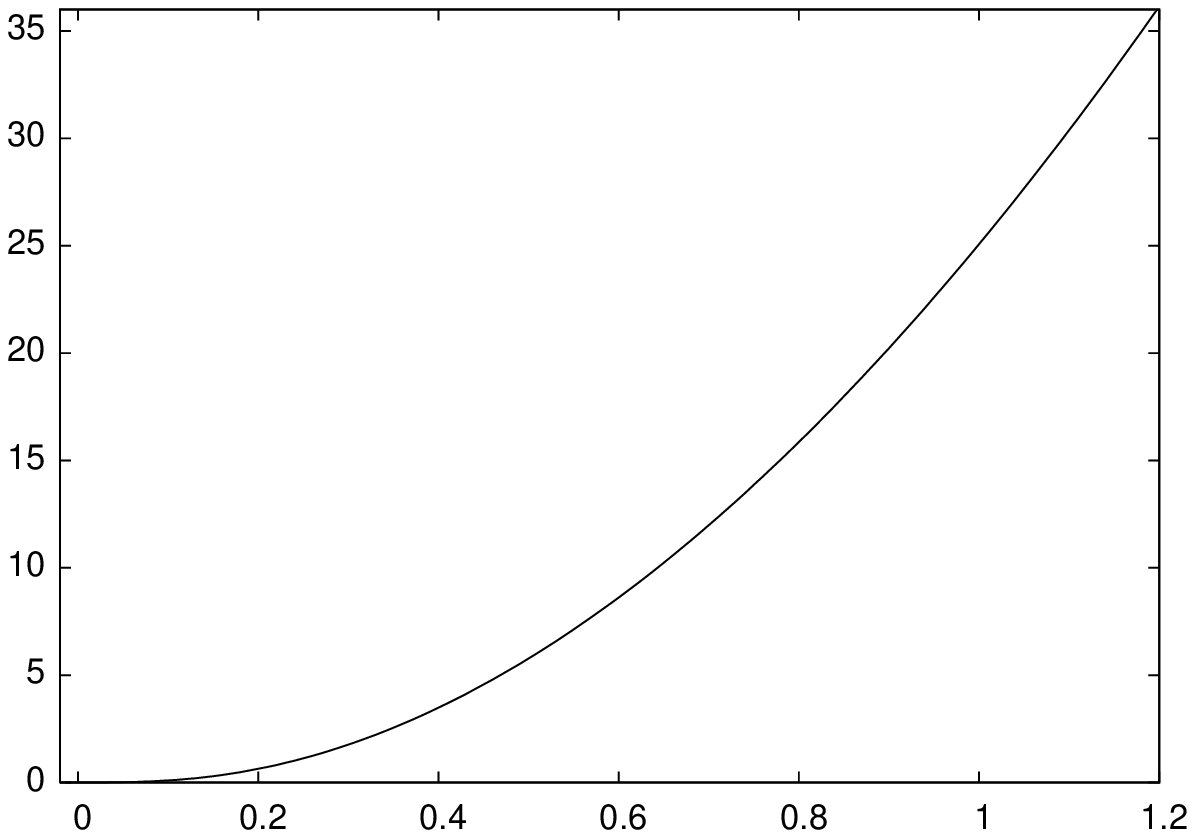}}

\vspace{-2.80in}  
\hspace{-1.20in} 
\resizebox{0.33\textwidth}{!}{\includegraphics{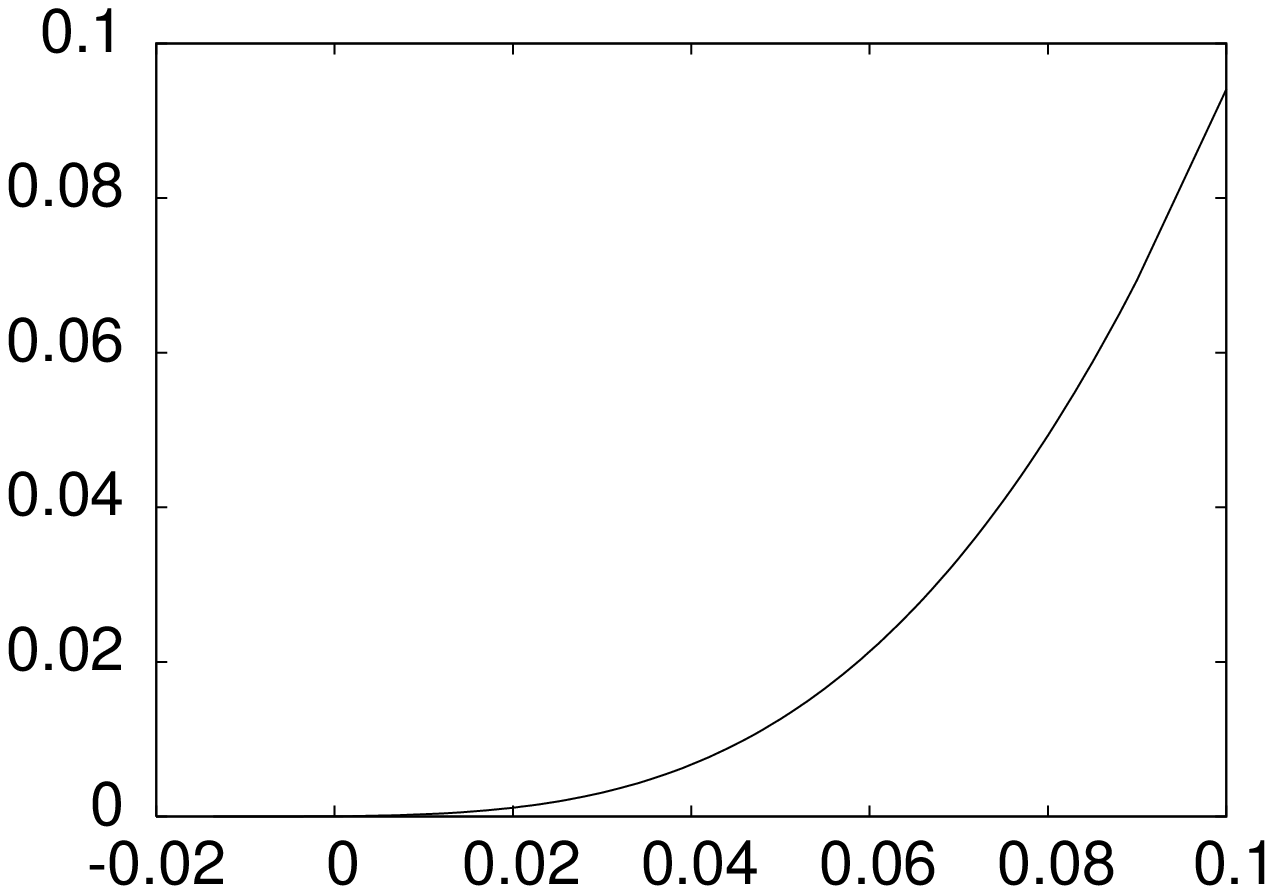}}

\vspace{1.20in}  
\hspace{0.25in}  
$h$ 

\vspace{-2.0in} 
\hspace{-4.6in} 
$M_{00}$  

\vspace{0.0in} 
\hspace{-2.50in} 
$h_{00}$ 

\vspace{1.8in} 
\resizebox{0.70\textwidth}{!}{\includegraphics{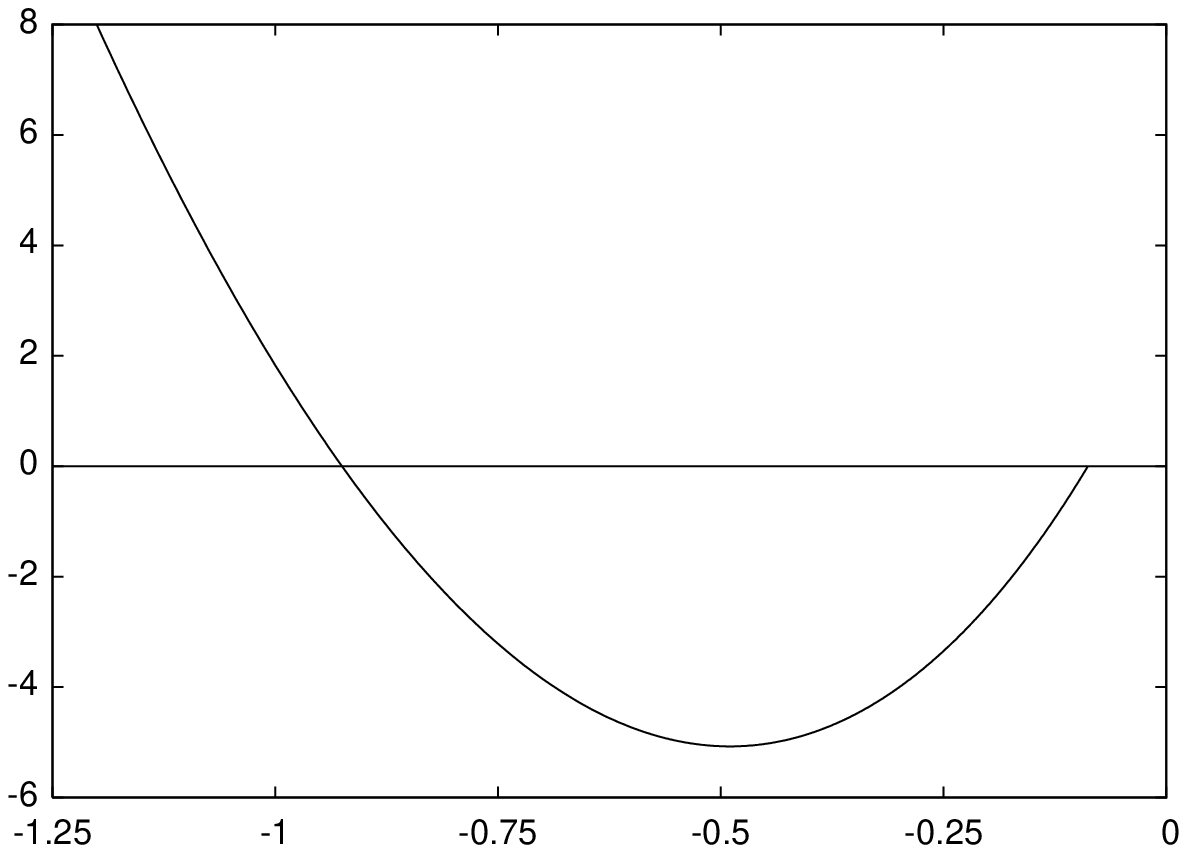}}

\vspace{-0.10in}  
\hspace{0.15in}
$h$ 

\vspace{-1.8in} 
\hspace{-4.6in} 
$M_{10}$  

\vspace{-0.27in} 
\hspace{1.05in} $h_1^*$ \hspace{2.35in} $h_{10}$  

\vspace{-4.95in} \hspace{3.20in} (a) 

\vspace{ 3.40in} \hspace{3.20in} (b) 

\vspace{2.8in} 

\caption{Functions $M_{00}(h) $ and $M_{10}(h) $ 
under the conditions 
$ \mu_{00} = \mu_{01} = \mu_{02} =0 $, $\, \mu_{03} \ne 0$ 
and $ \mu_{10} \ne 0$,   
for $ a_1 = -\,\frac{30}{7}$ and $ a_4 = \frac{1}{3}(a_1-5)
=- \frac{65}{21}$: (a) $\, M_{00}(h) > 0 $ 
for $h \in [h_{00}, \, +\infty)$, with 
$\, h_0  = -\,\frac{441}{32500} \approx -\,0.01357$; 
and (b) $M_{10}(h) $ for $h \in (-\infty, \, h_1]$, 
with $\, h_{10} = -\, \frac{33957}{747500}(\frac{23}{7})^{5/9} 
\approx -\,0.08797$,  
crossing the $h$-axis at $h=h_1^* \in (-\,0.9250363254,\,-\,0.9250363253)$.} 
\label{fig3} 
\end{center}
\vspace{0.30in} 
\end{figure}

\begin{figure}[!h]
\vspace{0.00in}
\begin{center}
\hspace{0.00in}

\resizebox{0.50\textwidth}{!}{\includegraphics{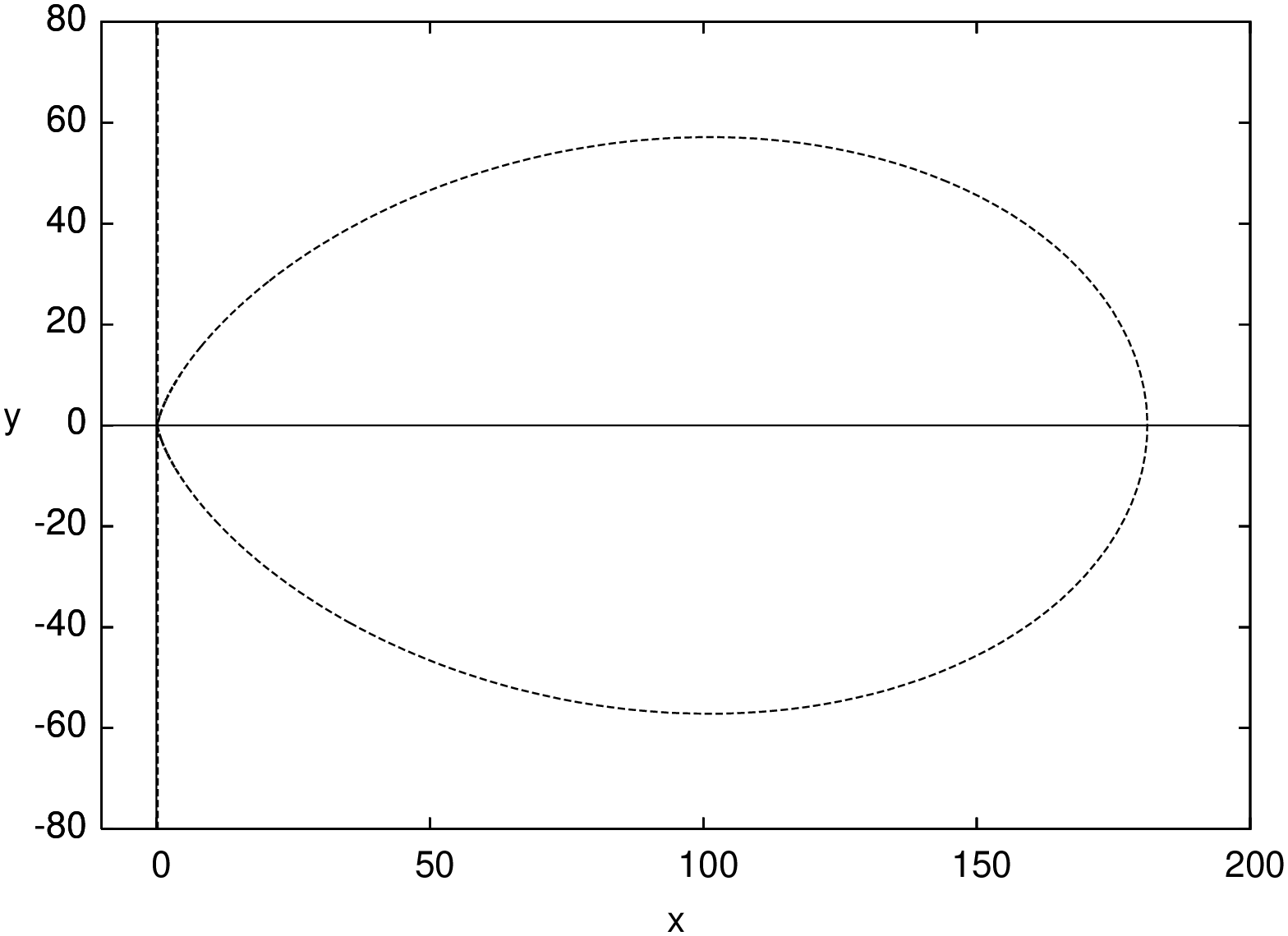}}

\vspace{-1.15in} 
\resizebox{0.50\textwidth}{!}{\includegraphics{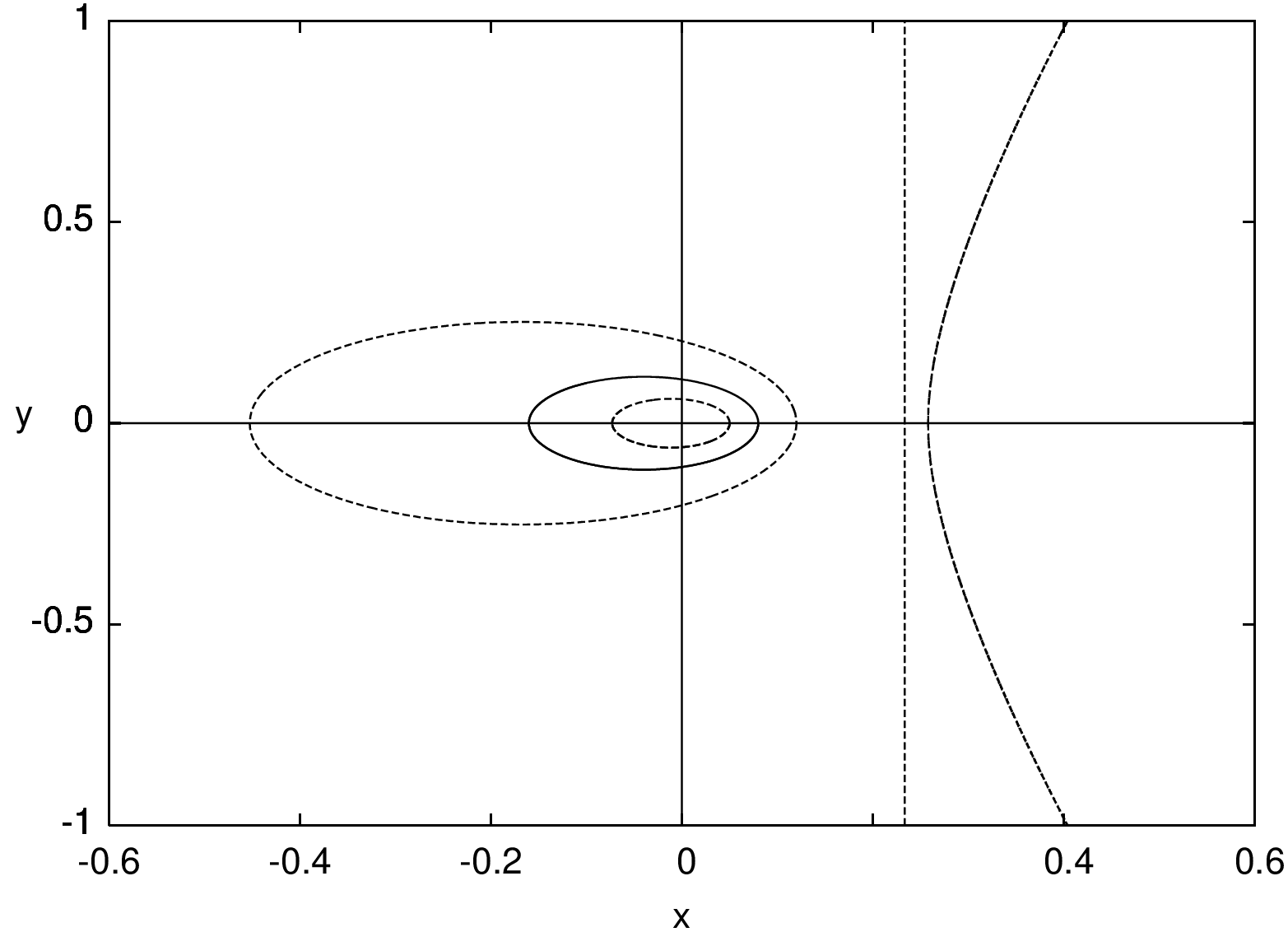}}

\vspace{-1.55in} 
\hspace{1.93in}${\small \frac{7}{30}}$

\vspace{-6.40in} \hspace{3.70in} (a) 

\vspace{ 3.30in} \hspace{3.70in} (b) 

\vspace{2.95in} 
\caption{Illustration of the existence of $4$ limit cycles 
when $\, a_1 = -\,\frac{30}{7}$, 
$\, a_4 = \frac{1}{3}(a_1-5) =- \frac{65}{21} - \varepsilon_1$, 
and $a_{10}= \frac{1}{2000}$, 
$ b_{11} = \frac{230}{21}\,a_{10} - \varepsilon_2$, 
$b_{01}=-\,a_{10} - \varepsilon_3$, 
where $\, 0 < \varepsilon_3 \ll 
\varepsilon_2 \ll  \varepsilon_1  \ll \varepsilon $: 
(a) An unstable large limit cycle enclosing the center $(1,0)$; 
and (b) Zoomed area around the center $(0,0)$ showing the existence of 
$3$ small limit cycles.}
\label{fig4}
\end{center} 
\vspace{0.20in} 
\end{figure}

\vspace{0.05in}
(A) First, consider the $(3,0)$-distribution. For this case, we have 
$$
\begin{array}{ll}  
a_4 = \dss\frac{1}{3}\, (a_1 - 5), \quad b_{01} = -\,a_{10}, \quad 
b_{11} = - \,10\,(1+a_1)\, a_{10}. 
\end{array}  
$$ 
Taking $\, a_1 = -\, \frac{30}{7}\,$ yields 
$\, a_4 = -\, \frac{65}{21}$, which denotes a point (a blank circle) 
on the line $\, a_4 = \frac{1}{3}\, (a_1 - 5)\,$ in the 
$a_1$-$a_4$ parameter plane (see Fig.~\ref{fig1}). 
Further, we have $\, \ b_{11} = \frac{230}{7}\, a_{10}$, and 
$$
\gamma = \Big(1- \frac{30}{7}\,x \Big)^{-\,\frac{22}{9}} 
\quad (x \ne \frac{7}{30}). 
$$  
Then, the Hamiltonian (\ref{b57}) becomes 
$$ 
H(x,y) =  \dss\frac{16807\, (16250\,y^2 + 13650\,x^2 +2730\,x -441)}
{32500\,(7-30\,x) \, ( 40353607-172944030\,x)^{4/9}} \quad 
{\rm for} \quad x \ne \dss\frac{7}{30},   
$$ 
with 
$$ 
h_{00} = -\,\dss\frac{441}{32500} > 
h_{10} = -\,\dss\frac{33957}{747500} \Big(\dss\frac{23}{7} \Big)^{5/9}. 
$$

The Melnikov functions $\,M_i(h,a_{10})\,$ 
can be expressed as 
\be  
M_i(h,a_{10}) = M_{i0}(h)\, a_{10}, \qquad i=0,\,1.  
\l{b80} 
\ee 
Without loss of generality, we may assume 
\vspace{-0.05in}  
\be  
a_{10} > 0, 
\l{b81} 
\ee  
 
\vspace{-0.05in}  
\noindent 
and thus $\, M_i(h,a_{10}) \,$ and $\, M_{i0}(h)\,$ have the same sign. 
It is noted that for the above chosen parameter values, we have 
$$ 
\mu_{03} = \dss\frac{139150000\, \pi}{453789}\, a_{10} > 0 \quad {\rm and} \quad 
\mu_{10} = - \dss\frac{2500 \sqrt{161}\, \pi}{3703} \, a_{10} < 0. 
$$ 

The computation results of $\, M_{00}(h) \,$ for $\, h \in (h_{00}, \infty) \,$ 
and $\, M_{10}(h) \,$ for $\, h \in (-\infty, h_{10}) \,$ are shown, 
respectively, in Figs.~\ref{fig3}(a) and \ref{fig3}(b).  
Figure~\ref{fig3}(a) shows that $\, M_{00}(h) > 0 \,$ for 
$\, h \in (h_{00}, \infty)$, and its sign agrees with that of  
$\, \mu_{03} > 0\,$ for $\, 0 < h-h_{00} \ll 1$, as expected. 
It is also noted, as shown in Fig.~\ref{fig3}(b), that 
the sign of $\, M_{10}(h)\,$ 
agrees with that of $\, \mu_{10} < 0 \,$ for 
$\, 0 < h_{10} - h \ll 1$. However, unlike the interval 
$\, h \in (h_{00}, \infty)$, this interval contains a critical value 
$\, h = h_1^* \in (-\,0.9250363254,\,-\,0.9250363253) \,$ at which 
$\, M_{10}(h_1^*) = 0 \,$ and the function $\, M_{10}(h) \,$ changes 
its sign as $\,h\,$ crosses this critical point. 
Thus, for this case, besides the $3$ small limit cycles, there exists 
at least one large limit cycle bifurcating from the closed orbit $L_{h_1^*}$ 
of (\ref{b59}). This large limit cycle is shown in Fig.~\ref{fig4}(a), 
which encloses the center $(1,0)$;  
and Fig.~\ref{fig4}(b) illustrates the existence of $3$ small limit 
cycles around the center $(0,0)$.

\begin{figure}[!h]
\vspace{0.05in}
\begin{center}
\hspace{0.00in}

\resizebox{0.70\textwidth}{!}{\includegraphics{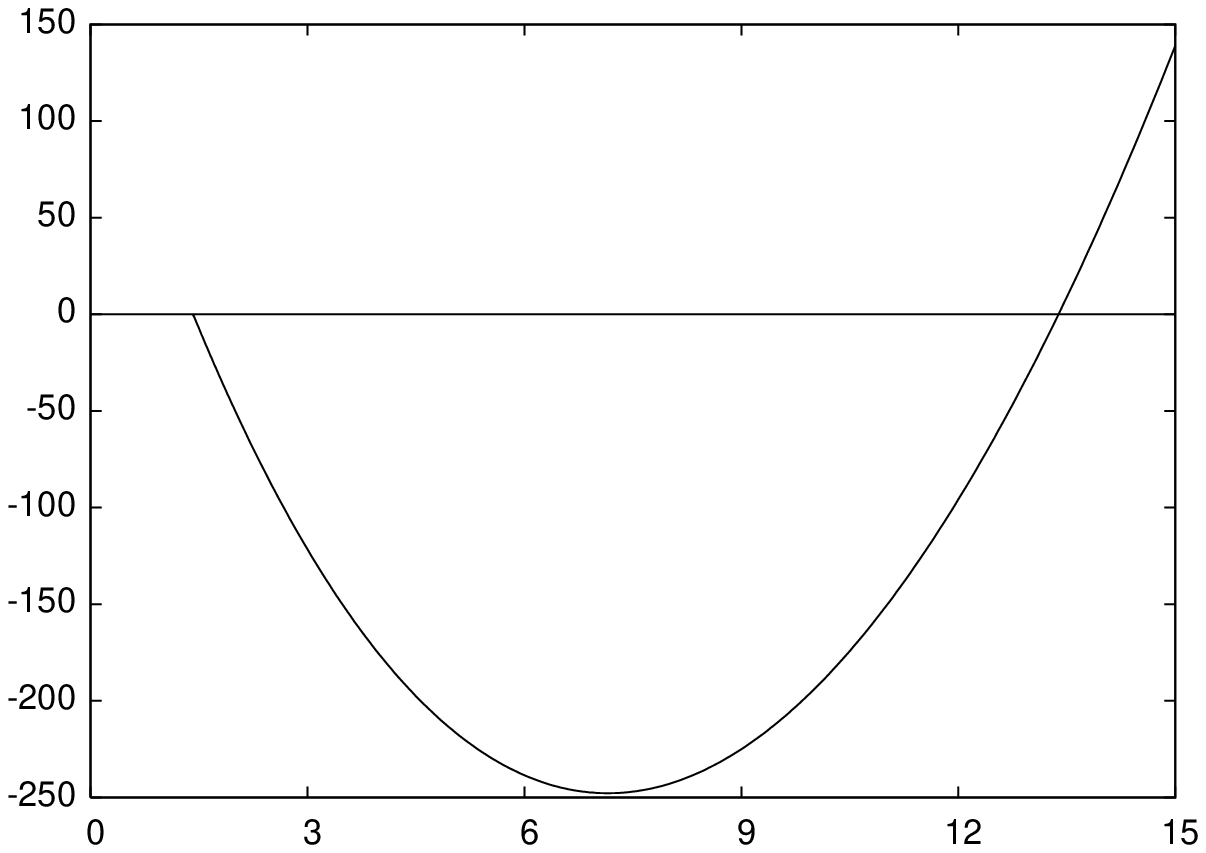}}

\vspace{-0.10in}  
\hspace{0.15in}
$h$ 

\vspace{-2.3in} 
\hspace{-4.6in} 
$M_{00}$  

\vspace{-0.30in} 
\hspace{0.10in} $h_{00}$ \hspace{2.7in} $h_2^*$

\vspace{2.4in} 

\resizebox{0.70\textwidth}{!}{\includegraphics{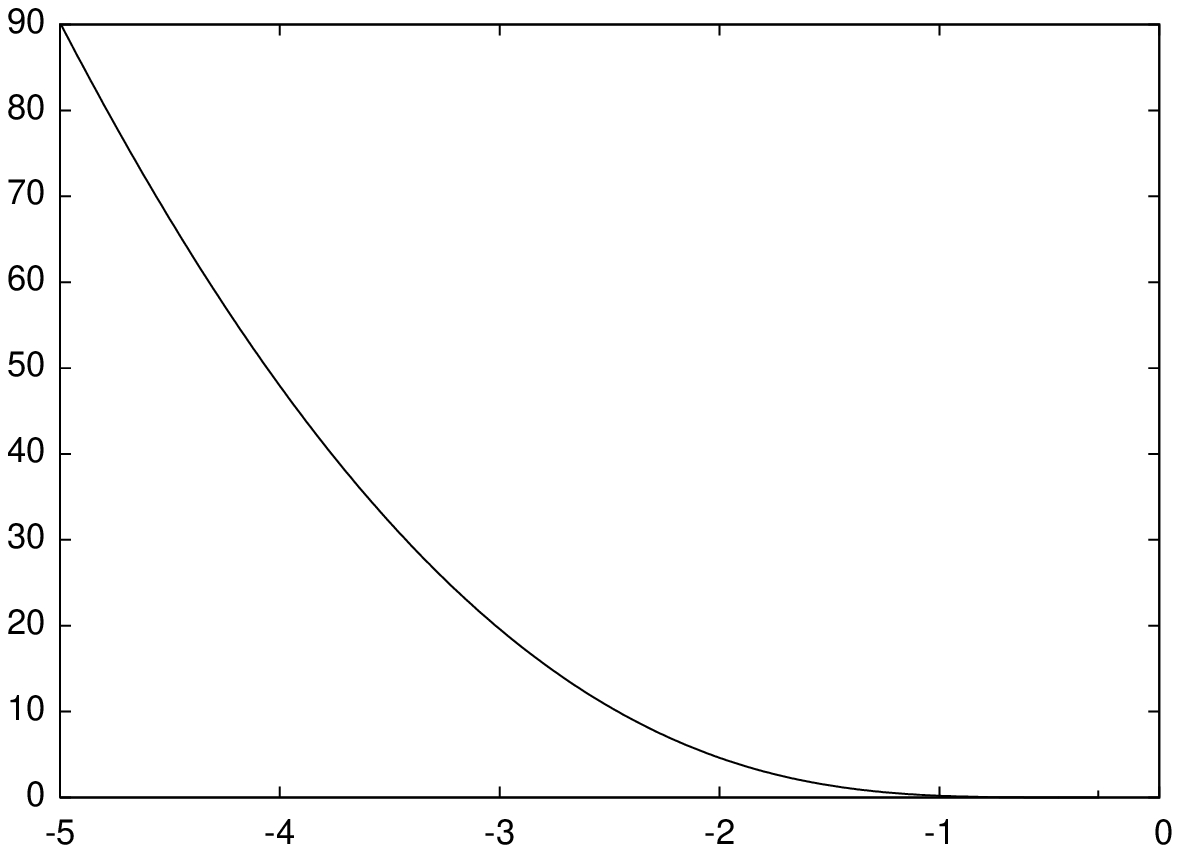}}

\vspace{-2.70in}  
\hspace{ 1.30in} 
\resizebox{0.33\textwidth}{!}{\includegraphics{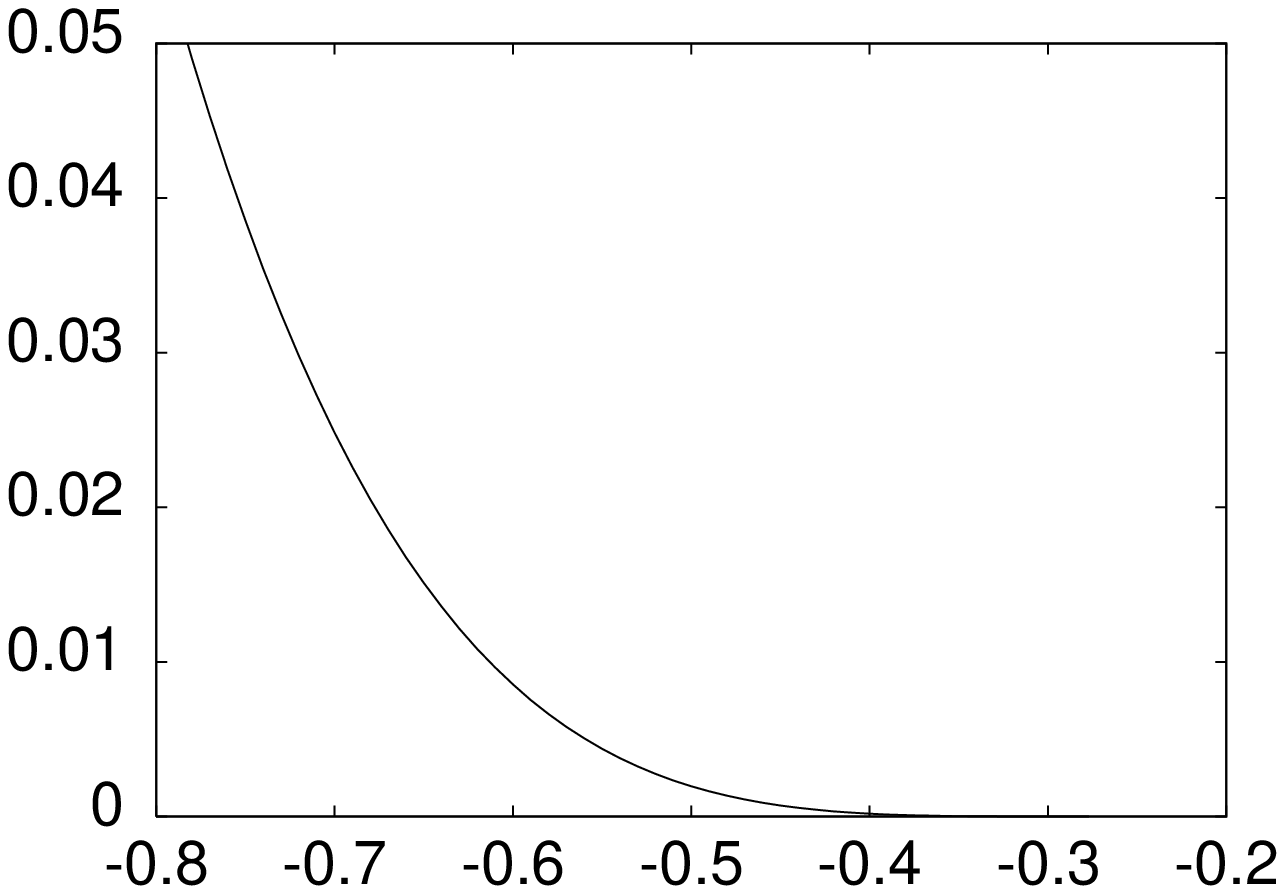}}

\vspace{-0.40in} 
\hspace{2.95in} 
$h_{10}$ 

\vspace{1.25in}  
\hspace{0.08in}  
$h$ 

\vspace{-2.1in} 
\hspace{-4.6in} 
$M_{10}$

\vspace{-4.75in} \hspace{3.20in} (a) 

\vspace{ 3.40in} \hspace{3.20in} (b) 

\vspace{2.8in}

\caption{Functions $M_{00}(h) $ and $M_{10}(h) $ 
under the conditions 
$ \mu_{10} = \mu_{11} = \mu_{12} =0 $, $ \mu_{13} \ne 0$ 
and $\, \mu_{00} \ne 0$,   
for $ a_1 = -\,\frac{70}{51}$ and $ a_4 = \frac{1}{3}(6a_1+5)
=- \frac{55}{51}$: (a) $M_{00}(h) $ for $h \in [h_0,\, +\infty)$, 
with $\, h_{00} =  \frac{7803}{5500} \approx 1.41873$, 
crossing the $h$-axis at $h=h_2^* \in (13.3847179116,\,13.3847179117)$; 
and (b) $\, M_{10}(h) > 0 \,$ 
for $h \in (-\infty,\,h_1]$, with 
$\, h_{10}  = -\,\frac{44217}{104500} ( \frac{19}{51})^{3/7} 
\approx -\,0.27714$.} 
\label{fig5}
\end{center}
\vspace{0.30in} 
\end{figure}

\begin{figure}[!h]
\vspace{0.00in}
\begin{center}
\hspace{0.00in}

\resizebox{0.50\textwidth}{!}{\includegraphics{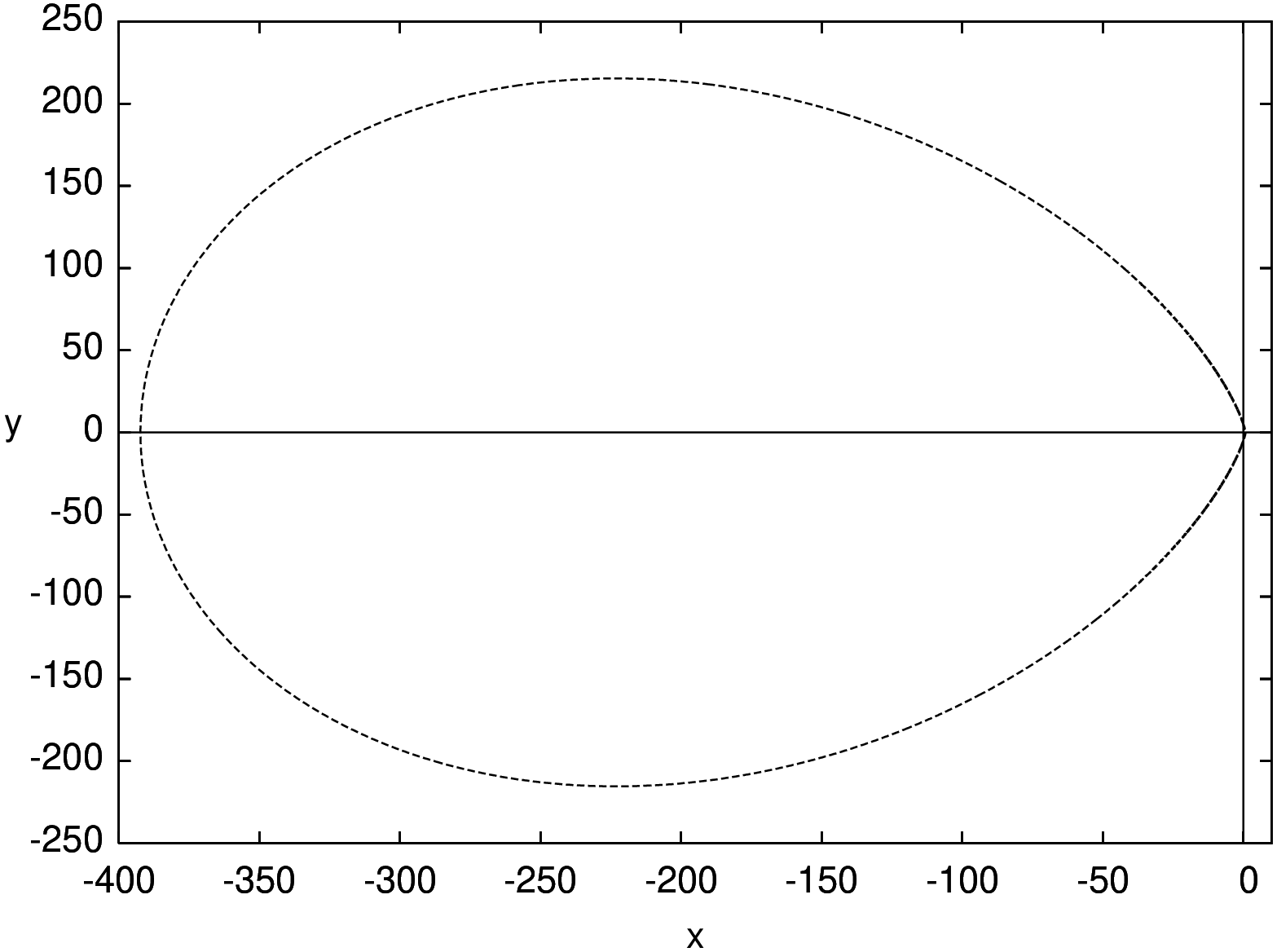}}

\vspace{-1.15in} 
\resizebox{0.50\textwidth}{!}{\includegraphics{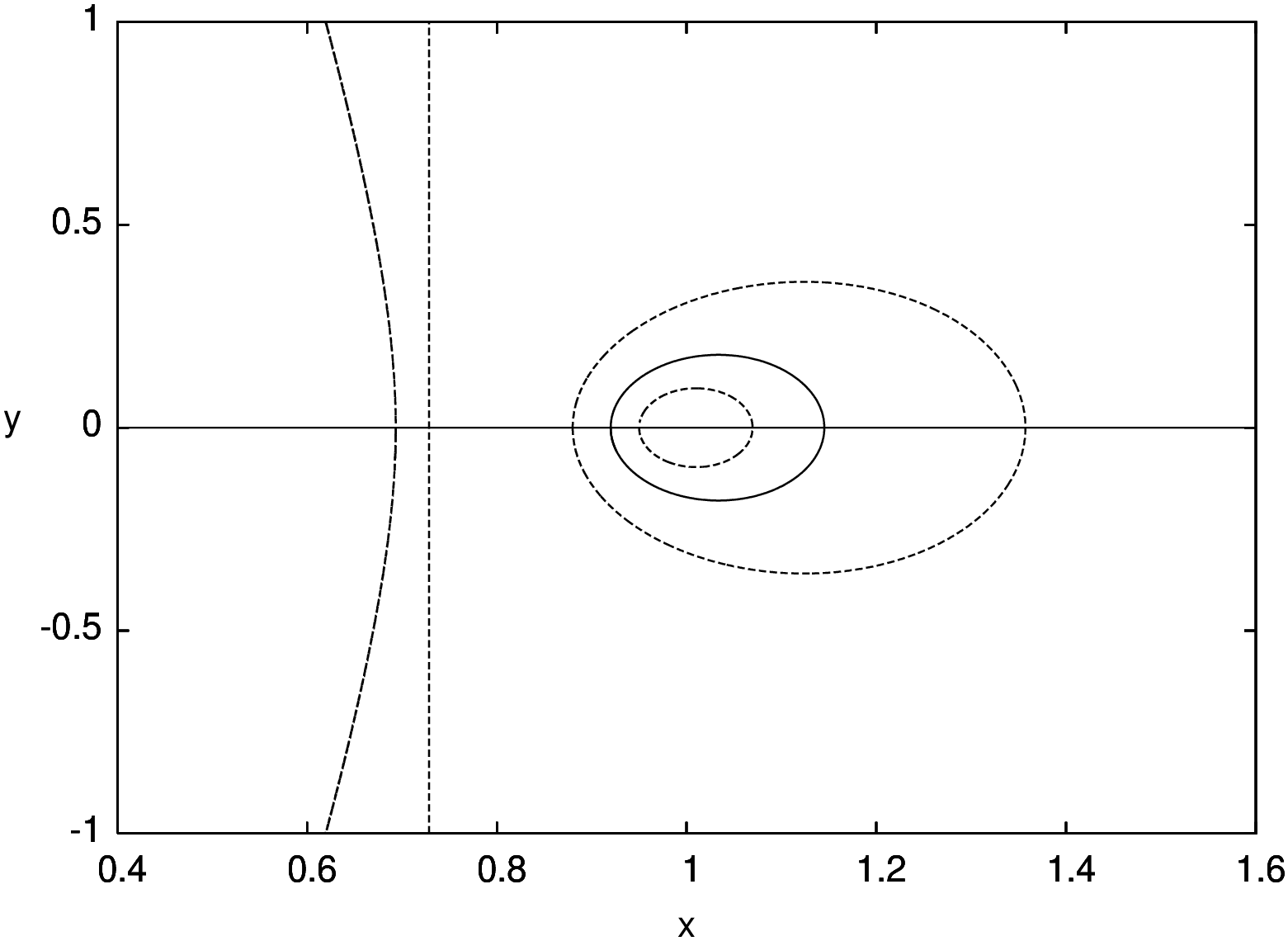}}

\vspace{-1.55in} 
\hspace{-1.43in}$\frac{51}{70}$

\vspace{-6.50in} \hspace{3.70in} (a) 

\vspace{ 3.30in} \hspace{3.70in} (b) 

\vspace{2.95in} 
\caption{Illustration of the existence of $4$ limit cycles 
when $\, a_1 = -\,\frac{70}{51}$, 
$\, a_4 = \frac{1}{3}(6 a_1+5) =- \frac{55}{51} - \varepsilon_1$, 
and $a_{10}= 10$, 
$ b_{11} = \frac{8670}{361}\,a_{10} - \varepsilon_2$, 
$b_{01}=-\,\frac{5611}{361}\,a_{10} + \varepsilon_3$, 
where $\, 0 < \varepsilon_3 \ll 
\varepsilon_2 \ll \varepsilon_1 \ll \varepsilon $:  
(a) An unstable large limit cycle enclosing the center $(0,0)$; 
and (b) Zoomed area around the center  $(1,0)$ showing the 
existence of $3$ small limit cycles.}
\label{fig6}
\end{center} 
\vspace{0.20in} 
\end{figure}

\vspace{0.05in}
(B) For the case of the $(0,3)$-distribution, we have 
$$
\begin{array}{ll}  
a_4 = \dss\frac{1}{3}\, (6\,a_1 + 5), \quad 
b_{01} = -\, b_{11} + \dss\frac{2\,a_4-1}{1+a_1} \,a_{10}, \quad 
b_{11} = \dss\frac{ ( a_1 + 2\, a_4) ( 2\,a_1 - a_4 + 1)}
{(1+a_1)^2 \, (a_1 - a_4 + 1)} \, a_{10}. 
\end{array}  
$$ 
By choosing $\, a_1 = -\, \frac{70}{51}$, we have 
$\, a_4 = -\, \frac{55}{51}, \ b_{01} = -\, \frac{5611}{361}\, a_{10} \,$ and 
$\, b_{11} = \frac{8670}{361}\, a_{10}$. 
The point $(a_1, a_4) = (-\, \frac{70}{51}, -\, \frac{55}{51}) \,$ 
is marked by a blank circle on the line $\, a_4 = \frac{1}{3}\, 
(6\, a_1 + 5) \,$ in the $a_1$-$a_4$ parameter plane 
(see Fig.~\ref{fig1}). Moreover, 
$$
\gamma = \Big(1- \frac{70}{51}\,x \Big)^{-\frac{18}{7}} 
\quad (x \ne \frac{51}{70}), 
$$  
and the Hamiltonian is  
$$ 
H(x,y) =  \dss\frac{345025251\, (2750\,y^2+9350\,x^2-16830\,x+7803)}
{5500\,(51-70\,x) \, ( 897410677851- 1231740146070\,x)^{4/7}} \quad 
{\rm for} \quad x \ne \dss\frac{51}{70},   
$$ 
with 
$$ 
h_{00} = \dss\frac{7803}{5500} > 
h_{10} = -\,\dss\frac{44217}{104500} \Big(\dss\frac{19}{51} \Big)^{3/7}. 
$$ 
For this case, $\mu_{00} \,$ and $\, \mu_{13} \,$ become 
$$ 
\mu_{00} = -\, \dss\frac{10500\, \pi}{361}\, a_{10} < 0 \quad {\rm and} \quad 
\mu_{13} = \dss\frac{4561235000 }{565036352721} \, 
\Big( \dss\frac{51}{19}\Big)^{2/7} \, \pi \, a_{10} > 0. 
$$ 

The computation results of $\, M_{00}(h) \,$ for $\, h \in (h_{00}, \infty) \,$ 
and $\, M_{10}(h) \,$ for $\, h \in (-\infty, h_{10}) \,$ are shown
in Figs.~\ref{fig5}(a) and \ref{fig5}(b), respectively. 
As shown in Fig.~\ref{fig5}(a), the sign of $\, M_{00}(h)\,$ 
agrees with that of $\, \mu_{00} < 0 \,$ for 
$\, 0 < h - h_{00} \ll 1$, and in addition the function $\, M_{00}(h)\,$ 
crosses a critical value at $\, h=h_2^* \in (13.3847179116,\,13.3847179117)$, 
at which it changes sign. 
Figure~\ref{fig5}(b) shows that $\, M_{10}(h) > 0 \,$ for 
$\, h \in (-\infty, h_{10})$, and its sign agrees with that of  
$\, \mu_{13} > 0\,$ for $\, 0 < h_{10}-h \ll 1$. 
Hence, for this case, in addition to the $3$ small limit cycles, 
there also exists 
at least one large limit cycle bifurcating from the  closed orbit $L_{h_2^*}$ of 
(\ref{b59}). This large limit cycle is depicted in Fig.~\ref{fig6}(a), 
which encloses the center $(0,0)$;  
and Fig.~\ref{fig6}(b) illustrates the existence of $3$ small limit 
cycles around the center $(1,0)$.

\begin{figure}[!h]
\vspace{0.00in}
\begin{center}
\hspace{0.00in}

\resizebox{0.70\textwidth}{!}{\includegraphics{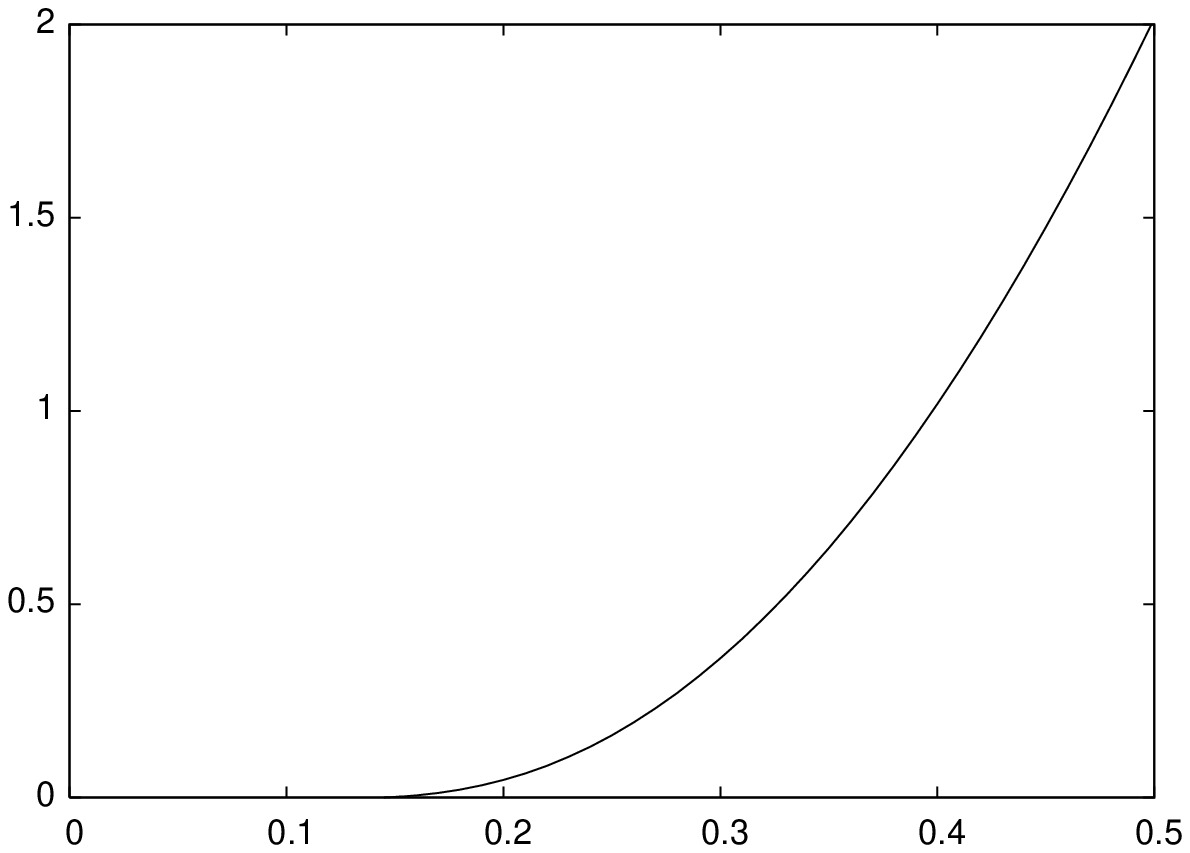}}

\vspace{-2.80in}  
\hspace{-0.70in} 
\resizebox{0.33\textwidth}{!}{\includegraphics{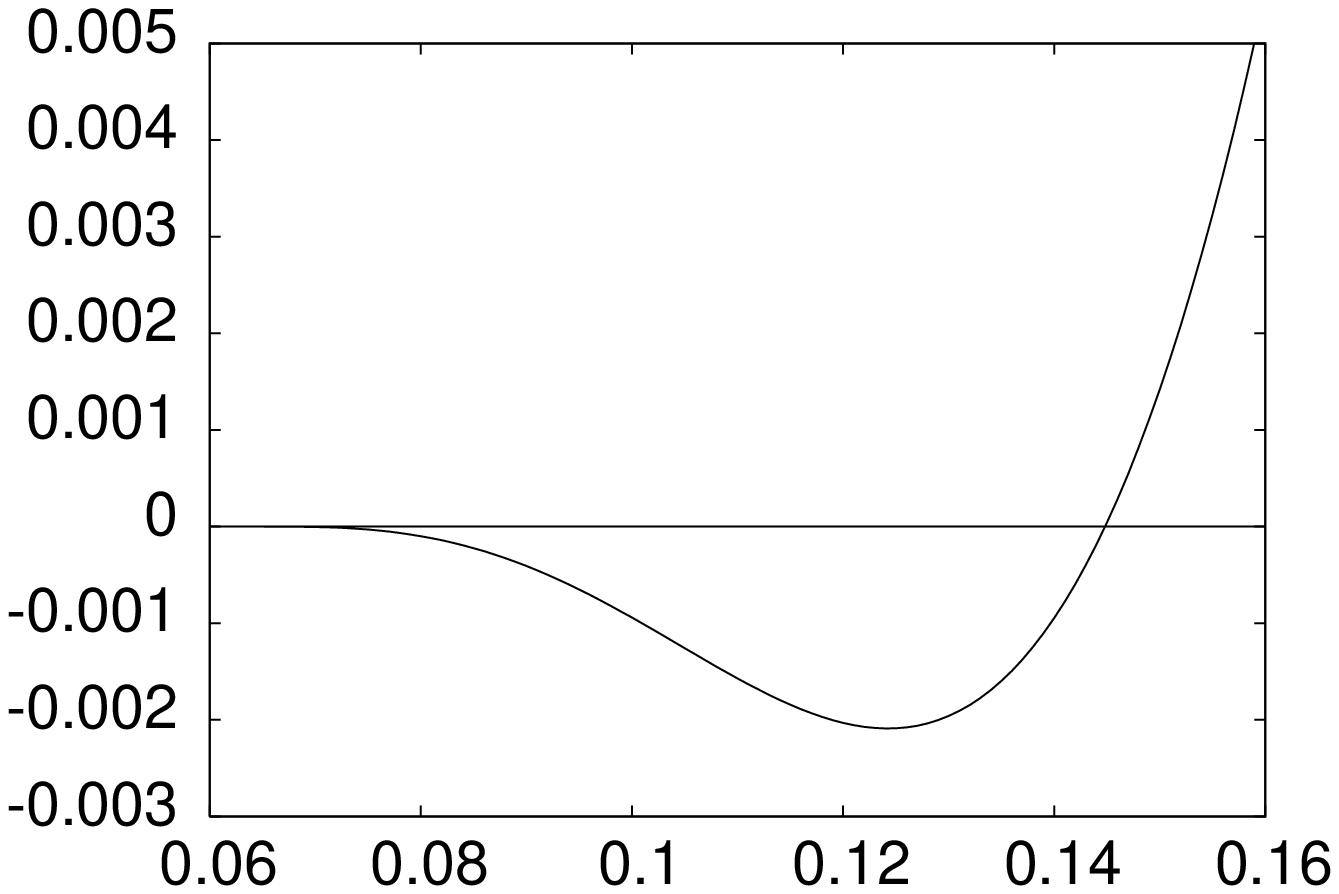}}

\vspace{1.15in}  
\hspace{0.08in}  
$h$ 

\vspace{-1.95in} 
\hspace{-4.5in} 
$M_{00}$  

\vspace{-0.52in} 
\hspace{-0.78in} $h_{00}$ \hspace{1.00in} $h_3^*$

\vspace{2.3in}

\resizebox{0.70\textwidth}{!}{\includegraphics{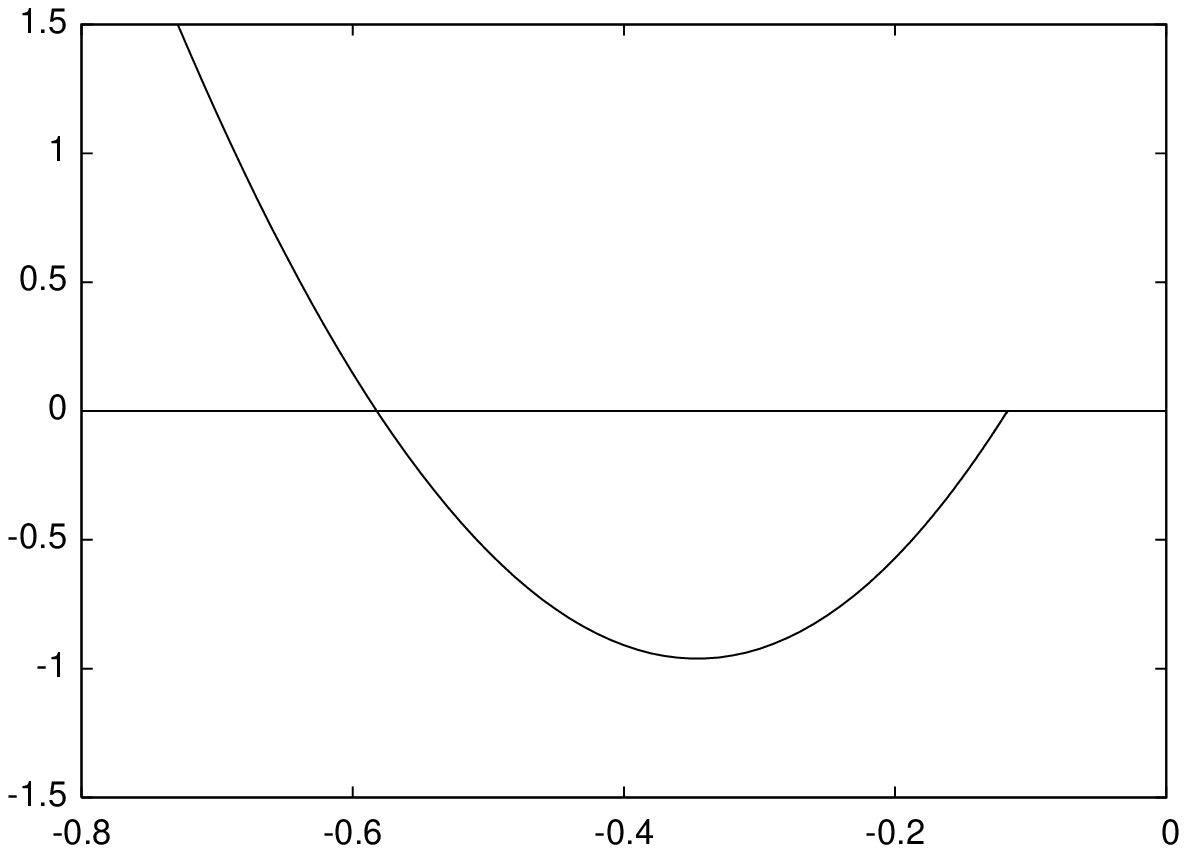}}

\vspace{-0.10in}  
\hspace{0.15in}
$h$ 

\vspace{-2.0in} 
\hspace{-4.6in} 
$M_{10}$  

\vspace{-0.28in} 
\hspace{0.90in} $h_4^*$ \hspace{1.9in} $h_{10}$

\vspace{-4.75in} \hspace{3.20in} (a) 

\vspace{ 3.40in} \hspace{3.20in} (b) 

\vspace{2.8in}

\caption{Functions $M_{00}(h) $ and $M_{10}(h) $ 
under the conditions 
$\, \mu_{00} = \mu_{01} = $, $\, \mu_{02} \ne 0$ 
and $\, \mu_{10} \ne 0$,   
for $\, a_1 = -\,4$ and $\, a_4 = - \frac{18}{5}$: 
(a) $\, M_{00}(h) \,$ 
for $h \in [h_{00}, \, +\infty)$, with 
$\, h_0  = \frac{25}{384} \approx 0.06510$, 
crossing the $h$-axis at $h=h_3^* \in (0.1448192224,\,0.1448192225)$; 
and (b) $M_{10}(h) $ for $h \in (-\infty, \, h_1]$, 
with $\, h_{10} = -\, \frac{325}{3456}3^{1/5} \approx -\,0.11715$, 
crossing the $h$-axis at $h=h_4^*  \in (-\,0.5822537644,\,-\,0.5822537643)$.} 
\label{fig7}
\end{center}
\vspace{0.30in} 
\end{figure}

\begin{figure}[!h]
\vspace{0.00in}
\begin{center}
\hspace{0.00in}

\resizebox{0.50\textwidth}{!}{\includegraphics{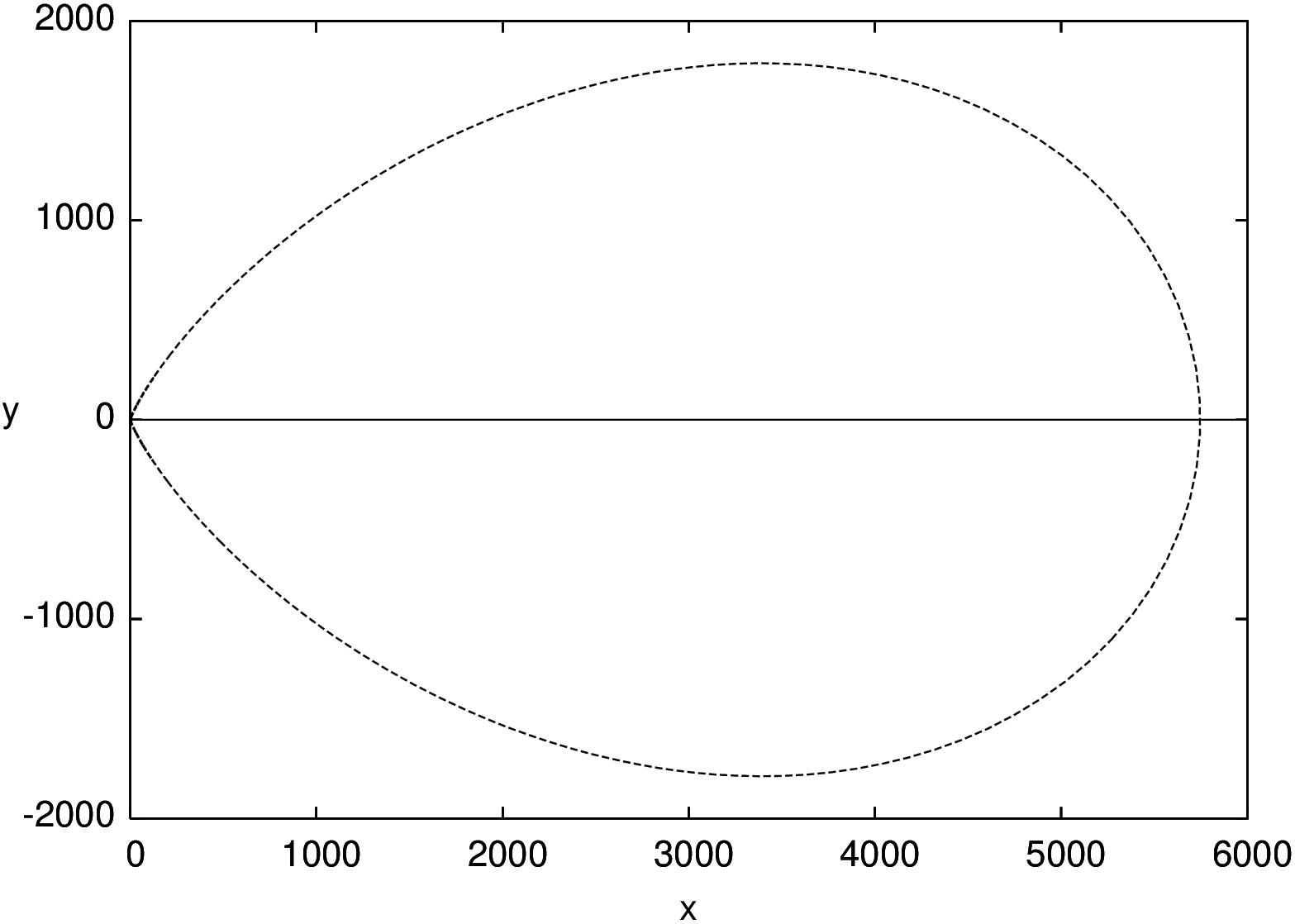}}

\vspace{-1.15in} 
\resizebox{0.50\textwidth}{!}{\includegraphics{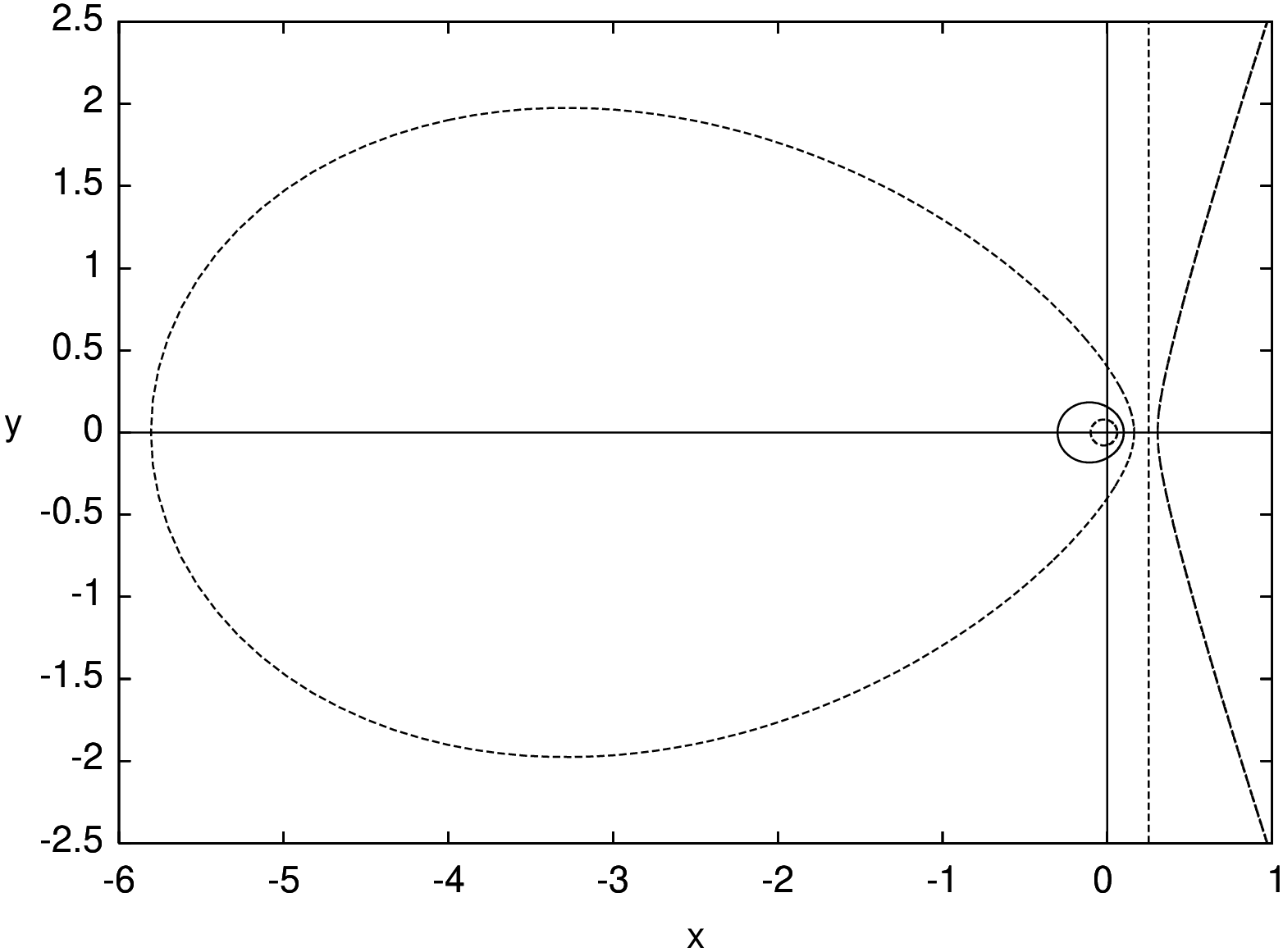}}

\vspace{-1.55in} 
\hspace{ 3.58in}${\small \frac{1}{4}}$

\vspace{-6.45in} \hspace{-3.10in} (a) 

\vspace{ 3.30in} \hspace{-3.10in} (b) 

\vspace{2.95in} 
\caption{Illustration of the existence of $4$ limit cycles 
when  $\, a_1 = -\,4$,  
$\, a_4 = -\,\frac{18}{5}$, 
and $a_{10}= \frac{1}{100}$, 
$\, b_{11} = \frac{392}{65}\, a_{10} -  \varepsilon_1$, and  
$\, b_{01} = -\, a_{10} - \varepsilon_2$, 
where $\, 0 < \varepsilon_2 \ll \varepsilon_1 \ll \varepsilon $:  
(a) An unstable large limit cycle enclosing the center $(1,0)$; 
and (b) Zoomed area around the center  $(0,0)$ showing the 
existence of $1$ large limit cycle and $2$ small limit cycles.}
\label{fig8}
\end{center}
\vspace{0.20in} 
\end{figure}

\vspace{0.05in}
(C) Now consider the $(2,0)$-distribution. For this case, the condition 
$\, a_4 = \frac{1}{3}\, (a_1 - 5) \,$ is not used. We need to 
determine the values for both $\,a_1 \,$ and $\,a_4$. 
We choose 
$$ 
a_1 = -\, 4, \quad a_4 = -\, \dss\frac{18}{5}, 
$$ 
which represents a point in the third quadrant of the 
$a_1$-$a_4$ parameter plane (see the dark circle in 
Fig.~\ref{fig1} near the line $a_4 = \frac{1}{3} (a_1-5)$). Thus,  
$$
\gamma = \Big(1- 4\,x \Big)^{-\frac{14}{5}} 
\quad (x \ne \frac{1}{4}). 
$$  
In addition, we have $\,
b_{01} = -\, a_{10}, \
b_{11} = \frac{392}{65}\, a_{10}$, 
and 
$$ 
H(x,y) =  \dss\frac{192\,y^2+ 480\,x^2-180\,x+25}
{384\,(1-4\,x)^{9/5}} \quad 
{\rm for} \quad x \ne \dss\frac{1}{4},   
$$ 
with 
$$ 
h_{00} = \dss\frac{25}{384} > 
h_{10} = -\,\dss\frac{325}{3456} \ 3^{1/5}. 
$$ 
For this case, $\mu_{02} \,$ and $\, \mu_{10} \,$ are reduced to 
$$ 
\mu_{02} = -\, \dss\frac{1344}{125} \,  \pi\, a_{10} < 0 \quad {\rm and} \quad 
\mu_{10} = -\, \dss\frac{40 \sqrt{3}}{9} \, \pi \, a_{10} < 0. 
$$ 

The computation results of $\, M_{00}(h) \,$ for $\, h \in (h_{00}, \infty) \,$ 
and $\, M_{10}(h) \,$ for $\, h \in (-\infty, h_{10}) \,$ are shown, 
respectively, in Figs.~\ref{fig7}(a) and \ref{fig7}(b).  
As shown in Fig.~\ref{fig7}(a), the sign of $\, M_{00}(h)\,$ 
agrees with that of $\, \mu_{02} < 0 \,$ for 
$\, 0 < h - h_{00} \ll 1$. Moreover, the function $\, M_{00}(h)\,$ 
crosses a critical value at $\, h=h_3^* \in (0.1448192224,\,0.1448192225)\,$  
at which it changes sign. 
Figure~\ref{fig7}(b) shows $\, M_{10}(h) \,$ for 
$\, h \in (-\infty, h_{10})$, whose sign agrees with that of  
$\, \mu_{10} < 0\,$ for $\, 0 < h_{10}-h \ll 1$. Also, 
$\, M_{10}(h) \,$ crosses a critical value at $\, h= h_4^* \in 
(-\,0.5822537644,\,-\,0.5822537643)\,$
at which it changes sign. 
Therefore, for this case, besides the two small limit cycles, 
there exist at least two large limit cycles bifurcating from the two 
different closed orbits $L_{h_3^*}$ and $L_{h_4^*}$ of (\ref{b59}). 
One large limit cycle surrounding the center $(1,0)$ 
is shown in Fig.~\ref{fig8}(a), while another large limit cycle 
enclosing the center $(0,0)$ with $2$ small limit 
cycles is depicted in Fig.~\ref{fig8}(b).

\begin{figure}[!h]
\vspace{0.00in}
\begin{center}
\hspace{0.00in}

\resizebox{0.70\textwidth}{!}{\includegraphics{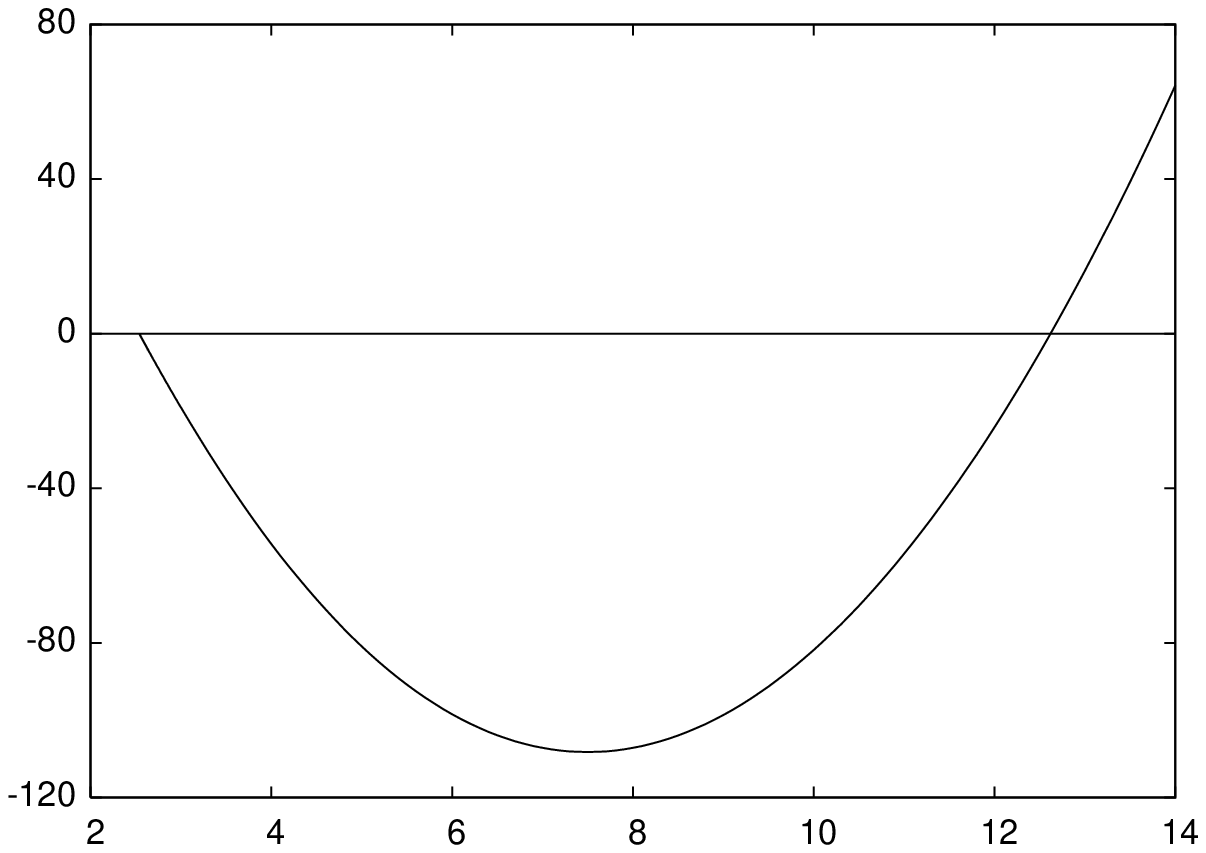}}

\vspace{-0.10in}  
\hspace{0.15in}
$h$ 

\vspace{-2.3in} 
\hspace{-4.6in} 
$M_{00}$  

\vspace{-0.30in} 
\hspace{-0.15in} $h_{00}$ \hspace{2.95in} $h_5^*$

\vspace{2.40in} 

\resizebox{0.70\textwidth}{!}{\includegraphics{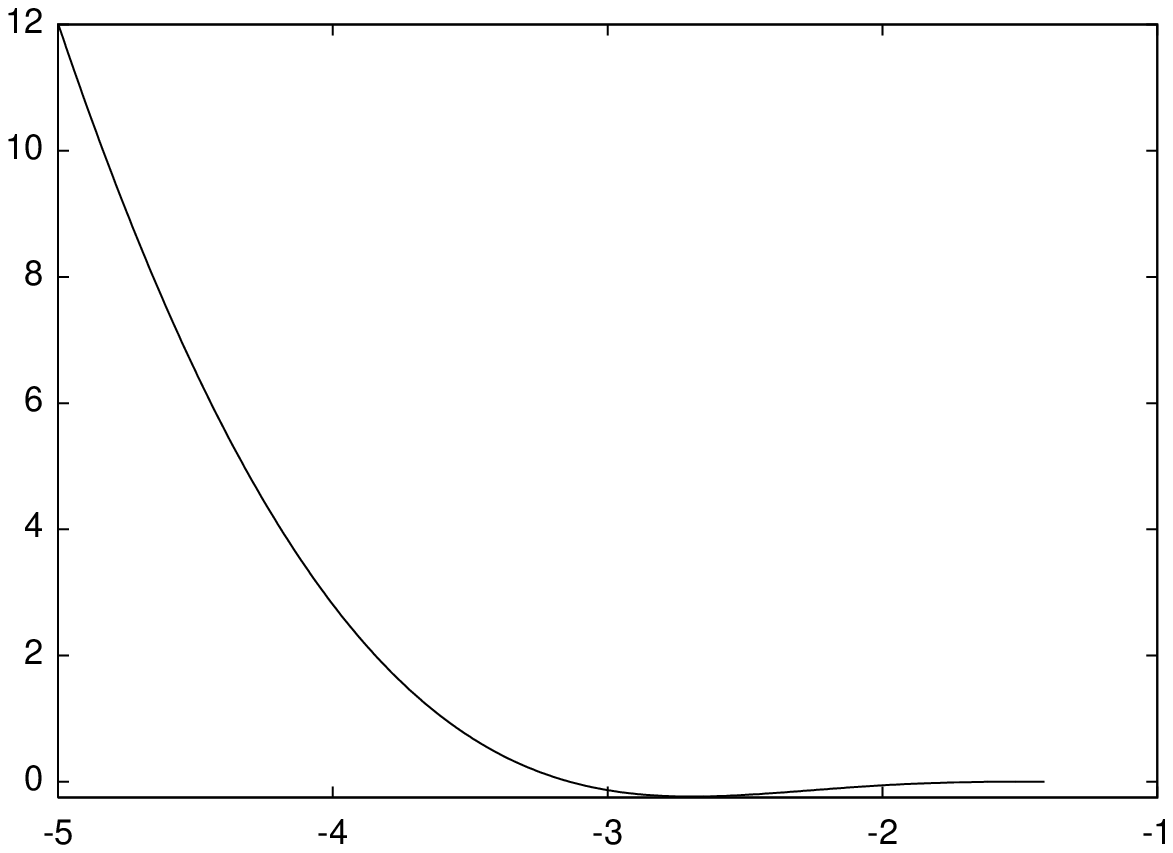}}

\vspace{-2.50in}  
\hspace{1.00in} 
\resizebox{0.33\textwidth}{!}{\includegraphics{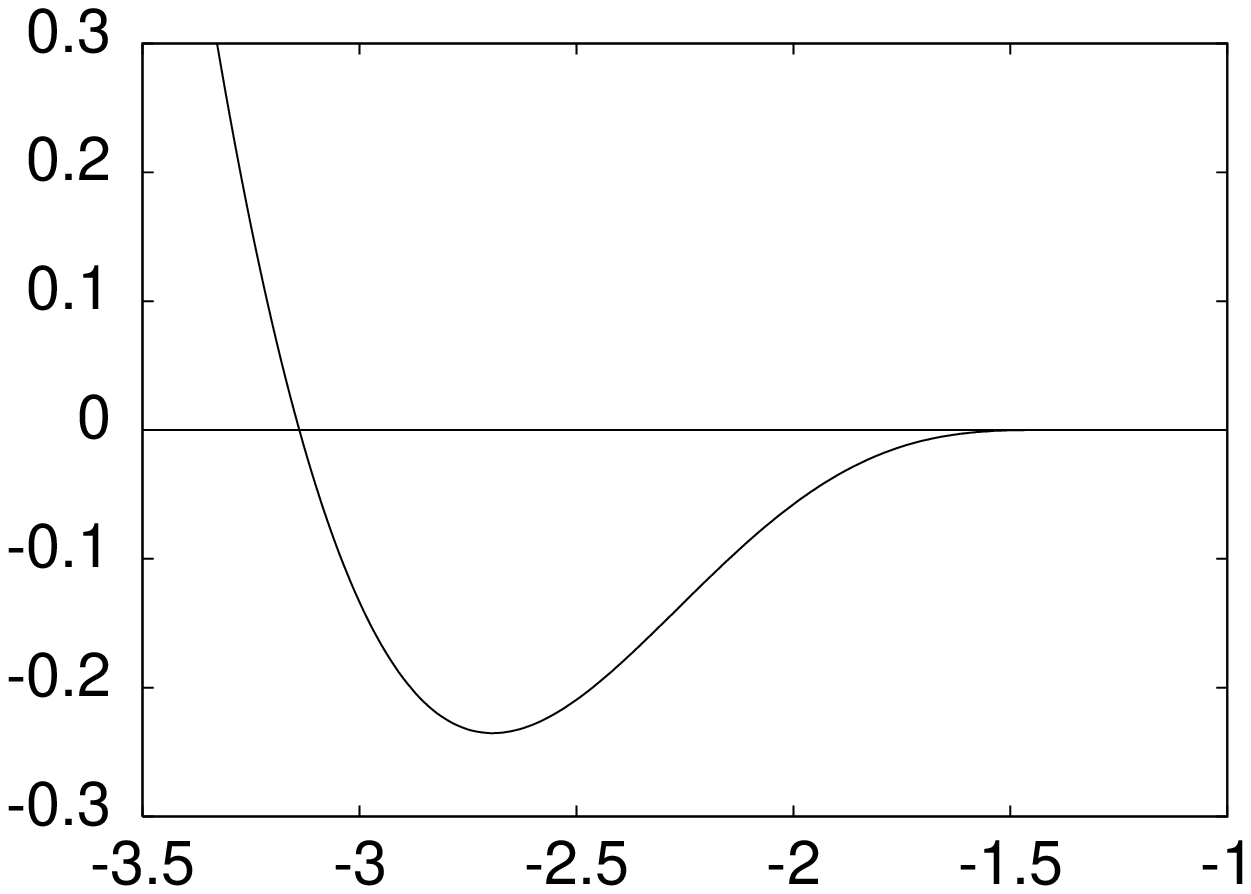}}

\vspace{0.92in}  
\hspace{0.22in}  
$h$ 

\vspace{-2.00in} 
\hspace{-4.6in} 
$M_{10}$  

\vspace{-0.40in} 
\hspace{1.25in} $h_6^*$ \hspace{0.93in} $h_{10}$  

\vspace{-4.65in} \hspace{3.20in} (a) 

\vspace{ 3.35in} \hspace{3.20in} (b) 

\vspace{2.9in} 

\caption{Functions $M_{00}(h) $ and $M_{10}(h) $ 
under the conditions 
$\, \mu_{10} = \mu_{11} = $, $\, \mu_{12} \ne 0$ 
and $\, \mu_{00} \ne 0$,   
for $\, a_1 = -\,\frac{4}{3}$ and $\, a_4 = - \frac{6}{5}$: 
(a) $\, M_{00}(h) \,$ 
for $h \in [h_{00}, \, +\infty)$, with 
$\, h_0  = \frac{325}{128} \approx 2.53096$, 
crossing the $h$-axis at $h=h_5^* \in (12.6197809949,\,12.6197809950)$; 
and (b) $M_{10}(h) $ for $h \in (-\infty, \, h_1]$, 
with $\, h_{10} = -\,\frac{75}{128}\ 3^{4/5} \approx -\,1.41107$, 
crossing the $h$-axis at $h=h_6^*  \in (-\,3.1388150376,\,-\,3.1388150375)$.} 
\label{fig9} 
\end{center}
\vspace{0.10in} 
\end{figure}

\begin{figure}[!h]
\vspace{0.00in}
\begin{center}
\hspace{0.00in}

\resizebox{0.50\textwidth}{!}{\includegraphics{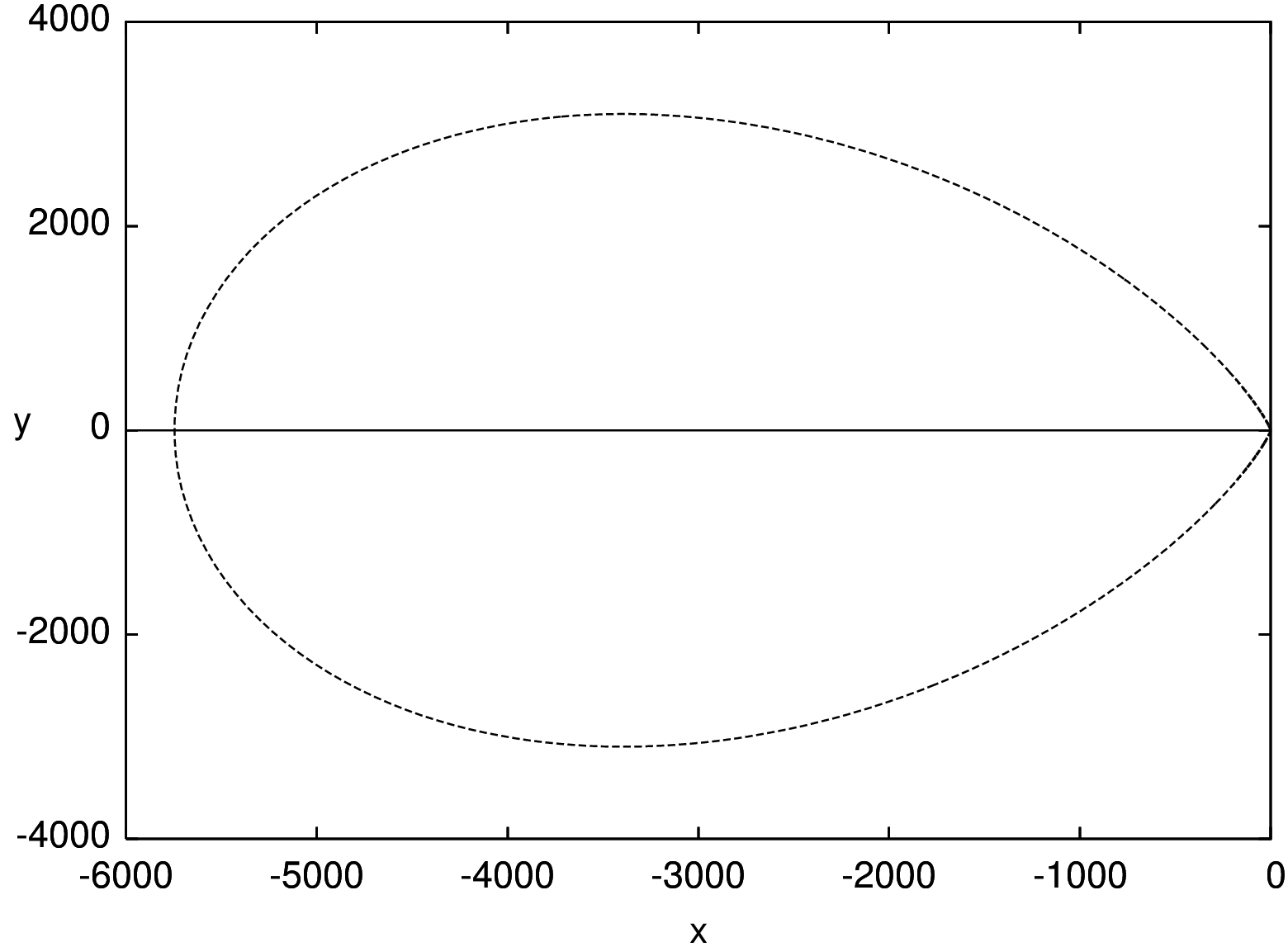}}

\vspace{-1.15in} 
\resizebox{0.50\textwidth}{!}{\includegraphics{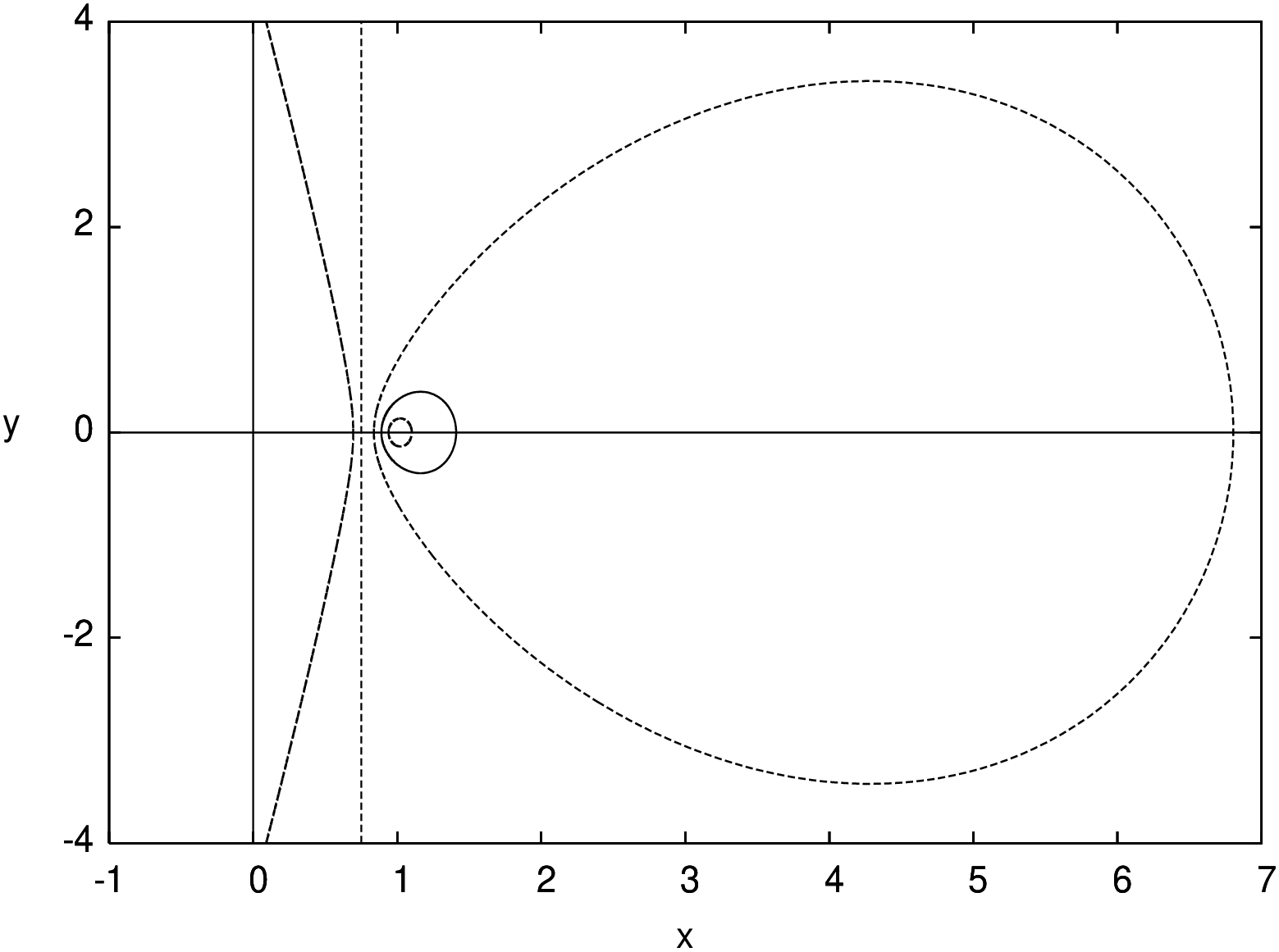}}

\vspace{-1.55in} 
\hspace{-1.88in}${\small \frac{3}{4}}$

\vspace{-6.45in} \hspace{3.70in} (a) 

\vspace{ 3.30in} \hspace{3.70in} (b) 

\vspace{2.95in} 
\caption{Illustration of the existence of $4$ limit cycles 
when $\, a_1 = -\,\frac{4}{3}$, $\, a_4 = - \frac{6}{5}$, 
and $a_{10}=1$, 
$ b_{11} = \frac{1176}{65}\, a_{10} -  \varepsilon_1$, 
$b_{01}=-\,\frac{513}{65} \, a_{10} +  \varepsilon_2$, 
where $\, 0 < \varepsilon_2 \ll \varepsilon_1 \ll \varepsilon $:  
(a) An unstable large limit cycle enclosing the center $(0,0)$; 
and (b) Zoomed area around the center  $(1,0)$ showing the 
existence of $1$ large limit cycle and $2$ small limit cycles.}
\label{fig10}
\end{center}
\vspace{0.00in} 
\end{figure}

\vspace{0.05in}
(D) Finally, consider the $(0,2)$-distribution. For this case, the condition 
$\, a_4 = \frac{1}{3}\, (6\,a_1 + 5) \,$ is not used. Taking 
$$ 
a_1 = -\, \dss\frac{4}{3}, \quad a_4 = -\, \dss\frac{6}{5}, 
$$ 
yields 
$$
\gamma = \Big(1- \dss\frac{4}{3}\,x \Big)^{-\frac{14}{5}} 
\quad (x \ne \frac{3}{4}). 
$$  
The point $(-\frac{4}{3}, -\frac{6}{5})$ is marked by a dark circle 
near the line $\, a_4 = \frac{1}{3}\, ( 6\, a_1 + 5)\,$ 
in the $a_1$-$a_4$ parameter plane (see Fig.~\ref{fig1}). 
Further, we have $\,
b_{01} = -\, \frac{513}{65} \, a_{10}, \
b_{11} = \frac{1176}{65}\, a_{10}$, 
and 
$$ 
H(x,y) =  \dss\frac{243\,(64\,y^2+ 480\,x^2-780\,x+325}
{((3-4\,x)\, (324\,x - 243)^{4/5}} \quad 
{\rm for} \quad x \ne \dss\frac{3}{4},   
$$ 
with 
$$ 
h_{00} = \dss\frac{325}{128} > 
h_{10} = -\,\dss\frac{75}{128} \ 3^{4/5}. 
$$ 
For this case, $\mu_{00} \,$ and $\, \mu_{12} \,$ are simplified as 
$$ 
\mu_{00} = -\, \dss\frac{896}{65} \,  \pi\, a_{10} < 0 \quad {\rm and} \quad 
\mu_{12} = -\, \dss\frac{448}{30375} \, 3^{9/10}\, \pi \, a_{10} < 0. 
$$ 

The computation results of $\, M_{00}(h) \,$ for $\, h \in (h_{00}, \infty) \,$ 
and $\, M_{10}(h) \,$ for $\, h \in (-\infty, h_{10}) \,$ are shown
in Figs.~\ref{fig9}(a) and \ref{fig9}(b), respectively.  
As shown in Fig.~\ref{fig9}(a), the sign of $\, M_{00}(h)\,$ 
agrees with that of $\, \mu_{00} < 0 \,$ for 
$\, 0 < h - h_{00} \ll 1$, and the function $\, M_{00}(h)\,$ 
crosses a critical value at $\, h=h_5^* \in (12.6197809949,\,12.6197809950)\,$ 
at which it changes sign. 
Figure~\ref{fig9}(b) shows $\, M_{10}(h) \,$ for 
$\, h \in (-\infty, h_{10})$, whose sign agrees with that of  
$\, \mu_{12} < 0\,$ for $\, 0 < h_{10}-h \ll 1$. Moreover, 
$\, M_{10}(h) \,$ crosses a critical value at $\, h=h_6^* \in 
(-\,3.1388150376,\,-\,3.1388150375)\,$ 
at which it changes sign. 
Therefore, for this case, in addition to the two small limit cycles, 
there also exist at least two large limit cycles bifurcating from the two 
different closed orbits $L_{h_5^*}$ and $L_{h_6^*}$ of (\ref{b59}). 
One large limit cycle surrounding the center $(0,0)$ 
is shown in Fig.~\ref{fig10}(a), while another large limit cycle 
enclosing the center $(1,0)$ with $2$ small limit 
cycles is depicted in Fig.~\ref{fig10}(b). 

It is noted that all the four sets of values of $a_1$ and $a_4$ chosen above 
in (A), (B), (C) and (D) satisfy 
\be  
\dss\frac{a_1 + 2\, a_4}{a_1} = \dss\frac{2\,n}{m}, \quad 
{\rm where} \ \, n \ \, {\rm is \ an \ integer \ and} \ \, m \ \, 
{\rm is \ an \ odd \ integer}, 
\l{b82} 
\ee  
so that a consistent integrating factor (and so a consistent 
Hamiltonian function for the whole transformed system) is obtained. 
However, this condition is not necessary 
since the singular line $\, 1+a_1\,x = 0\,$ 
divides the phase plane into two parts, and the analysis 
does not need the continuity on the singular line. 
To demonstrate this, in the following 
we present a case for which the condition
(\ref{b82}) is not satisfied. Consider the $(2,0)$-distribution, and choose 
$\, a_1 = -\,5\,$ and $\, a_4 = -\,4$. The point $(a_1, a_4) = (-5,\,-4)\,$ 
is marked by a dark circle in the $a_1$-$a_4$ parameter plane 
(see Fig.~\ref{fig1}). Then, 
$$ 
\dss\frac{a_1+2 a_4}{a_1} = \dss\frac{13}{5}, \quad 
b_{01}= a_{10}, \quad b_{11} = \dss\frac{26}{3}\, a_{10},  
$$ 
 
\vspace{-0.05in}  
\noindent 
and \vspace{0.10in}  
$$
H(x,y) = \left\{ 
\begin{array}{lll}  
\dss\frac{ x^2 + y^2}{2 \,( 1 - 5\,x)^{8/5}}, \ \ & 
\forall \ h \in (0, \infty), \ \ & {\rm when} \ \ x < \dss\frac{1}{5}, \\[2.0ex]  
-\, \dss\frac{ x^2 + y^2}{2 \,( 1 - 5\,x)^{8/5}},  & 
\forall \ h \in (-\infty, -\frac{1}{32}\, 2^{4/5}), & {\rm when} \ \ 
x > \dss\frac{1}{5}. 
\end{array}
\right.
$$ 
For this case, $\, \mu_{02} \,$ and $\, \mu_{10} \,$ become
$$ 
\mu_{02} = -\, \dss\frac{130}{3}\, \pi \, a_{10} < 0 \quad {\rm and} \quad 
\mu_{10} = -\, \dss\frac{65}{12}\, \pi \, a_{10} < 0. 
$$ 

The computation result of $\, M_{00} (h) \,$ shows that 
$\, M_{00}(h) < 0 \,$ for $\, 0 < h \ll 1$, agrees with the sign of 
$\, \mu_{02}$. Moreover, $\, M_{00}(0.1) = 0.0510077880 > 0$, implying that 
there exists $ h = h_7^* \in (0,\, 0.1) \,$ such that $\, M_{00}(h_7^*)=0$, 
and so a large limit cycle bifurcates from the closed orbit 
$\, L_{h_7^*}$ of (\ref{b59}). 
The result of $ \,M_{10}(h) \,$ also shows that $\, M_{10}(h) < 0 \,$ 
for $\, 0 < -\,\frac{1}{32}\, 2^{4/5} - h \ll 1 $, agreeing with the sign of 
$\, \mu_{10}$, and that $\, M_{10}(-\,\frac{1}{32}\, 2^{4/5}-0.8) 
= 7.4630743072 > 0$, implying the existence 
$\, h = h_8^* \in (-\,\frac{1}{32}\, 2^{4/5}-0.8,\, -\,\frac{1}{32}\, 2^{4/5}) 
= (-0.8544094102, \, -0.0544094102)\,$ such that 
$\, M_{10}(h_8^*)=0$. Thus, there exists another large limit cycle 
bifurcating from the closed orbit 
$\, L_{h_8^*}$ of (\ref{b59}). 
Therefore, this case exhibits $2$ small limit 
cycles and $2$ large limit cycles, leading to the existence of 
at least $4$ limit cycles.

Summarizing the above results with the continuity of parameters 
$\,a_1\,$ and $\, a_4\,$ shows that 
at least for some regions in the $a_1$-$a_4$ parameter plane
the reversible near-integrable system (\ref{b55}) can exhibit 
at least $4$ limit cycles around the two singular points $(0,0)$ and $(1,0)$ 
with distribution ether $(3,1)$ or $(1,3)$. 

The proof of Theorem 4.1 is finished. 
\put(10,0.5){\framebox(6,7.5)}

\section{Conclusion} 
In this paper, we have proved that a quadratic non-Hamiltonian integrable 
system with two centers can have at least $4$ limit cycles 
under quadratic perturbations, with distributions either 
$(3,1)$ or $(1,3)$.  
This result gives a new record, answering the open problem of the existence of 
limit cycles in near-integrable quadratic systems. 
It is shown that such systems 
can have at least $4$ limit cycles for some regions 
in the $2$-dimensional parameter plane, associated with the parameters 
of the integrable systems.
Further research is needed on global analysis for all possible parameter 
values in the parameter plane. 

\section*{Acknowledgments}
This work was supported by 
the Natural Sciences and Engineering Research Council of Canada (NSERC) 
and the National Natural Science Foundation of
China (NNSFC).


\end{document}